\documentclass{article}
\usepackage{epsfig}
\usepackage{graphicx}
\usepackage{color}
\usepackage{times}
\usepackage{amsmath}
\usepackage{amssymb}
\usepackage{mathtools}
\usepackage{xspace}
\usepackage{pb-diagram}
\usepackage{algorithm}
\usepackage{algorithmic}

\newcommand {\mm}[1] {\ifmmode{#1}\else{\mbox{\(#1\)}}\fi}

\newcommand{\eop}{\hfill\usebox{\smallProofsym}\bigskip}  %
\newsavebox{\smallProofsym}                            
\savebox{\smallProofsym}                               %
{
\begin{picture}(6,6)                                   %
\put(0,0){\framebox(6,6){}}                            %
\put(0,2){\framebox(4,4){}}                            %
\end{picture}                                          %
}                                                      %
\makeatletter
\long\def\@makecaption#1#2{%
  \vskip\abovecaptionskip
  \sbox\@tempboxa{\small #1: #2}%
  \ifdim \wd\@tempboxa >\hsize
    \small #1: #2\par
  \else
    \global \@minipagefalse
    \hb@xt@\hsize{\hfil\box\@tempboxa\hfil}%
  \fi
  \vskip\belowcaptionskip}
\makeatother

\newtheorem{theorem}{Theorem}[section]
\newtheorem{lemma}{Lemma}[section]

\newtheorem{definition}{Definition}[section]

\newcommand{\Rspace}        {\mm{{\mathbb R}}}
\newcommand{\Xspace}        {\mm{{\mathbb X}}}
\newcommand{\Yspace}        {\mm{{\mathbb Y}}}
\newcommand{\Zspace}        {\mm{{\mathbb Z}}}

\newcommand{\Sspace}        {\mm{{\mathbb S}}}

\newcommand{\Uspace}        {\mm{{U}}}

\newcommand{\Hgroup}        {\mm{\sf H}}
\newcommand{\Hgr}             {\mm{\sf H}}

\newcommand{\Cgroup}        {\mm{\sf C}}

\newcommand{\bdr}        {{\partial}}

\newcommand{\Dcal}        {\mm{\mathcal D}}
\newcommand{\Mcal}        {\mm{\mathcal M}}

\newcommand{\Fcal}        {\mm{\mathcal F}}
\newcommand{\Gcal}        {\mm{\mathcal G}}
\newcommand{\Acal}        {\mm{\mathcal A}}
\newcommand{\Ccal}        {\mm{\mathcal C}}

\newcommand{\Ddgm}[2]       {\mm{\rm Dgm}_{#1}{({#2})}}

\newcommand{\kernel}[1]     {\mm{\rm ker\,}{#1}}
\newcommand{\cokernel}[1]   {\mm{\rm cok\,}{#1}}
\newcommand{\image}[1]      {\mm{\rm im\,}{#1}}

\newcommand{\ep}            {\mm{\varepsilon}}

\newcommand{\vol}[1]       {\mm{\rm vol\,}({#1})}

\newcommand{\prb} {\xi}
\newcommand{\dgmD}[1]       {\mm{\rm Dgm}{({#1})}}
\newcommand{\cok}[1]   {\mm{\rm cok\,}{#1}}
\newcommand{\voi}[1]       {\mm{\rm Voi\,}{({#1})}}
\newcommand{\del}[1]       {\mm{\rm Del\,}{({#1})}}

\newcommand{\sd}[1]   {\mm{\rm Sd\,}{#1}}
\newcommand{\Nerve}[1]       {\mm{\rm Nerve\,}{#1}}
\newcommand{\A}        {\mm{{\mathcal A}}}
\newcommand{\Lune}        {\mm{{\mathsf L}}}
\newcommand{\Moon}       {\mm{{\mathsf M}}}
\newcommand{\Bisector}       {\mm{{\mathsf P}}}
\newcommand{\C}        {\mm{{\mathcal C}}}
\newcommand{\V}        {\mm{{\mathcal V}}}
\newcommand{\T}        {\mm{{\mathcal T}}}
\newcommand{\interior}[1]   {\mm{\rm int\,}{#1}}
\newcommand{\closure}[1]    {\mm{\rm cl\,}{#1}}
\newcommand{\distone}[1]   {\mm{\rm dist\,}{#1}}

\title{Towards Stratification Learning through Homology Inference}
\author{Paul Bendich\thanks{IST Austria, Klosterneuburg, Austria}, Bei Wang\thanks{Duke University, Durham, NC}, and Sayan Mukherjee\thanks{Duke University, Durham, NC}
}

\begin{document}
\maketitle
\abstract{A topological approach to stratification learning is developed
for point cloud data drawn from a stratified space. Given such data, our objective
is to infer which points belong to the same strata.
First we define a multi-scale notion of a stratified space, giving a stratification for each radius level.
We then use methods derived from kernel and cokernel persistent homology to cluster the data points
into different strata, and we prove a result which guarantees the correctness of our clustering, given
certain topological conditions; some geometric intuition for these topological conditions is also provided.
Our correctness result is then given a probabilistic flavor: we give bounds on the minimum number of sample points
required to infer, with probability, which points belong to the same strata.
Finally, we give an explicit algorithm for the clustering, prove its correctness, and apply it to some simulated data.}

\section{Introduction}
\label{sec:intro}

Manifold learning is a basic problem in geometry, topology, and statistical inference that has received a great deal of attention.
The basic idea is as follows: given a point cloud of data sampled from a manifold  in an ambient
space $\Rspace^k$, infer the underlying manifold. 
A  limitation of the problem statement is that it does not apply to sets that are not manifolds.
For example, we may consider the more general class of stratified spaces that can be decomposed into 
strata, which are manifolds of varying dimension, each of which fit together in some uniform way inside 
the higher dimensional space. 

In this paper, we study the following problem in stratification learning: given a point cloud sampled
from a stratified space, how do we cluster the points so that points in the same cluster are in the same stratum, while
points in different clusters are not? 
Intuitively, the strategy should be clear: two points belong in the same stratum if they ``look the same locally,'' meaning
that they have identical neighborhoods, within the larger space, at some very small scale.
However, the notion of ``local'' becomes unclear in the context of sampling uncertainty, since everything becomes quite noisy
at vanishingly small scale.
In response, we introduce a radius parameter $r$ and define a notion of local equivalence at each such $r$.

Our tools are derived from algebraic topology. In particular, we define local equivalence between points via maps between relative homology groups,
and we then attempt to infer this relation by using ideas coming from persistent homology \cite{EdeHar2010}.

\paragraph{Prior Work}

Consistency in manifold learning has often been recast as a homology inference statement:
as the number of points in a point cloud goes to infinity, the inferred homology
converges to the true homology of the underlying space.
Results of this nature have been given for manifolds  \cite{NiySmaWei2008,NiySmaWei2008b}
and a large class of compact subsets of Euclidean space \cite{ChaCohLie2009}. 
Stronger results in homology inference for closed subsets of a metric space are given in \cite{CohEdeHar2007}.

Geometric approaches to stratification inference have been developed.
These include inference
of a mixture of linear subspaces  \cite{LerZha2010}, mixture models for general stratified
spaces \cite{HarRanSap2005}, and generalized Principal Component Analysis (GPCA)
\cite{VidMaSas2005} which was developed for dimension reduction for mixtures of manifolds.

The study of stratified spaces has long
been a focus of pure mathematics; see, for example, \cite{GorMac1988,Wei1994}. 
The problem of inference for the
local homology groups of a sampled stratified space in a deterministic setting has been addressed in
\cite{BenCohEde2007}.

\paragraph{Results}
In this paper we propose an approach to stratification inference based on local homology inference; more specifically,
based on inference of the kernels and cokernels of several maps between groups closely related to the multi-scale local homology groups
for different pairs of points in the sample.
The results in this paper are: (1) a topological definition of two points belonging to the same strata by
assessing the multi-scale local structure of the points through kernel and cokernel persistent homology; (2) topological 
conditions on the point sample under which the topological characterization holds -- we call this topological 
inference; (3) a geometric intuition of these topological conditions based on quantities related to 
reach and to the gradient of a distance function; (4) finite sample bounds for the minimum number of points in the sample
required to state with high probability which points belong to the same strata; (5)
an algorithm that computes which points belong to the same strata and a proof of
correctness for some parts of this algorithm.

\paragraph{Outline}
We review the needed background in Section \ref{sec:Back}.
In Section \ref{sec:TIT}, we give the topological inference theorem and provide
some geometric intuition. The probabilistic statements are provided in Section \ref{sec:PIT}.
We describe the clustering algorithm in Section \ref{sec:Alg}; the correctness proof
of the algorithm is contained in three Appendices, \ref{app:phimap} through \ref{app:correctness}.
The main body of the paper closes with some further discussion in Section \ref{sec:Disc}.


 \section{Background}
\label{sec:Back}

We review necessary background on persistent homology and stratified spaces. 
The treatment of the former here is mostly adapted from \cite{ChaCohGli2009}, although we present the material in slightly simplified form.
We first describe general persistence modules, focusing mainly on those that arise from maps between homology groups induced by inclusions of topological spaces. 
We then discuss stratifications and their connection to the local homology groups of
a topological space. 
Basics on homology itself are assumed; for a readable background, see
\cite{Mun1984} or \cite{Hat2002}, or \cite{EdeHar2010} for a more computationally oriented treatment.

\subsection{Persistence Modules}

In \cite{ChaCohGli2009}, the authors define persistence modules over an arbitrary commutative ring $R$ with unity.
For simplicity, we restrict immediately to the case $R = \Zspace / 2 \Zspace$.
Let $A$ be some subset of $\Rspace$.
Then a \emph{persistence module} $\Fcal_A$ is a collection $\{F_{\alpha}\}_{\alpha \in A}$ of $
\mathbb{Z} / 2 \mathbb{Z}$-vector spaces,
together with a family $\{f_{\alpha}^{\beta}: F_{\alpha} \to F_{\beta}\}_{\alpha \leq \beta \in A}$ of
linear maps
such that $\alpha \leq \beta \leq \gamma$ implies $f_{\alpha}^{\gamma} = f_{\beta}^{\gamma} \circ f_{\alpha}^{\beta}$.
We will assume that the index set $A$ is either $\Rspace$ or
$\Rspace_{\geq 0}$ and not explicitly state indices unless necessary.

A real number $\alpha$ is said to be a \emph{regular value} of the persistence module $\Fcal$ if there exists some $\epsilon > 0$
such that the map $f_{\alpha - \delta}^{\alpha + \delta}$ is an isomorphisms for each $\delta < \epsilon$.
Otherwise we say that $\alpha$ is a \emph{critical value} of the persistence module; if $A = \Rspace_{\geq 0}$, then $\alpha = 0$ will always be considered to be a critical value.
We say that $\Fcal$ is \emph{tame} if it has a finite number of critical values and if all the vector spaces $F_{\alpha}$ are of finite rank.
Any tame $\Rspace_{\geq 0}$-module $\Fcal$ must have a smallest non-zero critical value $\rho(\Fcal)$; we call this number the \emph{feature size} of the persistence module.

\begin{figure}[tbp]
\begin{center}
\includegraphics[scale=0.35]{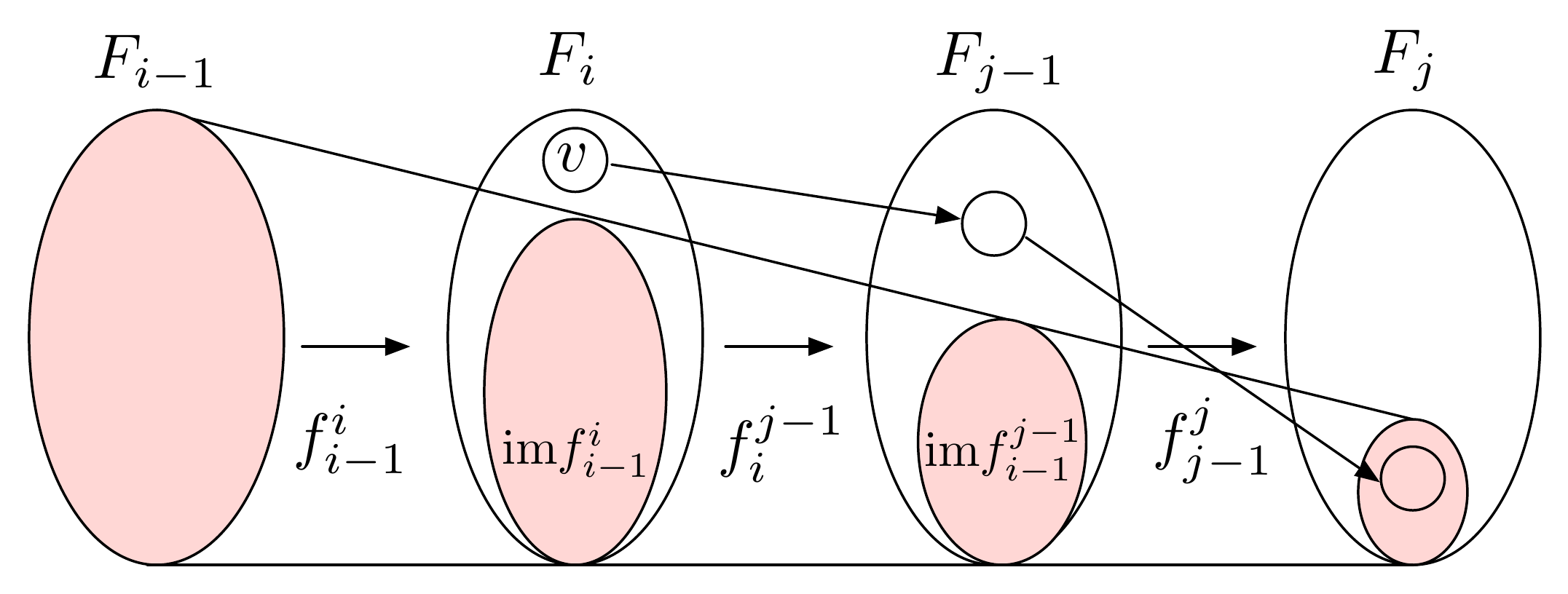} 
\end{center}
\caption[Persistence modules.]{The vector $v$ is born at level $i$ and then it dies at level $j$.}
\label{fig:pmodule}
\end{figure}

Assume $\Fcal$ is tame and so we have a finite ordered list of critical values $0 = c_0  < c_1 < \ldots < c_m$.
We choose regular values $\{a_i\}_{i=0}^{m}$ such that $c_{i-1} < a_{i-1} < c_i < a_i$ for all $1 \leq i  \leq m,$
and we adopt the shorthand notation $F_i \equiv F_{a_i}$ and $f_i^j : F_i \to F_j$, for $0 \leq i \leq j \leq m$.
A vector $v \in \Fcal_i$ is said to be \emph{born} at level $i$ if $v \not \in \image f_{i-1}^i$,
and such a vector \emph{dies} at level $j$ if $f_i^j(v) \in \image f_{i-1}^j$ but $f_i^{j-1}(v) \not \in \image f_{i-1}^{j-1}$.
This is illustrated in Figure \ref{fig:pmodule}.
We then define $P^{i,j}$ to be the vector space of vectors that are born at level $i$ and then subsequently die at level $j$, and $\beta^{i,j}$ denotes its rank.

\subsubsection{Persistence Diagrams}

\begin{figure}[tbp]
 \begin{center}
  \includegraphics[scale=0.35]{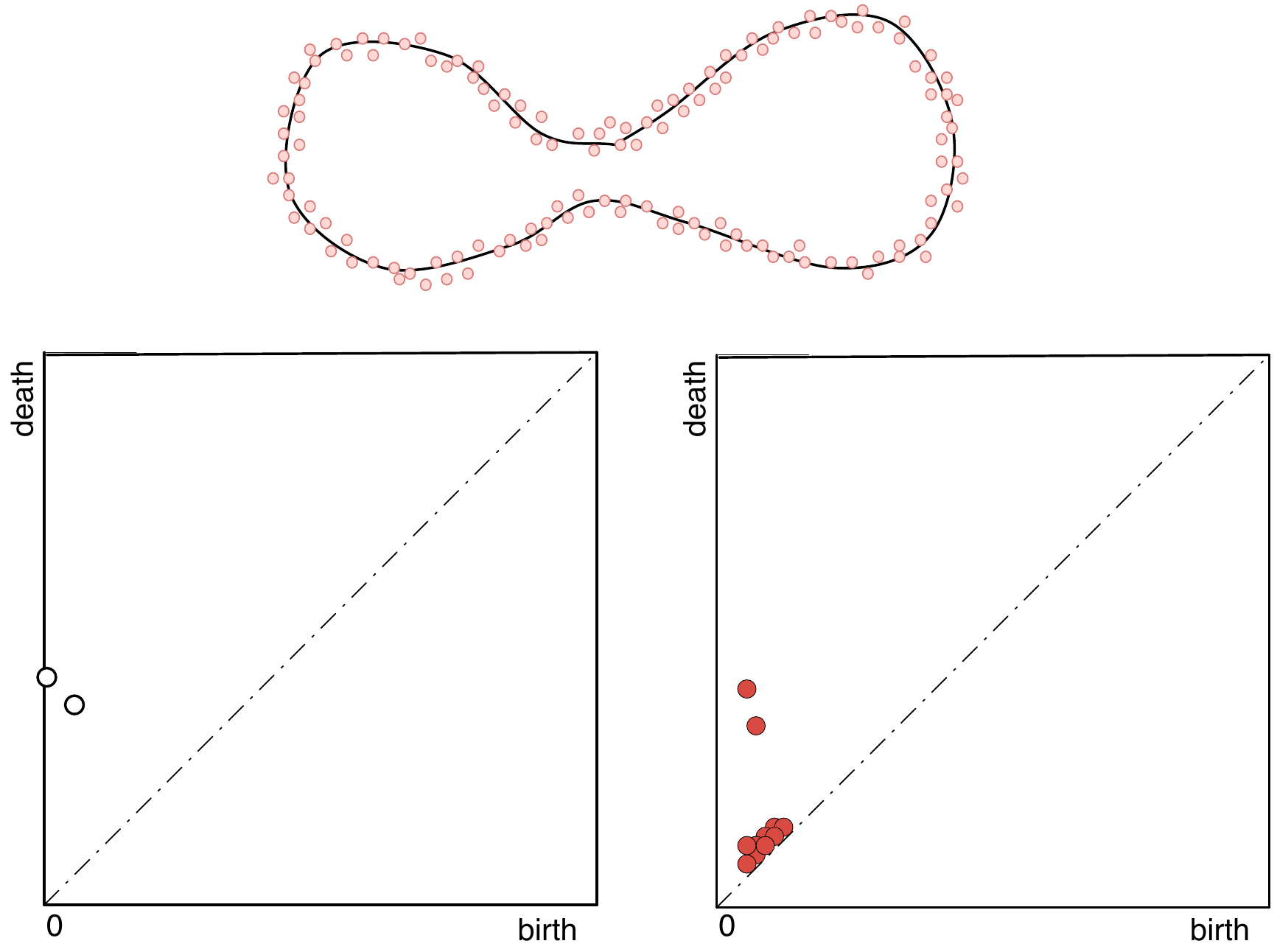}
 \end{center}
\caption[Illustration of a point cloud and its persistence diagram.]
{Illustration of a point cloud and its persistence diagram. 
Top: $\mathbb{X}$ is the curve embedded as shown in the plane and $U$ is the point cloud. 
Bottom left: the persistence diagram ${\rm Dgm}_1(d_{\mathbb{X}})$; Bottom right: the persistence diagram ${\rm Dgm}_1(d_{U})$.}
\label{fig:distance}
\end{figure}

The information contained within a tame module $\Fcal$ is often compactly represented by a \emph{persistence diagram}, $\Ddgm{}{\Fcal}$.
This diagram is a multi-set of points in the extended plane. It contains $\beta^{i,j}$ copies of the points $(c_i,c_j)$, as well
as infinitely many copies of each point along the major diagonal $y=x$. In Figure
 \ref{fig:distance} the persistence diagrams for a curve and a point cloud sampled from it
 are displayed; see Section \ref{subsec:examples} for a full explanation of this figure.

For any two points $u = (x,y)$ and $u' = (x',y')$ in the extended plane, we define $||u - u'||_{\infty} = \max\{|x - x'|, |y - y'|\}$.
We define the \emph{bottleneck distance} between any two persistence diagrams $D$
and $D'$ to be:
$$d_B(D,D') = \inf_{\Gamma: D \to D'} \sup_{u \in D} ||u - \Gamma(u)||_{\infty},$$
where $\Gamma$ ranges over all bijections from $D$ to $D'$.
Under certain conditions which we now describe, persistence diagrams will be stable under the bottleneck distance.

\begin{figure}[tbp]
\begin{center}
\includegraphics[scale=0.35]{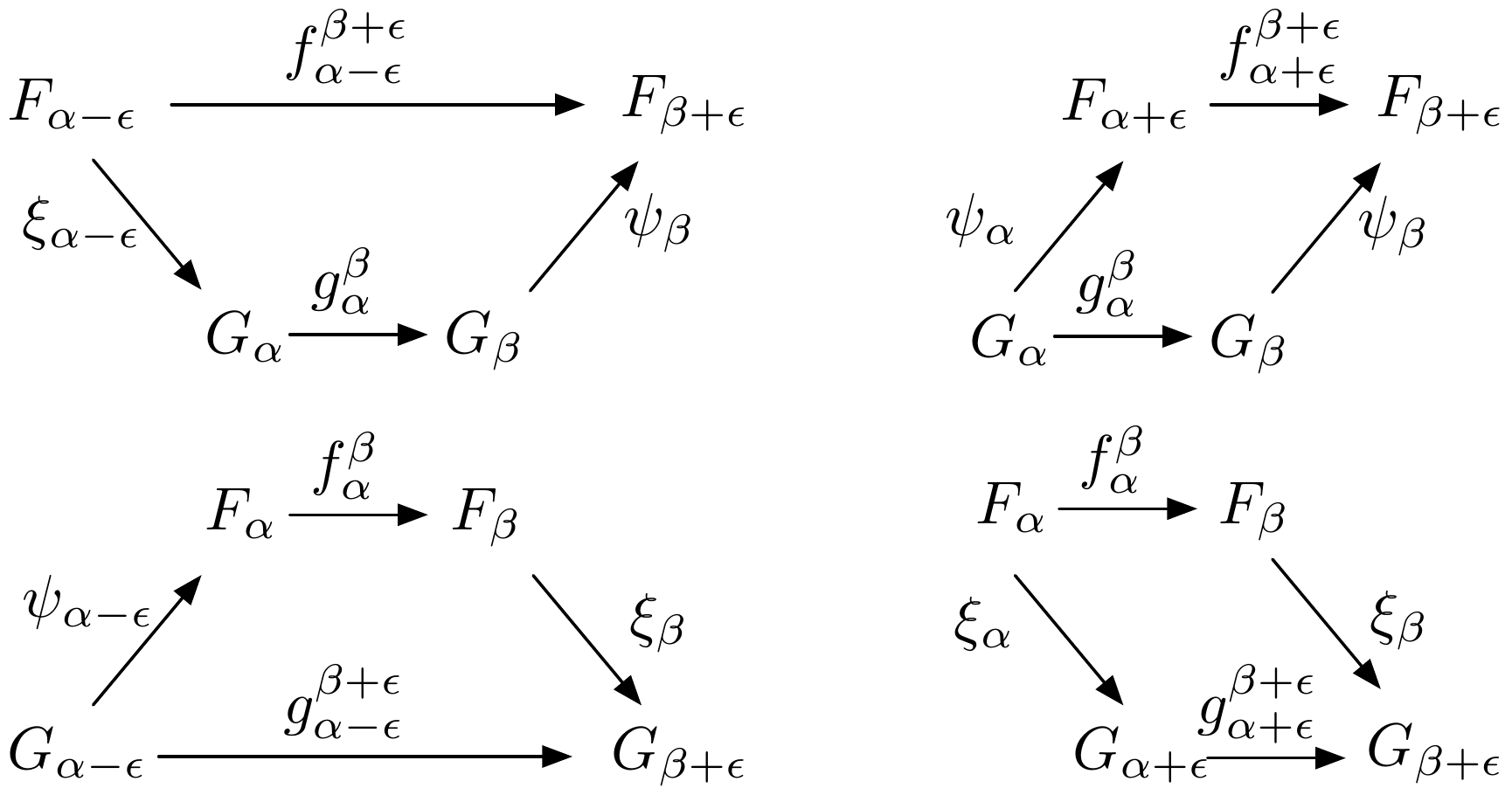} 
\end{center}
\caption[Commuting diagrams for strongly interleaving persistence modules.]{Commuting diagrams for strongly interleaving persistence modules.}
\label{fig:interleave}
\end{figure}

Two persistence modules $\Fcal$ and $\Gcal$ are said to be \emph{strongly $\epsilon$-interleaved} if, for some positive $\epsilon$,
there exist two families $\{\xi_{\alpha} : F_{\alpha} \to G_{\alpha + \epsilon}\}_{\alpha}$ and $\{\psi_{\alpha}: G_{\alpha} \to F_{\alpha + \epsilon}\}$ of
linear maps which commute with the module maps $\{f_{\alpha}^{\beta}\}$ and $\{g_{\alpha}^{\beta}\}$ in the appropriate manner.
More precisely, we require that, for each $\alpha \leq \beta$, the four diagrams in Figure \ref{fig:interleave}.all commute.

We can now state the diagram stability result (\cite{ChaCohGli2009}, Theorem 4.4), that we will need later in this paper.
\begin{theorem}[Diagram Stability Theorem]
Let $\Fcal$ and $\Gcal$ be tame persistence modules and $\epsilon > 0$. If $\Fcal$ and $\Gcal$ are strongly $\epsilon$-interleaved, then
$$d_B(\Ddgm{}{\Fcal}, \Ddgm{}{\Gcal}) \leq \epsilon.$$
\label{result:PDS}
\end{theorem}

When we wish to compute the persistence diagram associated to a module $\Fcal$, it is often
convenient to substitute another module $\Gcal$, usually one defined in terms of simplicial
complexes or other computable objects.
The following theorem (\cite{EdeHar2010}, p.159) gives a condition under which this is possible.
\begin{theorem}[Persistence Equivalence Theorem]
Given two persistence modules $\Fcal$ and $\Gcal$, suppose there exist for each $\alpha$ isomorphisms $F_{\alpha} \cong G_{\alpha}$
which commute with the module maps,
then $\Ddgm{}{\Fcal} = \Ddgm{}{\Gcal}$.
\label{result:PET}
 \end{theorem}

\subsubsection{(Co)Kernel Modules}
Suppose now that we have two persistence modules $\Fcal$ and $\Gcal$ along with a family of maps $\{\phi_{\alpha}: F_{\alpha} \to G_{\alpha}\}$
which commute with the module maps -- for every pair $\alpha \leq \beta$, we have $g_{\alpha}^{\beta} \circ \phi_{\alpha} = \phi_{\beta} \circ f_{\alpha}^{\beta}$.
In other words, every square commutes in the diagram below:
\begin{align*}
\ldots \to  &F_{\alpha}  \xrightarrow{f_{\alpha}^{\beta}} F_{\beta}  \to \ldots \\
   &~\downarrow   \phi_{\alpha}  ~~~\downarrow \phi_{\beta}\\
   \ldots \to & G_{\alpha}  \xrightarrow{g_{\alpha}^{\beta}}  G_{\beta}  \to \ldots
\end{align*}
Then, for each pair of real numbers $\alpha \leq \beta$, the restriction of $f_{\alpha}^{\beta}$ to $\kernel \phi_{\alpha}$ maps into $\kernel \phi_{\beta},$
giving rise to a new kernel persistence module, with persistence diagram denoted by $\Ddgm{}{\kernel \phi}$.
Similarly, we obtain a cokernel persistence module, with diagram $\Ddgm{}{\cokernel \phi}$.

\subsection{Homology}
\label{subsec:examples}

Our main examples of persistence modules all come from homology groups, either absolute or relative, and the various maps between them.
Homology persistence modules can arise from families of topological spaces $\{\Xspace_{\alpha}\}$, along
with inclusions $\Xspace_{\alpha} \hookrightarrow \Xspace_{\beta}$ for all $\alpha \leq \beta$.
Whenever we have such a family, the inclusions induce maps $\Hgroup_j(\Xspace_{\alpha}) \to \Hgroup_j(\Xspace_{\beta})$,
for each homological dimension $j \geq 0$, and hence we have persistence modules for each $j$.
Defining $\Hgroup(\Xspace_{\alpha}) = \bigoplus_j \Hgroup_j(\Xspace_{\alpha})$ and taking direct sums of maps in
the obvious way, will also give one large direct-sum persistence module
$\{\Hgroup(\Xspace_{\alpha})\}$.

\subsubsection{Distance Functions}

Here, the families of topological spaces will be produced by the sublevel sets of distance functions. 
Given a topological space $\Xspace$ embedded in some Euclidean space $\Rspace^N$, we define $d_{\Xspace}$
as the distance function which maps each point in the ambient space to the distance from its closest point in $\Xspace$.
More formally, 
for each $y \in \Rspace^N$,
$d_{\Xspace}(y) = \inf_{x \in \Xspace} \distone(x,y).$
We let $\Xspace_{\alpha}$ denote the sublevel set $d_{\Xspace}^{-1}[0,\alpha]$; each
sublevel set should be thought of as a thickening of $\Xspace$ within the ambient space.
Increasing the thickening parameter produces a growing family of sublevel sets, giving
rise to the persistence module $\{\Hgroup(\Xspace_{\alpha})\}_{\alpha \in \Rspace \geq 0}$;
we denote the persistence diagram of this module by $\Ddgm{}{d_{\Xspace}}$ and use
$\Ddgm{j}{d_{\Xspace}}$ for the diagrams of the individual modules for each homological dimension $j$.

In Figure \ref{fig:distance}, we see an example of such an $\Xspace$ embedded in the plane, along
with the persistence diagram $\Ddgm{1}{d_{\Xspace}}$.
We also have the persistence diagram $\Ddgm{1}{d_{\Uspace}}$, where $\Uspace$ is a dense point sample of $\Xspace$.
Note that the two diagrams are quite close in bottleneck distance.
Indeed, the difference between the two diagrams will always be upper-bounded by the Hausdorff distance between the space and its sample; this follows
from Theorem \ref{result:PDS}.

\begin{figure}[tbp]
 \begin{center}
  \includegraphics[scale=0.35]{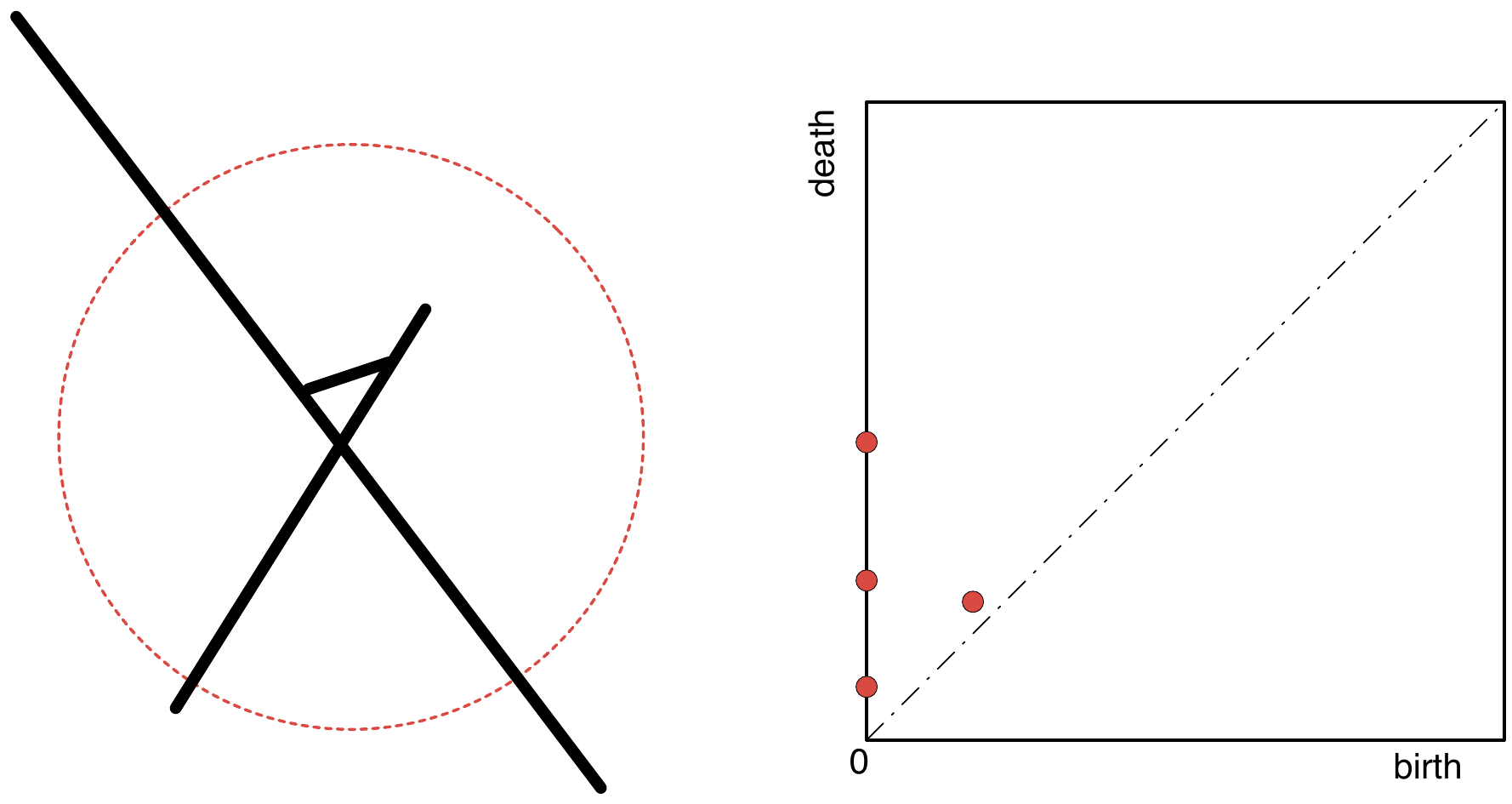}
 \end{center}
\caption[Persistence diagram for relative homology.]{Left: The space $\mathbb{X}$ is in solid line and the closed ball $B$ has dotted boundary.
        Right: the persistence diagram for the module $\{\sf{H}_1$$(\mathbb{X}_{\alpha} \cap B, \mathbb{X}_{\alpha} \cap \partial B)\}$.}
\label{fig:relhom}
\end{figure}

Persistence modules of relative homology groups also arise from families of pairs of spaces, as the next example shows.
Referring to the left part of Figure \ref{fig:relhom}, we let $\Xspace$ be the space drawn in solid lines and $B$ the closed ball whose boundary is drawn as a dotted circle.
By restricting $d_{\Xspace}$ to $B$ and also to $\partial B$, we produce pairs of sublevel sets $(\Xspace_{\alpha} \cap B, \Xspace_{\alpha} \cap \partial B)$.
Using the maps induced by the inclusions of pairs, we obtain the persistence module
$\{\Hgroup(\Xspace_{\alpha} \cap B, \Xspace_{\alpha} \cap \partial B)\}_{\alpha \in \Rspace_{\geq 0}}$ of relative homology groups. The persistence diagram, for homological dimension $1$, appears
on the right half of Figure \ref{fig:relhom}.

\subsection{Stratified Spaces}

\begin{figure}[tbp]
 \begin{center}
  \includegraphics[scale=0.25]{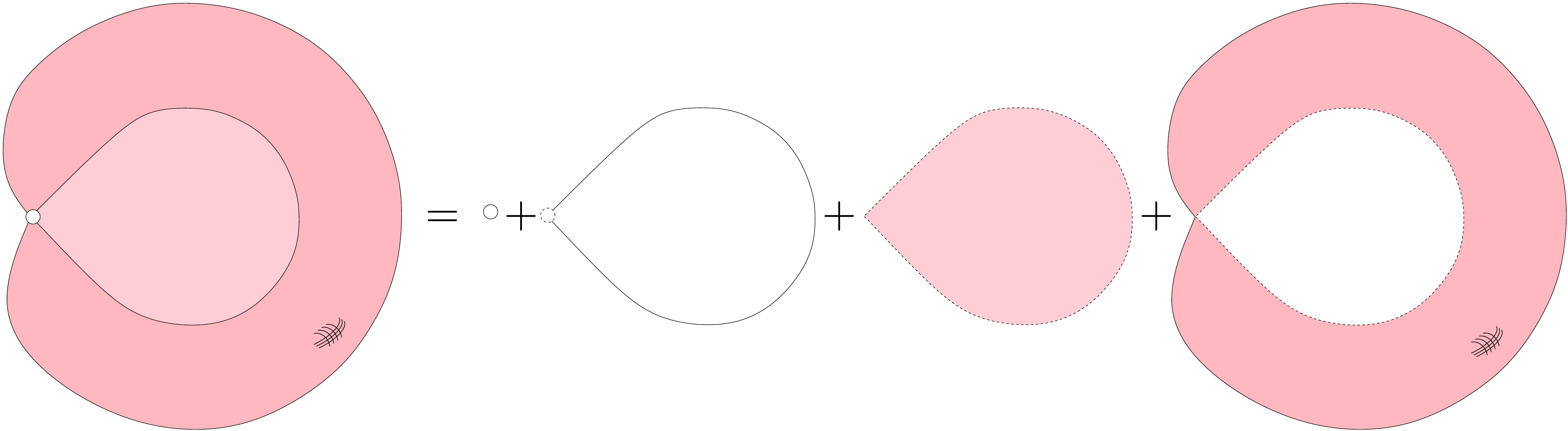}
 \end{center}
\caption[Example of the stratification of a pinched torus with a spanning disc stretched across the hole.]{The coarsest stratification of a pinched torus with a spanning disc stretched across the hole.}
\label{fig:ptorus}
\end{figure}

We assume that we have a topological space $\Xspace$ embedded in some Euclidean space $\Rspace^N$.
A (purely) \emph{$d$-dimensional stratification of $\Xspace$} is
a decreasing sequence of closed subspaces
$$\Xspace = \Xspace_d \supseteq \Xspace_{d-1} \supseteq \ldots \Xspace_0 \supseteq \Xspace_{-1} = \emptyset,$$
such that for each $i$, the $i$-dimensional \emph{stratum} $\Sspace_i = \Xspace_i - \Xspace_{i-1}$ is a (possibly empty) $i$-manifold.
The connected components of $\Sspace_i$ are called $i$-dimensional \emph{pieces}.
This is illustrated in Figure \ref{fig:ptorus}, where the space $\Xspace$ is a pinched torus
with a spanning disc stretched across the hole.

\begin{figure}[tbp]
 \begin{center}
  \includegraphics[scale=0.25]{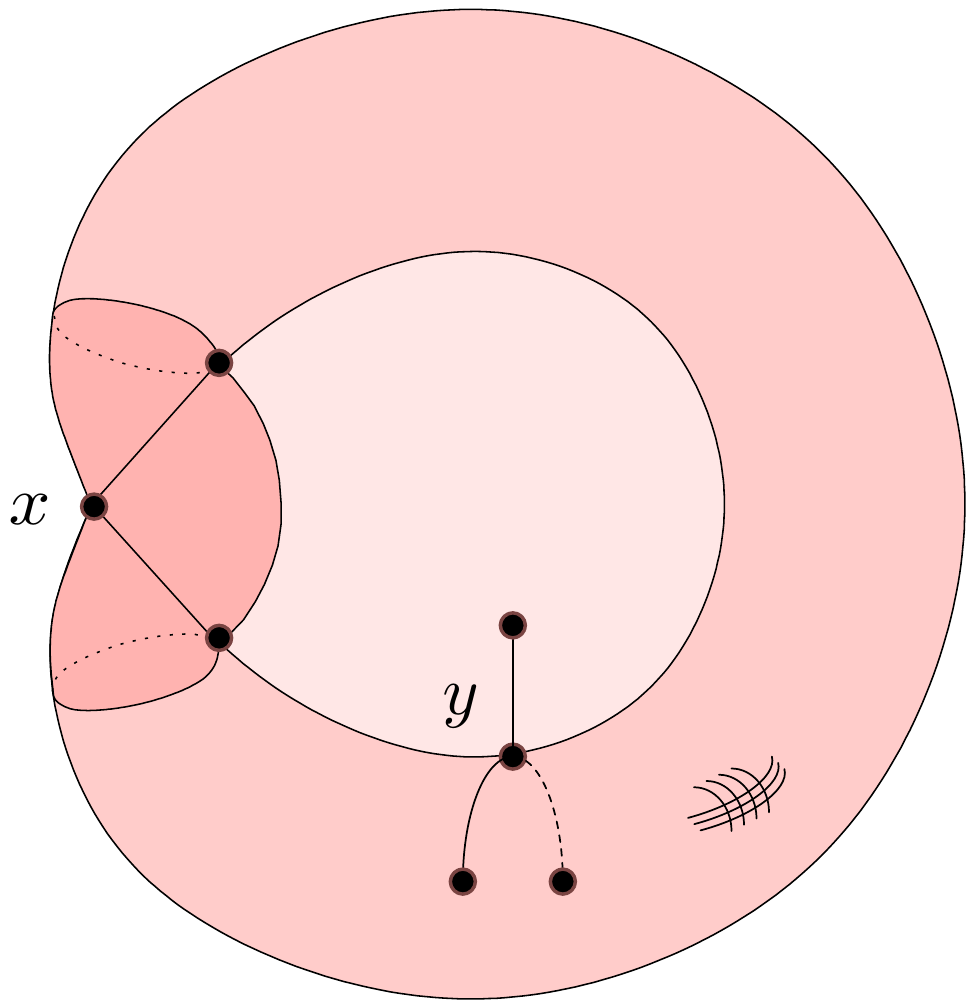}
 \end{center}
\caption[A $2$-dimensional stratified space and corn construction.] {The cones $c(L_x)$ and $c(L_y)$, where $x$ and $y$ are respectively in the $0$-stratum and the $1$-stratum,
are highlighted.}
\label{fig:stratified}
\end{figure}

One usually also imposes a requirement to ensure that the various pieces fit together uniformly.
There are a number of different ways this can be done
(see \cite{HugWei2000} for an extensive survey).
For example, one might assume that for each $x \in \Sspace_i$, there exists a small enough neighborhood $N(x) \subseteq \Xspace$
and a $(d-i-1)$-dimensional stratified space $L_x$ such that $N(x)$ is stratum-preserving homeomorphic to the product
of an $i$-ball and the cone on $L_x$; one can then show that the space $L_x$ depends only on the particular piece containing $x$.
This definition is illustrated in Figure \ref{fig:stratified}.

Since the topology on $\Xspace$ is that inherited from the ambient space, this neighborhood $N(x)$ will
take the form $\Xspace \cap B_r(x)$, where $B_r(x)$ is a small enough ball around $x$ in the ambient space.

We note that the above definition requires all strata to be contained within the closure of the top-dimensional stratum.
It is also possible, of course, to have spaces where this is not the case: for example, a two-dimensional plane that has been punctured by a line.
In this case, a slight adjustment to the above definitions can be made in order to impose similar notions of uniformity.

\subsubsection{Local Homology and Homology Stratifications}

Recall (\cite{Mun1984}) that the local homology groups of a space $\Xspace$ at a point $x \in \Xspace$ are the groups
$\Hgroup_i(\Xspace, \Xspace - x)$ in each homological dimension $i$.
If $\Xspace$ happens to be a $d$-manifold, or if $x$ is simply a point in the top-dimensional stratum of a $d$-dimensional stratification,
then these groups are rank one in dimension $d$ and trivial
in all other dimensions.
On the other hand, the local homology groups for lower-stratum points can be more interesting; for example
if $x$ is the crossing point in Figure \ref{fig:crosskernel}, then $\Hgroup_1(\Xspace, \Xspace - x)$ has rank three.

If $x$ and $y$ are close enough points in a particular piece of the same stratum, then there is a natural
isomorphism between their local homology groups $\Hgroup(\Xspace, \Xspace - x) \cong \Hgroup(\Xspace, \Xspace - y)$,
which can be understood in the following manner.
Taking a small enough radius $r$ and using excision, we see that the two local homology groups in question are in fact just
$\Hgroup(\Xspace \cap B_r(x), \Xspace \cap \partial B_r(x))$ and $\Hgroup(\Xspace \cap B_r(y), \Xspace \cap \partial B_r(y))$.
Both of these groups will then map, via intersection of chains, isomorphically into the group
$\Hgroup(\Xspace \cap B_r(x) \cap B_r(y), \partial(B_r(x) \cap B_r(y))$, and the isomorphism above is then derived from these
two maps.
See the points in Figure \ref{fig:crosskernel} for an illustration of this idea.

In \cite{RouSan1999}, the authors define the concept of a homology stratification of a space $\Xspace$.
Briefly, they require a decomposition of $\Xspace$ into pieces such that the locally homology groups are locally constant across each piece;
more precisely, that the maps discussed above be isomorphisms for each pair of close enough points in each piece.


\section{Topological Inference Theorem}
\label{sec:TIT}

From the discussion above, it is easy to see that
any stratification of a topological space will also be a homology stratification.
The converse is unfortunately false. However, we can build a useful analytical tool
based on the contrapositive: given two points in a point cloud, we can hope to state, based on
their local homology groups and the maps between them, that the two points
should not be placed in the same piece of any stratification. To do this, we first adapt
the definition of these local homology maps into a more multi-scale and robust framework.
More specifically, we introduce a radius parameter $r$ and a notion of local equivalence, $\sim_r$, which allows
us to group the points of $\Xspace$, as well as of the ambient space, into strata at this radius scale.
We then give the main result of this section: topological conditions under which the point cloud $U$ can be used to infer
the strata at different radius scales.

\subsection{Local Equivalence}
We assume that we are given some topological space $\Xspace$ embedded in some Euclidean space in $\Rspace^N$.
For each radius $r \geq 0$, and for each pair of points $p,q \in \Rspace^N$, we define the following homology map $\phi^{\Xspace}(p,q,r)$:
\begin{align}
\label{eq:intersection}
\Hgroup(\Xspace \cap B_r(p), \Xspace \cap \partial B_r(p)) \to \Hgroup(\Xspace \cap B_r(p) \cap B_r(q), \Xspace \cap \partial(B_r(p) \cap B_r(q))).
\end{align}
Intuitively, this map can be understood as taking a chain, throwing away the parts that lie outside the smaller
range, and then modding out the new boundary. Alternatively, one may think of it as being induced by a
combination of inclusion and excision. 
A formal definition is given in Appendix A. 

\begin{figure}[tbp]
 \begin{center}
  \includegraphics[scale=0.35]{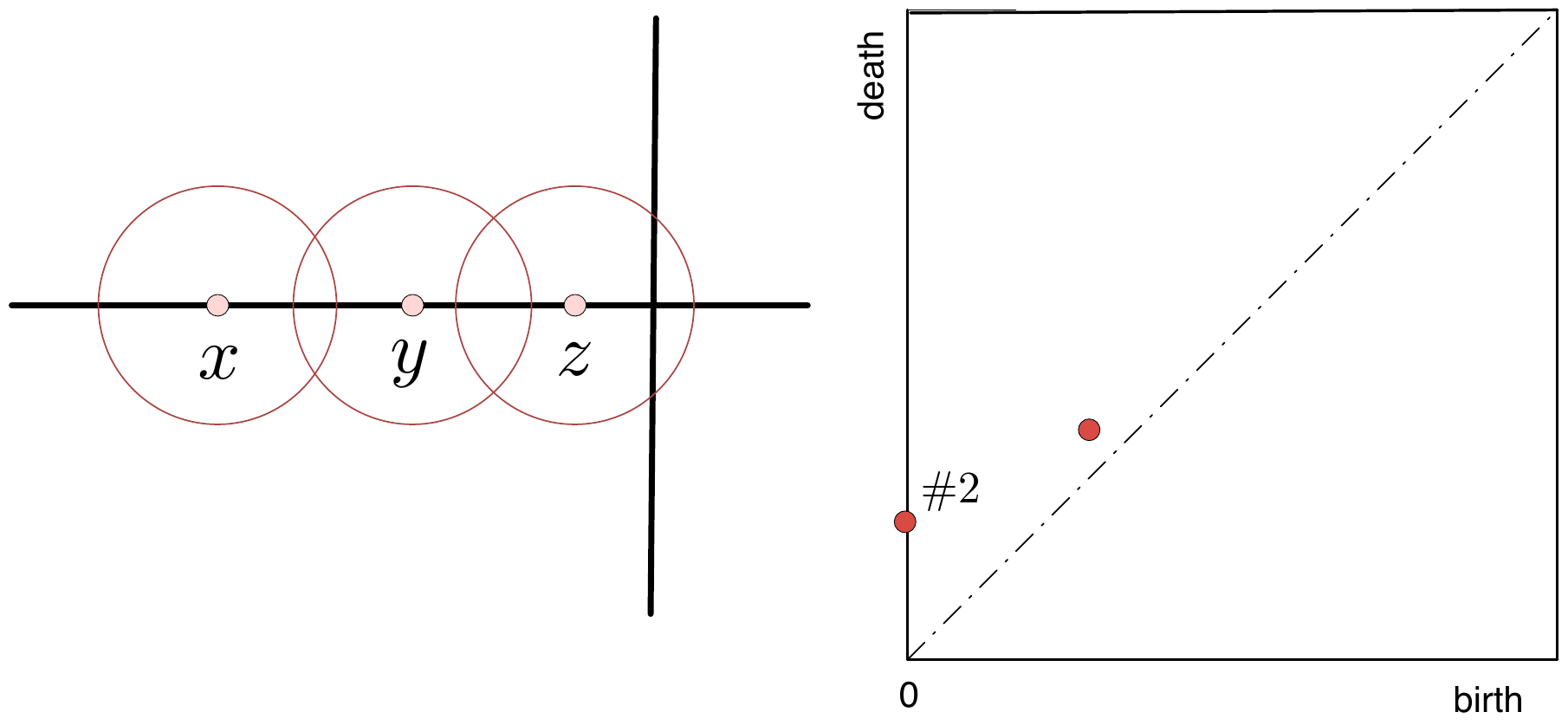}
 \end{center}
\caption[Illustration of equivalence relation.]{Left: $x \sim_r y$, $y \nsim_r z$. 
Right: the $1$-dim persistence diagram, for the kernel of the map going from the $z$ ball into its intersection with the $y$ ball.
A number, i.e., $\#2$, labeling a point in the persistence diagram indicates its multiplicity. 
}
\label{fig:crosskernel}
\end{figure}

Using these maps, we impose an equivalence relation on $\Rspace^N$.
\begin{definition}[Local equivalence]
\label{eqn:locequiv}
Two points $x$ and $y$ are said to have equivalent local structure at radius $r$,
denoted $x \sim_r y$, iff there exists a chain of points $x = x_0, x_1, \ldots, x_m = y$ from
$\Xspace$ such that, for each $1 \leq i \leq m$,
the maps $\phi^{\Xspace}(x_{i-1},x_i,r)$ and $\phi^{\Xspace}(x_i,x_{i-1},r)$ are both isomorphisms.
\end{definition}
In other words, $x$ and $y$ have the same local structure at this radius iff they can be connected by a chain of points which are pairwise close enough
and whose local homology groups at radius $r$ map into each other via intersection.
Different choices of $r$ will of course lead to different equivalence classes.
For example, consider the space $\Xspace$ drawn in the plane as shown in the left half of Figure \ref{fig:crosskernel}.
At the radius drawn, point $z$ is equivalent to the cross point and is not equivalent to either the point $x$ or $y$.
Note that some points from the ambient space will now be considered equivalent to $x$ and $y$, and some others will be equivalent
to $z$.

On the other hand, a smaller choice of radius would result in all three of $x, y,$ and $z$ belonging to the same equivalence class.

\subsubsection{(Co)Kernel Persistence}

In order to relate the point cloud $U$ to the equivalence relation $\sim_r$, we
must first define a multi-scale version of the maps $\phi^{\Xspace}(p,q,r)$; we do so by gradually thickening the space $\Xspace$.
Let $d_{\Xspace}: \Rspace^N \to \Rspace$ denote the function which maps each point in the ambient space
to the distance from its closest point on $\Xspace$. For each $\alpha \geq 0$, we define $\Xspace_{\alpha} = d_{\Xspace}^{-1}[0,\alpha]$.
For each $p,q,$ and $r,$ we will consider the intersection map $\phi^{\Xspace}_{\alpha}(p,q,r)$,
which is defined by substituting $\Xspace_{\alpha}$ for $\Xspace$ in (\ref{eq:intersection}).
Note of course that $\phi^{\Xspace}(p,q,r) = \phi^{\Xspace}_0(p,q,r)$.

For the moment, we fix a choice of $p,q,$ and $r$, and we use the following shorthand:
\begin{align*}
B_{p}^{\Xspace}(\alpha) &= \Xspace_{\alpha} \cap B_r(p), \\
\bdr B_{p}^{\Xspace}(\alpha) &= \Xspace_{\alpha} \cap \bdr B_r(p), \\
B_{pq}^{\Xspace}(\alpha) &= \Xspace_{\alpha} \cap B_r(p) \cap B_r(q), \\
\bdr B_{pq}^{\Xspace}(\alpha) &= \Xspace_{\alpha} \cap \bdr (B_r(p) \cap B_r(q)).
\end{align*}
and we also often write $B_{p}^{\Xspace} = B_{p}^{\Xspace}(0)$ and $B_{pq}^{\Xspace} = B_{pq}^{\Xspace}(0)$.
By replacing $\Xspace$ with $U$ in this shorthand, we also write $B_p^U(\alpha) = U_{\alpha} \cap B_r(p)$, and so forth.

For any pair of non-negative real values $\alpha \leq \beta$ the inclusion $\Xspace_{\alpha} \hookrightarrow \Xspace_{\beta}$ gives rise to the following
commutative diagram:
\begin{align}
\label{diag:K}
\Hgr(B_{p}^{\Xspace}(\alpha), \bdr B_{p}^{\Xspace}(\alpha))
 &\xrightarrow  {\phi_{\alpha}^{\Xspace}}
\Hgr(B_{pq}^{\Xspace}(\alpha), \bdr B_{pq}^{\Xspace}(\alpha))\notag\\
 \downarrow ~~~~~~~~~~~~&~~~~~~~~~~~~\downarrow~~~~~~~\notag\\
 \Hgr(B_{p}^{\Xspace}(\beta), \bdr B_{p}^{\Xspace}(\beta))
 &\xrightarrow  {\phi_{\beta}^{\Xspace}}
\Hgr(B_{pq}^{\Xspace}(\beta), \bdr B_{pq}^{\Xspace}(\beta))
\end{align}

Hence there are maps $\kernel \phi^{\Xspace}_{\alpha} \to \kernel \phi^{\Xspace}_{\beta}$ and $\cokernel \phi^{\Xspace}_{\alpha} \to \cokernel \phi^{\Xspace}_{\beta}$.
Allowing $\alpha$ to increase from $0$ to $\infty$ gives rise to two persistence modules,
$\{\kernel \phi_{\alpha}^{\Xspace}\}$ and $\{\cokernel \phi_{\alpha}^{\Xspace}\}$,
with diagrams $\Ddgm{}{\kernel \phi^{\Xspace}}$
and $\Ddgm{}{\cokernel \phi^{\Xspace}}$.
Recall that a homomorphism is an isomorphism iff its kernel and cokernel are both zero.
In our context then, the map $\phi^{\Xspace}$ is an isomorphism iff neither $\Ddgm{}{\kernel \phi^{\Xspace}}$
nor $\Ddgm{}{\cokernel \phi^{\Xspace}}$ contain any points on the $y$-axis above $0$.

\paragraph{Example.}
As shown in the left part of Figure \ref{fig:crosskernel}, $x$, $y$ and $z$ are points sampled from a cross embedded in the plane. Taking $r$ as drawn,
we note that the right part of the figure displays $\Ddgm{1}{\kernel \phi^{\Xspace}}$, where $\phi^{\Xspace} = \phi^{\Xspace}(z,y,r)$; we now explain this diagram in some detail.
The group $\Hgroup_1(B_{z}^{\Xspace}, \bdr B_{z}^{\Xspace})$ 
has rank three; as a possible basis we might take
the three classes represented by the horizontal line across the ball, the vertical line across the ball, and the two short segments defining
the northeast-facing right angle.
Under the intersection map $\phi^{\Xspace} = \phi^{\Xspace}_0$, the first of these classes maps to the generator of
$\Hgroup_1(B_{zy}^{\Xspace}, \bdr B_{zy}^{\Xspace})$,
while the other two map to zero.
Hence $\kernel \phi^{\Xspace}_0$ has rank two.
Both classes in this kernel eventually die, one at the $\alpha$ value which fills in the northeast corner of the larger ball, and
the other at the $\alpha$ value which fills in the entire right half; these two values are the same here due to symmetry in the picture.
At this value, the map $\phi^{\Xspace}_{\alpha}$ is an isomorphism and it remains so until the intersection of the two balls
fills in completely.
This gives birth to a new kernel class which subsequently dies when the larger ball finally fills in.
The diagram $\Ddgm{1}{\kernel \phi^{\Xspace}}$ thus contains three points; the leftmost two show that the map
$\phi^{\Xspace}$ is not an isomorphism.

\subsection{Inference Theorem}

Given a point cloud $\Uspace$ sampled from $\Xspace$ consider the following question: for a radius $r$, how
can we infer whether or not any given pair of points in $\Uspace$ has the same local structure at this radius? 
In this subsection, we prove a theorem which describes the circumstances under which we can
make the above inference. Naturally, any inference will require that we use $\Uspace$ to judge whether or not the
maps $\phi^{\Xspace}(p,q,r)$ are isomorphisms. The basic idea is that if $\Uspace$ is a dense
enough sample of $\Xspace$, then the (co)kernel diagrams defined by $\Uspace$ will be good enough
approximations of the diagrams defined by $\Xspace$.

\subsubsection{(Co)Kernel Stability}

Again we fix $p,q,$ and $r$, and write $\phi^{\Xspace} = \phi^\Xspace(p,q,r)$.
For each $\alpha \geq 0$, we let $\Uspace_{\alpha} = d_{\Uspace}^{-1}[0,\alpha]$.
We consider $\phi^{\Uspace}_{\alpha} = \phi^{\Uspace}_{\alpha}(p, q, r)$, defined by replacing
$\Xspace$ with $U_{\alpha}$ in (\ref{eq:intersection}).
Running $\alpha$ from $0$ to $\infty$, we obtain two more persistence modules,
$\{\kernel \phi_{\alpha}^{\Uspace}\}$ and $\{\cokernel \phi_{\alpha}^{\Uspace}\}$, 
with diagrams
$\Ddgm{}{\kernel \phi^{\Uspace}}$ and $\Ddgm{}{\cokernel \phi^{\Uspace}}$.

If $\Uspace$ is a dense enough sample of $\Xspace$, then the (co)kernel diagrams defined by $\Uspace$ will be good
approximations of the diagrams defined by $\Xspace$.
More precisely, we have the following easy consequence of Theorem \ref{result:PDS}:
\begin{theorem}[(Co)Kernel Diagram Stability] The bottleneck distances between the (co)kernel diagrams
of $\phi^U$ and $\phi^{\Xspace}$ are upper-bounded by the Hausdorff distance between $U$ and $\Xspace$:
\begin{align*}
d_B(\Ddgm{}{\kernel \phi^{\Uspace}}, \Ddgm{}{\kernel \phi^{\Xspace}}) & \leq d_H(\Uspace, \Xspace),\\
d_B(\Ddgm{}{\cokernel \phi^{\Uspace}}, \Ddgm{}{\cokernel \phi^{\Xspace}}) & \leq d_H(\Uspace, \Xspace).
\end{align*}
\label{result:CKDS}
\end{theorem}
\begin{proof}
We prove the first inequality; the proof of the second is identical.
Put $\epsilon = d_H(\Uspace,\Xspace)$.
Then, for each $\alpha \geq 0$, the inclusions $\Uspace_{\alpha} \hookrightarrow \Xspace_{\alpha + \epsilon}$ and
$\Xspace_{\alpha} \hookrightarrow \Uspace_{\alpha + \epsilon}$ induce maps $\kernel \phi^{\Uspace}_{\alpha} \to \kernel \phi^{\Xspace}_{\alpha + \epsilon}$
and $\kernel \phi^{\Xspace}_{\alpha} \to \kernel \phi^{\Uspace}_{\alpha + \epsilon}$.
These maps clearly commute with the module maps in the needed way, and hence we have
the required $\epsilon$-interleaving and can thus appeal to Theorem \ref{result:PDS}.
\eop
\end{proof}

\subsubsection{Main Inference Result}

We now suppose that we have a point sample $U$ of a space $\Xspace$, where the Hausdorff distance between the two is no more than some $\epsilon$;
in this case, we call $\Uspace$ an \emph{\ep-approximation} of $\Xspace$.
Given two points $p,q \in U$ and a fixed radius $r$, we set $\phi^{\Xspace} = \phi^{\Xspace}(p,q,r)$, and we wish to determine
whether or not $\phi^{\Xspace}$ is an isomorphism.
Since we only have access to the point sample $U$, we instead compute the diagrams $\Ddgm{}{\kernel \phi^{\Uspace}}$ and $\Ddgm{}{\cokernel \phi^{\Uspace}}$;
we provide an algorithm for doing this in Section \ref{sec:Alg}.
The main Theorem of this section, Theorem \ref{result:IIT}, gives conditions under which these diagrams enable us to answer the isomorphism question for $\phi^{\Xspace}$.
To state the theorem we first need some more definitions.

Given any persistence diagram $\Dcal$, which we recall is a multi-set of points in the extended plane, and two positive real numbers $a < b$, we let $\Dcal(a,b)$ denote
the intersection of $\Dcal$ with the portion of the extended plane which lies above $y = b$ and to
the left of $x = a$; note that these points correspond to classes which are born no later than $a$ and die no
earlier than $b$.

For a fixed choice of $p,q,r$, we consider the following two persistence modules:
$\{\Hgroup(B_{p}^{\Xspace}(\alpha), \bdr B_{p}^{\Xspace})\}$ and 
$\{\Hgroup(B_{pq}^{\Xspace}(\alpha), \bdr B_{pq}^{\Xspace})\}$.
We let $\sigma(p,r)$ and $\sigma(p,q,r)$ denote their respective feature sizes and then
set $\rho(p,q,r)$ to their minimum. 

We now give the main theorem of this section, which states that we can use $\Uspace$ to decide whether or not $\phi^{\Xspace}(p,q,r)$ is an isomorphism as long as 
$\rho(p,q,r)$ is large enough relative to the sampling density.
\begin{figure}[tbp]
 \begin{center}
  \includegraphics[scale=0.35]{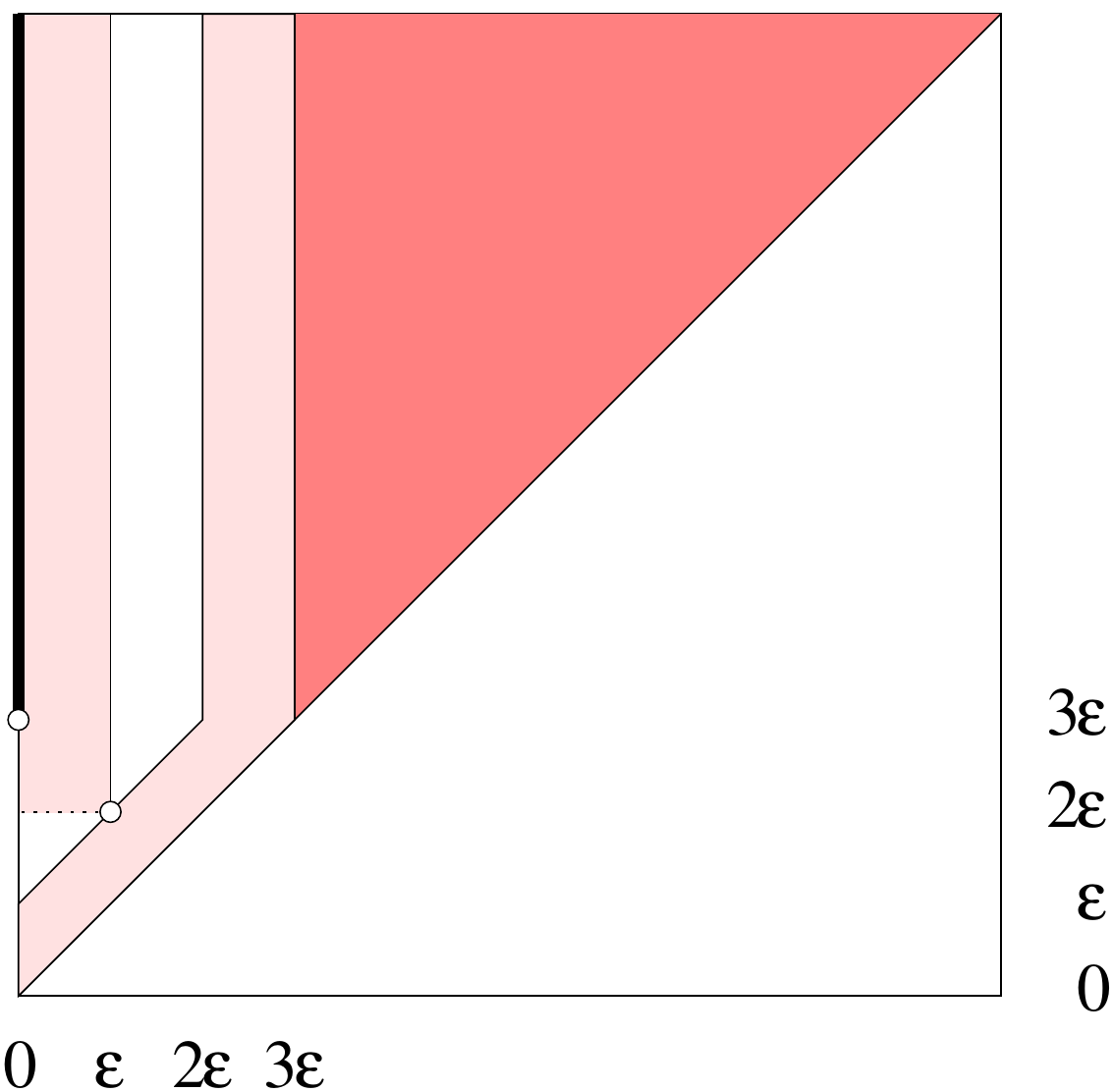}
 \end{center}
\caption[Regions in $\mathbb{X}$-diagrams and $U$-diagrams.]{The point in the $\mathbb{X}$-diagrams lie either along the solid black line or in the darkly shaded region.
Adding the lightly shaded regions, we get the region of possible points in the $U$-diagrams.}
\label{fig:IIT}
\end{figure}

\begin{theorem}[Topological Inference Theorem]
Suppose that we have an $\epsilon$-sample $\Uspace$ from $\Xspace$.
Then for each pair of points $p,q \in \Rspace^N$ such that $\rho = \rho(p,q,r) \geq 3 \epsilon$,
the map $\phi^{\Xspace} = \phi^{\Xspace}(p, q, r)$ is an isomorphism iff
\begin{align*}
\Ddgm{}{\kernel \phi^{\Uspace}}(\epsilon, 2 \epsilon) \cup \Ddgm{}{\cokernel \phi^{\Uspace}}(\epsilon,2\epsilon) = \emptyset.
\end{align*}
\label{result:IIT}
\end{theorem}
%
 \begin{proof}

To simplify exposition, we will refer to points in $\Ddgm{}{\kernel \phi^{\Xspace}} \cup \Ddgm{}{\cokernel \phi^{\Xspace}}$
and $\Ddgm{}{\kernel \phi^{\Uspace}} \cup \Ddgm{}{\cokernel \phi^{\Uspace}}$ as $\Xspace$-points and $\Uspace$-points, respectively.

Whenever $0 < \alpha < \beta < 3 \epsilon < \rho$, the two vertical maps in diagram (\ref{diag:K}) will by definition both be isomorphisms.
Hence the maps $\kernel \phi^{\Xspace}_{\alpha} \to \kernel \phi^{\Xspace}_{\beta}$
and $\cokernel \phi^{\Xspace}_{\alpha} \to \cokernel \phi^{\Xspace}_{\beta}$ must also be isomorphisms, and so, 
as $\alpha$ increases from $0$ to $\infty$,
any element of the (co)kernel of $\phi^{\Xspace}$ must live until at least $3 \epsilon$,
and any (co)kernel class which is born after $0$ must in fact be born after $3 \epsilon$.
In other words, any $\Xspace$-point
must lie either to the right of the line $x = 3 \epsilon$,
or along the $y$-axis and above the point $(0, 3 \epsilon)$; see Figure \ref{fig:IIT}.
Recall that $\phi^{\Xspace}$ is an isomorphism iff $\kernel \phi^{\Xspace} = 0 = \cokernel \phi^{\Xspace}.$
Thus $\phi^{\Xspace}$ is an isomorphism iff the black line in Figure \ref{fig:IIT} contains no $\Xspace$-points.

On the other hand, Theorem \ref{result:CKDS} requires that every $\Uspace$-point must lie within $\epsilon$ of an $\Xspace$-point.
That is, all $\Uspace$-points are contained
within the two lightly shaded regions drawn in Figure \ref{fig:IIT}.
Since the rightmost such region is more than $\epsilon$ away from the thick black line, there will be a $\Uspace$-point
in the left region iff there is an $\Xspace$-point on the thick black line.
But the $\Uspace$-points within the left region
are exactly the members of  $\Ddgm{}{\kernel \phi^{\Uspace}}(\epsilon, 2 \epsilon) \cup \Ddgm{}{\cokernel \phi^{\Uspace}}(\epsilon,2\epsilon)$. 
\eop
\end{proof}


\paragraph{Examples.}
Here we give two examples illustrating the topological inference theorem.

\begin{figure}[tbp]
 \begin{center}
  \includegraphics[scale=0.35]{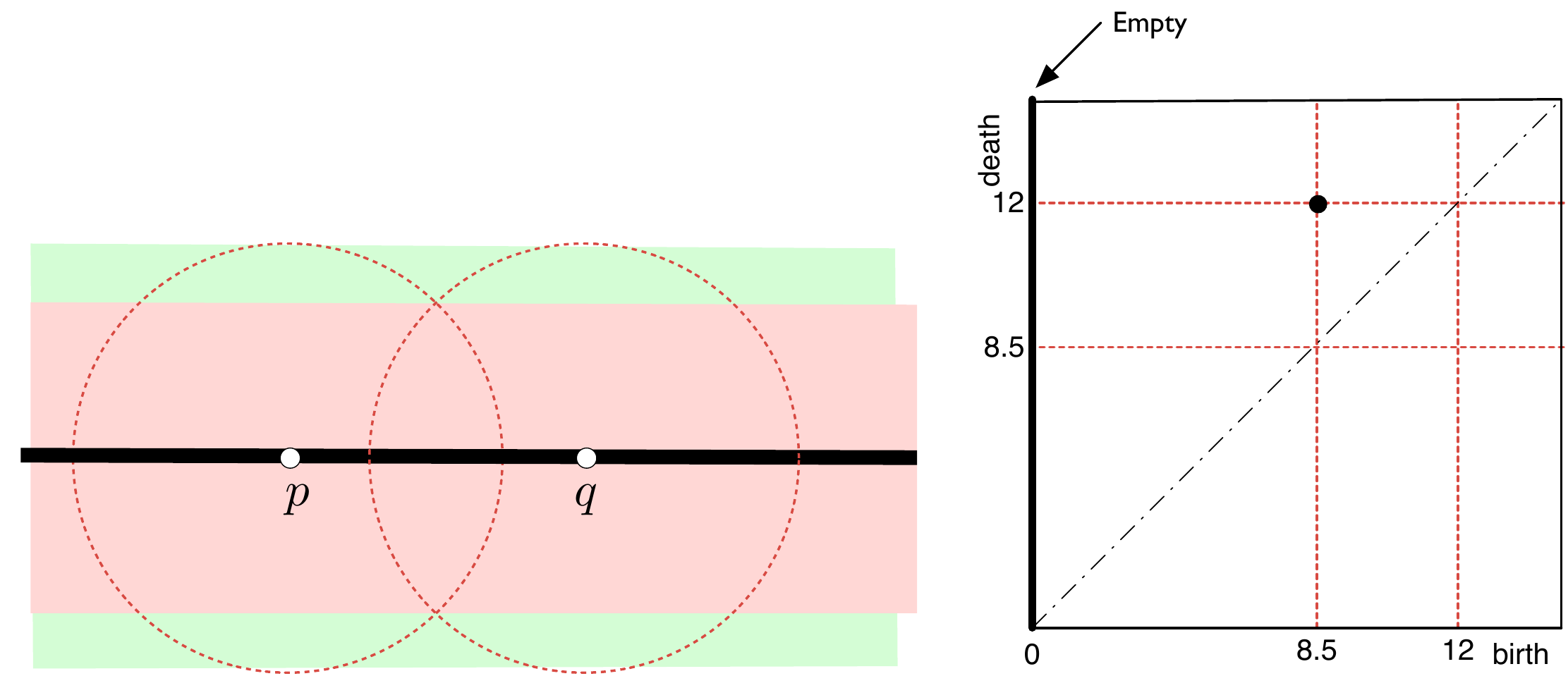}
 \end{center}
\caption[Kernel persistence diagram of two locally equivalent points, given $\mathbb{X}$.]{Kernel persistence diagram of two local equivalent points, given $\mathbb{X}$.}
\label{fig:tii-iso}
\end{figure}

For the first example, suppose we have the space $\Xspace$ in the left half of Figure \ref{fig:tii-iso}, and we take the labelled points $p$ and $q$ and the radius $r$ as drawn; 
in this case, one can show that $\rho(p,q,r) = 8.5$, which here is the distance between the line segment and the boundary of the intersection of the two $r$-balls.
First we compute the (co)kernel persistence diagrams for $\phi^{\Xspace}$, showing the kernel diagram in the right half of Figure \ref{fig:tii-iso}.
Since the $y$-axis of this diagram is free of any points (and the same holds for the un-drawn cokernel diagram), $p$ and $q$ have the
same local structure at this radius level.

\begin{figure}[tbp]
 \begin{center}
  \includegraphics[scale=0.35]{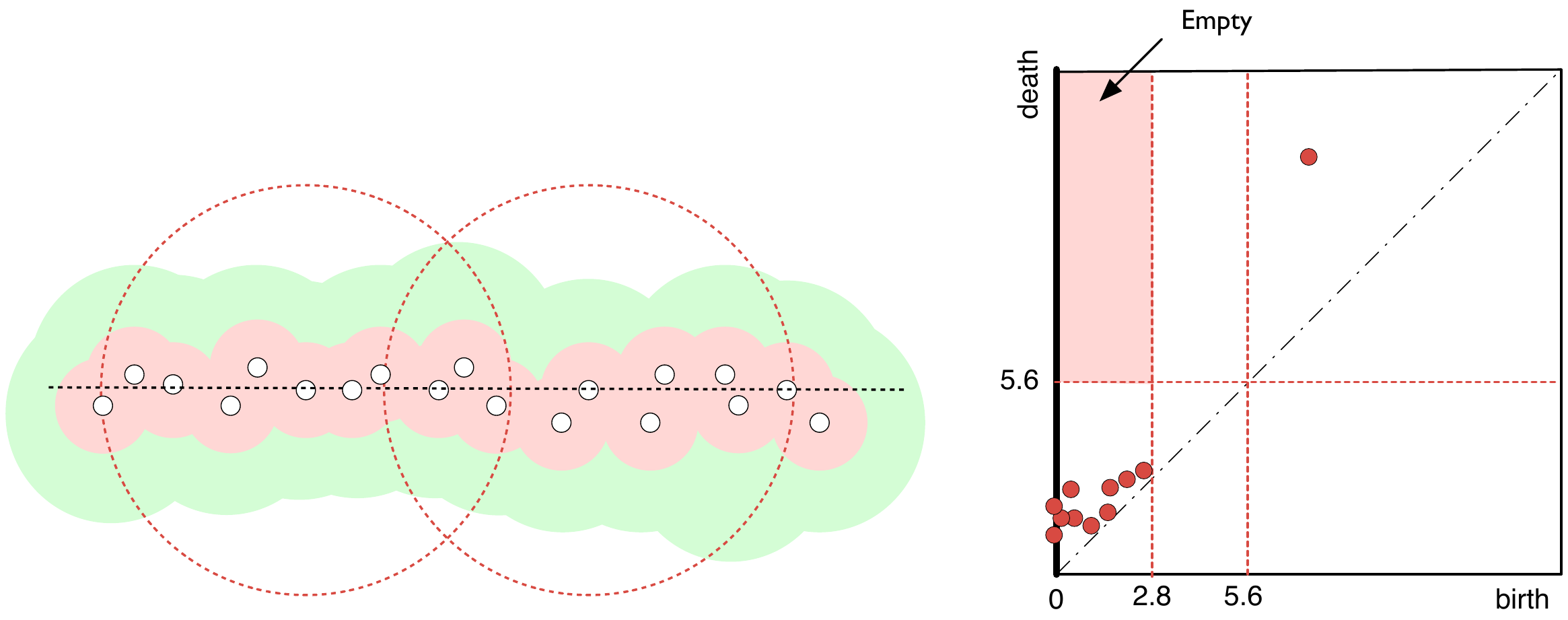}
 \end{center}
\caption[Kernel persistence diagram of two locally equivalent points, given $U$.]{Kernel persistence diagram of two local equivalent points, given $U$.}
\label{fig:tii-iso-pcd}
\end{figure}

On the other hand, suppose that we have an $\epsilon$-sample $U$ of $\Xspace$, with $\epsilon = 2.8 < \rho/3$, as drawn in the left half of Figure \ref{fig:tii-iso-pcd}.
We can compute the analogous $U$-diagrams, with the kernel diagram drawn in the right half of the same figure.
Noting that the two rectangles defined by $(\epsilon, 2\epsilon)$ in the two diagrams are indeed empty, and that the same holds for the cokernel diagrams,
we can apply Theorem \ref{result:IIT} to infer that the points have the same local structure at radius level $r$.

\begin{figure}[tbp]
 \begin{center}
  \includegraphics[scale=0.35]{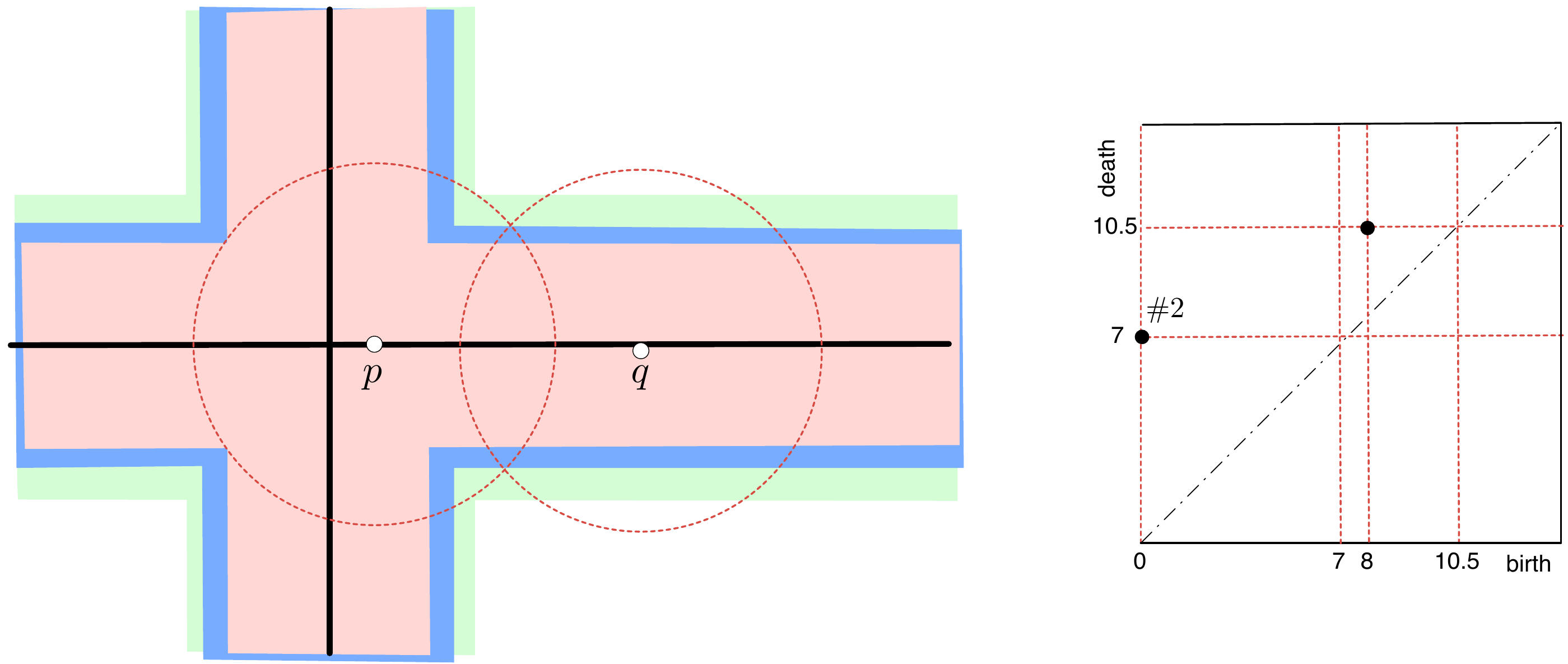}
 \end{center}
\caption[Kernel persistence diagram of two points that are not local equivalent, given $\mathbb{X}$.]
{Kernel persistence diagram of two points that are not locally equivalent, given $\mathbb{X}$. A number, i.e., $\#2$, labeling a point in the persistence diagram indicates its multiplicity.}
\label{fig:tii-noiso}
\end{figure}

\begin{figure}[tbp]
 \begin{center}
  \includegraphics[scale=0.35]{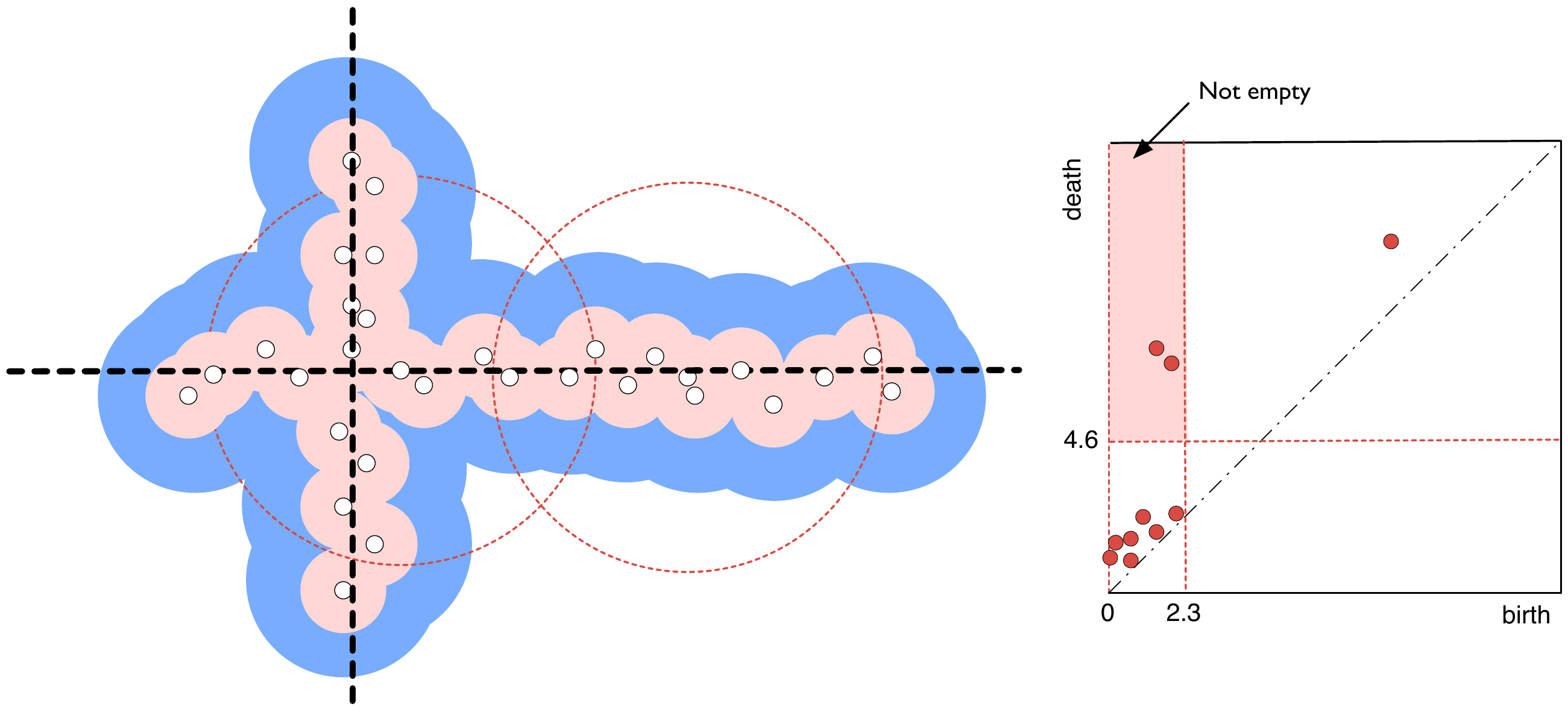}
 \end{center}
\caption[Kernel persistence diagram of two points that are not local equivalent, given $U$.]
{Kernel persistence diagram of two points that are not locally equivalent, given $U$.}
\label{fig:tii-noiso-pcd}
\end{figure}

For a second example, suppose $\Xspace$ is the cross on the left half of Figure \ref{fig:tii-noiso}, with $p,q,r$ as drawn.
Then $p$ and $q$ are locally different at this radius level, as shown by the presence of two points on the $y$-axis of the kernel
In the left half of Figure \ref{fig:tii-noiso-pcd}, we show an $\epsilon$-sample $U$ of $\Xspace$, with $3 \epsilon < \rho(p,q,r)$.
Note that the kernel diagram for $\phi^U$ does indeed have two points in the relevant rectangle.

\subsection{Geometric Intuition}

Theorem \ref{result:IIT} is stated in terms of a topological parameter, $\rho = \rho(p,q,r)$.
We can relate $\rho$ to more geometrically-flavored quantities. Specifically, we will
show that $\rho$ is lower bounded by a parameter derived from local variants of \emph{reach}, as well as from parameters related
to the gradient of $d_{\Xspace}$. Unfortunately, this lower bound can be quite loose, zero in certain cases, limiting its practical utility. It does provide a geometric intuition to the topological constraints
on the point cloud.

Recall that the medial axis $\Mcal$ of an embedded space $\Xspace$ is the subset of the ambient space consisting of all points which have
at least two nearest neighbors on $\Xspace$, and that the reach $\tau$ of $\Xspace$ is defined
by $\tau = \inf_{x \in \Xspace} \distone(x,\Mcal).$
We fix notation for the following four intersections of $\Mcal$ with different subsets of the ambient space:
$\Mcal(p,r) = \Mcal \cap B_r(p),
\Mcal_0(p,r) = \Mcal \cap \partial B_r(p),
 \Mcal(p,q,r) = \Mcal \cap B_r(p) \cap B_r(q),$
and  $\Mcal_0(p,q,r) = \Mcal \cap \partial(B_r(p) \cap B_r(q)),$
and we define a variant of reach for each such space:
\begin{align*}
& \tau(p, r) = \inf_{x \in \Xspace} \distone(x,\Mcal(p, r))\\
& \tau_0(p, r) =  \inf_{x \in \Xspace} \distone(x,\Mcal_0(p, r))\\
& \tau(p, q, r) =  \inf_{x \in \Xspace} \distone(x,\Mcal(p, q, r))\\
& \tau_0(p, q, r) =  \inf_{x \in \Xspace} \distone(x,\Mcal_0(p, q, r)).
\end{align*}
Note that all four of these quantities are of course upper bounds on $\tau$ itself.

Letting $\nabla_{\Xspace}$ be shorthand for the gradient of $d_{\Xspace}$, we define the following subset
of $\partial B_r(p):$
\begin{align*}
G(p,r) = \{y \in \partial B_r(p) \mid \nabla_{\Xspace}(y) \perp \partial B_r(p)\},
\end{align*}
and then set $\eta(p,r) = \inf_{x \in \Xspace} \distone(x,G(p, r)).$
We similarly define $G(p,q,r)$ and $\eta(p,q,r)$,
\begin{align*}
&G(p, q, r) = \{y \in \partial (B_r(p) \cap B_r(q)) \mid \nabla_{\Xspace}(y) \perp \partial (B_r(p) \cap B_r(q))\},\\
&\eta(p, q, r) = \inf_{x \in \Xspace} \distone (x, G(p, q, r)).
\end{align*}

Given the above quantities the following lower bound holds.
\begin{theorem}[Geometric lower bound]
\label{result:GIIT}
If we define
$$\gamma = \gamma(p,q,r) = \min\{\tau(p,r), \tau(p,q,r), \eta(p,r), \eta(p,q,r)\},$$
then
$\rho(p,q,r) \geq \gamma(p,q,r)$.
\end{theorem}
The proof appears in Appendix \ref{app:geoproof}.


\section{Probabilistic Inference Theorem}
\label{sec:PIT}

The topological inference of Section \ref{sec:TIT} states conditions under which the point sample $\Uspace$
can be used to infer stratification properties of the space $\Xspace$. The basic condition is that the Hausdorff
distance between the two must be small. In this section we describe two probabilistic models for generating the
point sample $\Uspace$, and we provide an estimate of how large this point sample should be to 
infer stratification properties of the space $\Xspace$ with a quantified measure of confidence.
More specifically, we provide a local estimate, based on $\rho(p,q,r)$ and $\rho(q,p,r)$, of how many sample points are needed to 
infer the local relationship at radius level $r$ between two fixed points $p$ and $q$; this same theorem can be used to give a global estimate
of the number of points needed for inference between any pair of points whose $\rho$-values are above some fixed low threshold.

\subsection{Sampling Strategies}

We assume $\Xspace$ to be compact. Since the stratified space  $\Xspace$ can contain singularities and
maximal strata of varying dimensions, some care is required in the sampling design.
Consider for example a sheet of area one, punctured by a line of length one. In this case, sampling from a naively  constructed uniform
measure on this space would result in no points being drawn from the line. This same issue arose and was
dealt with in \cite{NiySmaWei2008b}, although in a slightly different approach than we will develop.

The first sampling strategy is to remove the problems of singularities and varying dimension by
replacing $\Xspace$ by a slightly thickened version $\Xspace \equiv \Xspace_{\delta}$.
We assume that $\Xspace$ is embedded in $\Rspace^k$ for some $k$.
This new space is a smooth manifold with boundary and our point sample is a set of
 $n$ points drawn identically and independently from the uniform measure $\mu(\Xspace)$ on $\Xspace$,
 $\Uspace = \{x_1,...,.x_n\} \stackrel{iid}{\sim} \mu(\Xspace)$. This model can be thought of as placing
 an appropriate measure on the highest dimensional strata to ensure that lower dimensional strata will
 be sampled from. We call this model $M_1$.

The second sampling strategy is to deal with the problem of varying dimensions using a mixture
model. In the example of the sheet and line, a uniform measure would be placed on the sheet, while
another uniform measure would be placed on the line, and a mixture probability is placed
on the two measures; for example, each measure could be drawn with probability
$1/2$. We now formalize this approach. Consider each (non-empty) $i$-dimensional stratum
$\Sspace_i = \Xspace_i - \Xspace_{i-1}$ of $\Xspace$. All strata that are included
in the closure of some higher-dimensional strata, in other words all non-maximal strata, 
are not considered in the model.
A uniform measure is assigned to the closure of each maximal stratum,
$\mu_i(\Sspace_i)$, this is possible since each such closure is compact.  We assume a finite number of
maximal strata $K$ and assign to the closure of each such stratum a probability $p_i = 1/K$. This implies the following
density
$$f(x) = \frac{1}{K} \sum_{j=1}^K  \nu_i(X=x),$$
where $\nu_i$ is the density corresponding to measure $\mu_i$.
The point sample is generated from the following model:
$\Uspace = \{x_1,...,.x_n\} \stackrel{iid}{\sim} f(x)$. We call this model $M_2$.

The first model replaces a stratified space with its thickened version, which enables us to place
a uniform measure on the thickened space. Although this replacement makes it convenient for sampling,
it does not sample directly from the actual space. The second model samples from the actual space,
however the sample is not uniform on $\Xspace$ with respect to Lebesgue measure.

\subsection{Lower bounds on the sample size of the point cloud}

Our first main theorem is the probabilistic analogue of Theorem \ref{result:IIT}.
An immediate consequence of this theorem is that, for two points $p,q \in \Uspace$,
we can infer  with probability at least $1-\prb$  whether $p$ and $q$ are locally equivalent,
$p \sim_r q$. The confidence level  $1-\prb$ will be a monotonic function of the size of the point sample.

The theorem involves a parameter $v(\rho),$ for each positive $\rho$, which is based
on the volume of the intersection of $\rho$-balls with $\Xspace$.
First we note that each maximal stratum of $\Xspace$ comes with its own notion of volume: in the plane
punctured by a line example, we measure volume in the plane and in the line as area and length, respectively.
The volume $\vol{\Yspace}$ of any subspace $\Yspace$ of $\Xspace$ is the sum
of the volumes of the intersections of $\Yspace$ with each maximal stratum.
For $\rho > 0,$ we define 
\begin{equation}
 v(\rho) = \inf_{x \in \Xspace} \frac{\vol{B_{\rho/24}(x) \cap \Xspace}}{\vol{\Xspace}}
\label{eqn:volume}
\end{equation}
We can then state:
\begin{theorem}[Local Probabilistic Sampling Theorem]
\label{locprob}
Let $U = \{x_1, x_2, ..., x_n\}$ be drawn from either model $M_1$ or $M_2$.
Fix a pair of points $p, q \in \Rspace^N$ and a positive radius $r$, and put
$\rho = \min \{\rho(p, q, r), \rho(q, p, r)\}$.
If
\begin{align*}
n \geq \frac{1}{v(\rho)}\left(\log \frac{1}{v(\rho)} + \log \frac{1}{\prb}\right),
\end{align*}
then, with probability greater than $1-\prb$ we can correctly infer whether or not $\phi^{\Xspace}(p,q,r)$ and $\phi^{\Xspace}(q,p,r)$ are both isomorphisms.
\end{theorem}

\begin{proof}

A finite collection $U = \{x_1, x_2, ..., x_n\}$ of points in $\Rspace^N$ is \emph{$\ep$-dense}
with respect to $\Xspace$ if $\Xspace \subseteq U^{\ep}$; equivalently, $U$ is an $\ep$-cover of
$\Xspace$. Let $C(\ep)$ be the \emph{$\ep$-covering number} of $\Xspace$, the minimum number of sets
$B_{\ep} \cap \Xspace$ that cover $\Xspace$. Let $P(\ep)$ be the \emph{$\ep$-packing number} of
$\Xspace$, the maximum number of sets $B_{\ep} \cap \Xspace$ that can be packed into
$\Xspace$ without overlap.

We consider a cover of $\Xspace$ with balls of radius $\rho/12$.
If there is a sample point in each $\rho/12$-ball, then $U$ will be an $\ep$-approximation of $\Xspace$, with $\ep \leq 4 (\rho/12) = \rho/3$.
This satisfies the condition of the topological inference theorem, and therefore we can infer the local structure
between $p$ and $q$.

The following two results from \cite{NiySmaWei2008} will be useful in computing the number of
sample points $n$ needed to obtain, with confidence, such an $\ep$-approximation.

\begin{lemma}[Lemma 5.1 in \cite{NiySmaWei2008}]
\label{lemma:estsample}
Let $\{A_1, A_2, ..., A_l\}$ be a finite collection of measurable sets with probability measure $\mu$ on $\cup^l_{i=1} A_i$, such that for all $A_i$, $\mu(A_i) > \alpha$. Let $\Uspace = \{x_1, x_2, ..., x_n\}$ be drawn iid according to $\mu$. If
$n \geq \frac{1}{\alpha} (\log l + \log \frac{1}{\prb})$, then, with probability $1 - \prb$, $\forall i$, $\Uspace \cap A_i \neq \emptyset$.
\end{lemma}

\begin{lemma}[Lemma 5.2 in \cite{NiySmaWei2008}]
\label{lemma:boundcovernumber}
Let $C(\ep)$ be the covering number of an $\ep$-cover of $\Xspace$ and $P(\ep)$ be the packing number of an $\ep$-packing, then
$$P(2\ep) \leq C(2\ep) \leq P(\ep).$$
\end{lemma}

Again, we consider a cover of $\Xspace$ by balls of radius $\rho/12$.
Let $\{y_i\}_{i=1}^{l} \in \Xspace$ be the centers of the balls contained in a minimal sub-cover.
Put $A_i = B_{\rho/12}(y_i) \cap \Xspace$.
Applying Lemma \ref{lemma:estsample}, we obtain the estimate
\begin{align*}
n \geq \frac{1}{\alpha} \left(\log l + \log \frac{1}{\prb}\right),
\end{align*}
where $l$ is the $\rho/12$-covering number, and $\alpha = \min_{i} \frac{\vol {A_i}}{\vol \Xspace} $.

Applying Lemma \ref{lemma:boundcovernumber},
\begin{align*}
l = C(\rho/12) \leq P(\rho/24) \leq \frac{\vol{\Xspace}}{\vol{B_{\rho/24} \cap \Xspace}} \leq \frac{1}{v(\rho)}.
\end{align*}
On the other hand, $\frac{1}{\alpha} \leq \frac{1}{v(\rho)}$ by definition, and the result follows. \eop

\end{proof}

To extend the above theorem to a more global result, one can pick a 
positive $\rho$ and radius $r$, and consider the set of all pairs of points 
$(p,q)$ such that $\rho \leq \min\{\rho(p,q,r),\rho(q,p,r\}$. Applying 
Theorem \ref{locprob} uniformly to all pairs of points will give the minimum number of sample
points needed to settle the isomorphism question for all of the intersection maps between all pairs.


\section{Algorithm}
\label{sec:Alg}

The theorems in the last sections give conditions under which a point cloud $U$, sampled from a stratified space $\Xspace$, 
can be used to infer the local equivalences between points on $\Xspace$ and its surrounding ambient space.
We now switch gears slightly, and imagine clustering the $U$-points themselves into strata.
The basic strategy is to build a graph on the point set, with edges coming from positive isomorphism judgements.
The connected components of this graph will then be our proposed strata.
We begin by describing this strategy in Section \ref{subsec:Clust}.
Some of its potential limitations are discussed in Section \ref{subsec:LD}, where we also describe a more robust variant
based on graph diffusion.

A crucial subroutine in the clustering algorithm is the computation of the diagrams 
$\Ddgm{}{\kernel \phi^U}$ and $\Ddgm{}{\cokernel \phi^U}$, for $\phi^U = \phi^U(p,q,r)$ between
all pairs $(p,q) \in U \times U$. The algorithm for this sub-routine is quite complicated, we describe
it in detail in Section \ref{subsec:DC}.
The correctness proof is even more complicated; we give a proof sketch in Section \ref{subsec:correctness}, deferring
all major details to Appendix \ref{app:correctness}.

We would like to make clear that we consider the algorithm in this paper a first step and several
issues both statistical and computational can be improved upon.

\subsection{Clustering}
\label{subsec:Clust}

We imagine that we are given the following input: a point cloud $U$ sampled from some stratified space $\Xspace$, and a fixed radius $r$.
We make the assumption that $d_H(U,\Xspace) \leq \epsilon \leq \frac{\rho_{\min}}{3}$, where $\rho_{\min}$ is the minimum of $\rho(p,q,r)$ for all
pairs $(p,q) \in U \times U$.
Later we discuss the consequences when this assumption does not hold and a possible solution.

We build a graph where each node
in the graph corresponds uniquely to a point from $U$. 
Two points $p, q \in U$ (where $||p - q|| \leq 2r$) are connected by an edge
iff both $\phi^{\Xspace}(p,q,r)$ and $\phi^{\Xspace}(q, p, r)$ are isomorphisms,
equivalently iff $\Ddgm{}{\kernel \phi^{\Uspace}}(\epsilon, 2 \epsilon)$ and $\Ddgm{}{\cokernel \phi^{\Uspace}}(\epsilon, 2 \epsilon)$ are empty.
The connected components of the resulting graph are our clusters. 
A more detailed statement of this procedure is giving in pseudo-code, see Algorithm \ref{alg:global}. 
Note that the connectivity of the graph is encoded by a weight matrix, and our clustering strategy is based on a $0/1$-weight assignment. 

\begin{algorithm}
\caption{Strata-Inference({$\Uspace, r, \ep$})}
\begin{algorithmic}
	\FORALL {$p, q \in \Uspace$}
     	\IF {$||p - q|| > 2r$}
      	\STATE  $W(p, q) = 0$
     \ELSE
     		\STATE Compute $\dgmD{ \ker \phi^{U}(p, q, r)}$ and $\dgmD{ \cok \phi^U(p, q, r)}$
		\STATE Compute $\dgmD{ \ker \phi^{U}(q, p, r)}$ and $\dgmD{ \cok \phi^U(q, p, r)}$
		\IF {$\Ddgm{}{\kernel \phi^{\Uspace}(p, q, r)}(\epsilon, 2 \epsilon) \cup \Ddgm{}{\cokernel \phi^{\Uspace}(p, q, r)}(\epsilon,2\epsilon) \neq \emptyset$}
			\STATE $W(p,q) = 0$
		\ELSIF	 {$\Ddgm{}{\kernel \phi^{\Uspace}(q, p, r)}(\epsilon, 2 \epsilon) \cup \Ddgm{}{\cokernel \phi^{\Uspace}(q, p, r)}(\epsilon,2\epsilon) \neq \emptyset$}
			\STATE $W(p,q) = 0$
		\ELSE
			\STATE $W(p, q) = 1$
		\ENDIF
	\ENDIF
	\ENDFOR
	\STATE Compute connected components based on $W$.
\end{algorithmic}
\label{alg:global}
\end{algorithm}

\subsection{Robustness of clustering}
\label{subsec:LD}

Two types of errors in the clustering can occur: false positives where the algorithm connects
points that should not be connected and false negatives where points that should be connected
are not. The current algorithm we state is somewhat brittle with respect to both false positives
as well as false negatives. We will suggest a very simple adaptation of our current algorithm
that should be more stable with respect to both false positives and false negatives.

The false positives are driven by the condition in Theorem \ref{result:IIT} that $\rho_{\min} < 3 \epsilon$,
so if the point cloud is not sampled fine enough we can get incorrect positive isomorphisms
and therefore incorrect edges in the graph. If we use transitive closure to define the connected
components this can be very damaging in practice since a false edge can collapse disjoint components
into one large cluster. 

The false negatives occur because our point sample $U$ is not fine enough to capture
chains of points that connect pairs in $U$ through isomorphisms, there may be other points
in $\Xspace$ which if we had sampled then the chain would be observed. The probability
of these events in theory decays exponentially as the sample size increases and the
confidence parameter $\xi$ in Theorem \ref{result:IIT}  controls these errors.

We now state a simple adaptation of the algorithm that will make it more robust. It is natural
to think of the $0/1$-weight assignment on pairs of points $p,q \in \Uspace$ as an association
matrix $\mathbf{W}$. A classic approach for robust partitioning is via spectral graph theory
\cite{MeilaShi,KannanVempala,Chung}. This approach is based an eigen-decomposition of the the graph Laplacian, 
$\mathbf{L} = \mathbf{D} -\mathbf{W}$ with the diagonal matrix $ \mathbf{D}_{ii} = \sum_j \mathbf{W}_{ij}.$ 
The smallest nontrivial eigenvalue $\lambda_1$ of $\mathbf{W}$ is called the Fiedler constant
and estimates of how well the vertex set can be partitioned \cite{Fiedler}. The corresponding
eigenvector $v_1$ is used to partition the vertex set. There are strong connections between
spectral clustering and diffusions or random walks on graphs \cite{Chung}. 

The problems of spectral clustering and lower dimensional embeddings
have been examined from a manifold learning perspective \cite{BelNiy2002,BelNiy2008,GineKol}.
The idea central to these analyses is given a point sample from a manifold construct an appropriate graph
Laplacian and use its eigenvectors to embed the point cloud in a lower dimensional space.
A theoretical analysis of this idea involves proving convergence of the graph Laplacian
to the Laplace-Beltrami operator on the manifold and the convergence of the eigenvectors
of the graph Laplacian to the eigenvalues of the  Laplace-Beltrami operator. A key quantity
in this analysis is the Cheeger constant which is the first nontrivial eigenvalue of the Laplace-Beltrami 
operator \cite{Cheeger}. An intriguing question is whether the association matrix we construct from
the point cloud can be related to the Laplacian on high forms.

\subsection{Diagram Computation}
\label{subsec:DC}

We now describe the computation of the diagrams $\Ddgm{}{\kernel \phi^U}$ and $\Ddgm{}{\cokernel \phi^U}$.
To do this, we need for each $\alpha \geq 0$ a simplicial analogue of the map
\begin{align*}
\phi^{U}_{\alpha}: \Hgroup(B_{p}^{U}(\alpha), \bdr B_{p}^{U}(\alpha)) \to \Hgroup(B_{pq}^{U}(\alpha), \bdr B_{pq}^{U}(\alpha)).
\end{align*}
To produce this, we first define, for each $\alpha \geq 0$, two pairs of simplicial complexes $L_0(\alpha) \subseteq L(\alpha)$
and $K_0(\alpha) \subseteq K(\alpha)$, and z relative homology map
\begin{align*}
\psi_{\alpha}: \Hgroup(L(\alpha),L_0(\alpha)) \to \Hgroup(K(\alpha),K_0(\alpha)
\end{align*}
between them.
We will then show that 
$$\Ddgm{}{\kernel \phi^U} = \Ddgm{}{\kernel \psi} \mbox{ and } \Ddgm{}{\cokernel \phi^U} = \Ddgm{}{\cokernel \psi}.$$
To compute the diagrams involving $\psi$, we reduce various boundary matrices;  since we follow very closely the (co)kernel persistence 
algorithm described in \cite{CohEdeHar2009b}, we omit any further details here.



\subsubsection{Complexes}

To construct the simplicial complexes in our algorithm, we take the nerves of several collections of sets which are derived from a variety of Voronoi
diagrams of different spaces. Here we briefly review these concepts.

\paragraph{Nerves.}
The \emph{nerve} $N(\Ccal)$ of a finite collection of sets $\Ccal$ is defined to be the abstract simplicial complex
with vertices corresponding to the sets in $\Ccal$ and 
with simplices corresponding to all non-empty intersections among these sets;
that is, $N(\Ccal) = \{S \subseteq \Ccal \mid \bigcap S \neq \emptyset\}.$
Every abstract simplicial complex can be geometrically realized, and therefore the concept of homotopy type makes sense.
Under certain conditions, for example whenever the sets in $\Ccal$ are all closed and convex
subsets of Euclidean space (\cite{EdeHar2010}, p.59), the nerve of $\Ccal$ has the same
homotopy type, and thus the same homology groups, as the union of sets in $\Ccal$.

\paragraph{Voronoi diagram.}
If $U$ is a finite collection of points in $\Rspace^k$ and $u_i \in U$, then the \emph{Voronoi cell} of $u_i$
is defined to be:
$$V_i = V(u_i) = \{x \in \Rspace^k \mid ||x - u_i|| \leq ||x - u_j||, \forall u_j \in U\}.$$
The set of cells $V_i$ covers the entire space and forms the \emph{Voronoi diagram} of $\Rspace^k$, 
denoted as $\voi{U | \Rspace^k}$.
If we restrict each $V_i$ restricted to some subset $X \subseteq \Rspace^k$, then the
set of cells $V_i \cap X$ forms a \emph{restricted Voronoi diagram}, denoted as $\voi{U | X}$.  
For a simplex $\sigma \in U$,  we set $V_{\sigma} = \cap_{u_i \in \sigma} V_i$.

The nerve of the restricted Voronoi diagram $\voi{U | X}$ is called the \emph{restricted Delaunay triangulation},
denoted as $\del{U | X}$. It contains the set of simiplices $\sigma$ for which $V_{\sigma} \cap X \neq \emptyset$.

\paragraph{Power cells.}

An important task in our algorithm is the computation of the relative homology
groups $\Hgroup(B_p^U(\alpha),\bdr B_p^U(\alpha))$ and $\Hgroup(B_{pq}^U(\alpha),\bdr B_{pq}^U(\alpha))$.
Now to compute $\Hgroup(U_{\alpha})$, the absolute homology of the thickened point cloud,
we would need only to compute the nerve of the collection of sets $V_i \cap U_{\alpha}$.
This is because each such set is convex and their union obviously equals the space $U_{\alpha}$.
Such a direct construction will not work in our context, for the simple reason that
the sets $V_i \cap \bdr B_p^U(\alpha)$ and $V_i \cap \bdr B_{pq}^U(\alpha)$ need not be convex.

\begin{figure}[tbp]
 \begin{center}
  \includegraphics[scale=0.25]{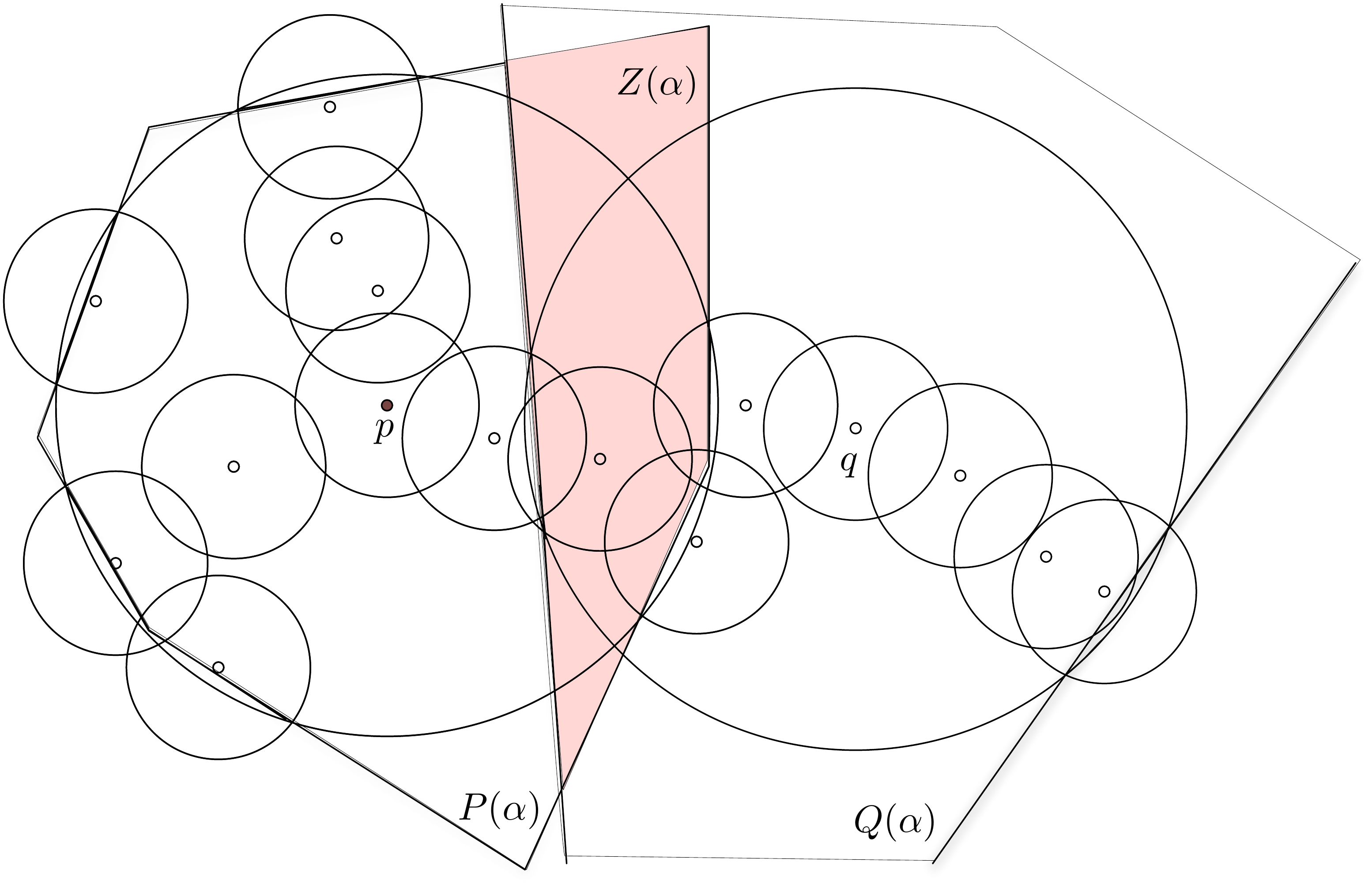}
 \end{center}
\caption[Illustration of $Z(\alpha)$.]
{Illustration of intersection power cell $Z(\alpha)$, as the shaded region. The unshaded convex regions are $P(\alpha)$ and $Q(\alpha)$ respectively.}
\label{fig:ipc}
\end{figure}

To get around this problem, we first define $P(\alpha)$, the \emph{power cell} with
respect to $B_r(p)$, to be:
\begin{equation}
 P(\alpha) = \{x \in \Rspace^k \mid ||x-p||^2 - r^2 \leq ||x - u||^2 - \alpha^2, \forall u \in U\},
\label{eq:power}
\end{equation}
and we set $P_0(\alpha) = B_r(p) - \interior P(\alpha)$.
To define $Q(\alpha)$, the power cell with respect
to $B_r(q)$, we replace $p$ with $q$ in (\ref{eq:power}).
Finally, we set $Z(\alpha) = P(\alpha) \cap Q(\alpha)$, and $Z_0(\alpha) = (B_r(p) \cap B_r(q)) - \interior Z(\alpha)$.
This is illustrated in Figure \ref{fig:ipc}. 
We note that $P_0(\alpha)$ and $Z_0(\alpha)$ are both contained in  $U_{\alpha}$, as can be seen by manipulating the inequalities
in their definitions.

Replacing $\bdr B_p^U(\alpha)$ with $P_0(\alpha)$ and $\bdr B_{pq}^U(\alpha)$ with $Z_0(\alpha)$
has no effect on the relative homology groups in question, as is implied by the following two lemmas. The first
lemma was proven in \cite{BenCohEde2007}; a proof of the second appears in Appendix \ref{app:correctness}.
\begin{lemma}[Power Cell Lemma] 
Assume $B_r(p) - P_0(\alpha) \neq 0$.
The identity on $B_p^U(\alpha)$ is a homotopy equivalence of 
$(B_p^U(\alpha),\bdr B_p^U(\alpha))$ and
$(B_p^U(\alpha), P_0(\alpha))$, as a map of pairs. 
\label{result:PCL}
\end{lemma}
\begin{lemma}[Intersection Power Cell Lemma] 
Assume $B_r(p) \cap B_r(q) - Z_0(\alpha) \neq 0$.
The identity on $B_{pq}^U(\alpha)$ is a homotopy equivalence of 
$(B_{pq}^U(\alpha),\bdr B_{pq}^U(\alpha))$ and
$(B_{pq}^U(\alpha), Z_0(\alpha))$, as a map of pairs. 
\end{lemma}

\paragraph{Lune and moon.}
It can be shown (\cite{BenCohEde2007}) that the sets $V_i \cap P_0(\alpha)$ are convex.
Unfortunately, it is still possible for some set $V_i \cap Z_0(\alpha)$ to be non-convex.
To fix this, we must further divide the Voronoi cells in a manner we now describe.

\begin{figure}[tbp]
\begin{center}
\includegraphics[scale=0.35]{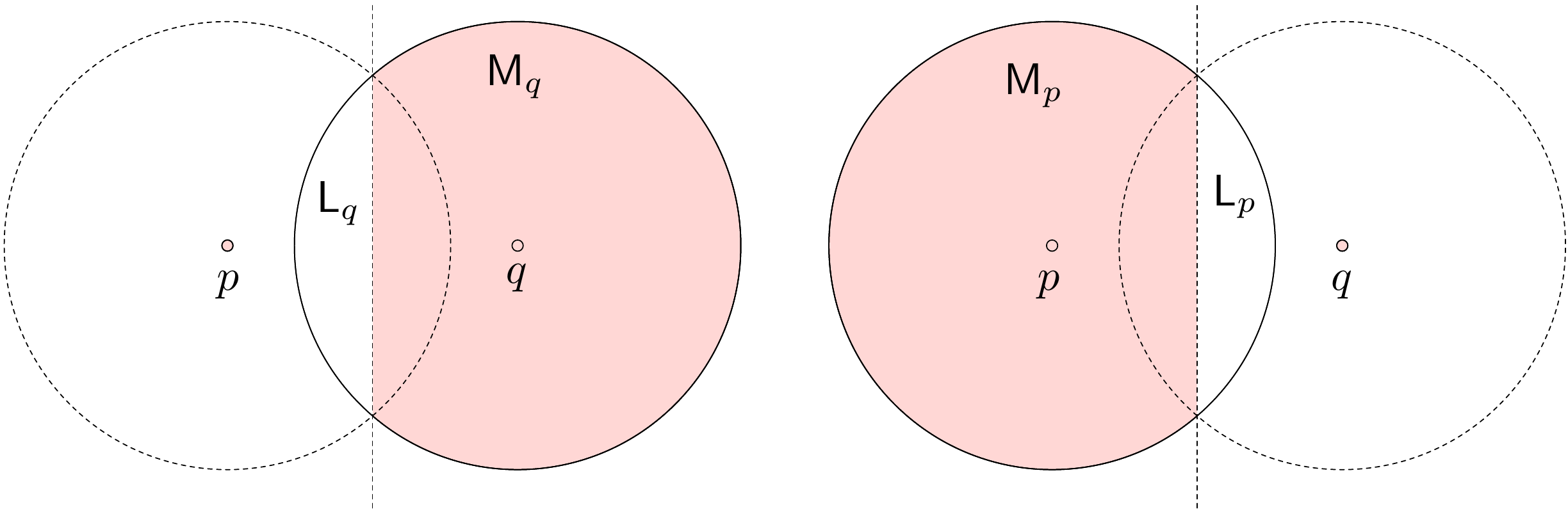} 
\end{center}
\caption[Illustration of the lune and the moon.]{Illustration of the lune and the moon. The shaded regions are the respective moons. The white regions within solid circles are the respective lunes.}
\label{fig:strata-lune}
\end{figure}

We consider the hyperplane $\Bisector$ of points in $\Rspace^k$ which are equidistant
to $p$ and $q$; we often refer to this hyperplane as the \emph{bisector}.
This will divide $\Rspace^k$ into two half-spaces; let $\Bisector_p$
and $\Bisector_q$  
denote the half-spaces containing $p$ and $q$, respectively..
We also define the \emph{$p$-lune}, $\Lune_p = \Bisector_q \cap B_r(p)$,
and the \emph{$p$-moon}, $\Moon_p =  \Bisector_p \cap B_r(p)$, as illustrated in 
Figure \ref{fig:strata-lune}.
The lune and the moon divide each Voronoi cell into two parts,
$V_i^{\Lune} = V_i \cap \Lune_p$ and $V_i^{\Moon} = V_i \cap \Moon_p$.
These sets are convex, assuming they are non-empty, since they are each the intersection
of two convex sets.
Furthermore, we have the following lemma whose simple but technical proof we omit:
\begin{lemma}[Convexity Lemma]
The sets $V_i^{\Lune} \cap Z_0(\alpha)$ and $V_i^{\Moon} \cap Z_0(\alpha)$ are
all convex, assuming they are non-empty.
\label{result:ConvL} 
\end{lemma}
Of course the nonempty sets among $V_i^{\Lune} \cap P_0(\alpha)$ and $V_i^{\Moon} \cap P_0(\alpha)$
will also be convex.

To construct the simplicial complexes needed in our algorithm,
we define $\Acal$ to be the collection of the nonempty sets
among $V_i^{\Lune} \cap B_p^U(\alpha)$ and $V_i^{\Moon} \cap B_p^U(\alpha)$,
and we define $\Acal_0$ to be the collection of the nonempty sets among
$V_i^{\Lune} \cap P_0(\alpha)$ and $V_i^{\Moon} \cap P_0(\alpha)$.
Note that $\cup \Acal = B_p^U(\alpha)$ and $\cup \Acal_0 = P_0(\alpha)$.
Taking the nerve of both collections, we define the simplicial
complexes $L(\alpha) = N(\Acal)$ and $L_0(\alpha) = N(\Acal_0)$.

Similarly, we define $\Ccal$ and $\Ccal_0$ to be the collections of the nonempty sets among,
respectively, $V_i^{\Lune} \cap B_{pq}^U(\alpha)$ and $V_i^{\Moon} \cap B_{pq}^U(\alpha)$,
and $V_i^{\Lune} \cap Z_0(\alpha)$ and $V_i^{\Moon} \cap Z_0(\alpha)$.
And we define $K(\alpha) = N(\Ccal)$ and $K_0(\alpha) = N(\Ccal_0)$.

\begin{figure}[tbp]
\begin{center}
\includegraphics[scale=0.35]{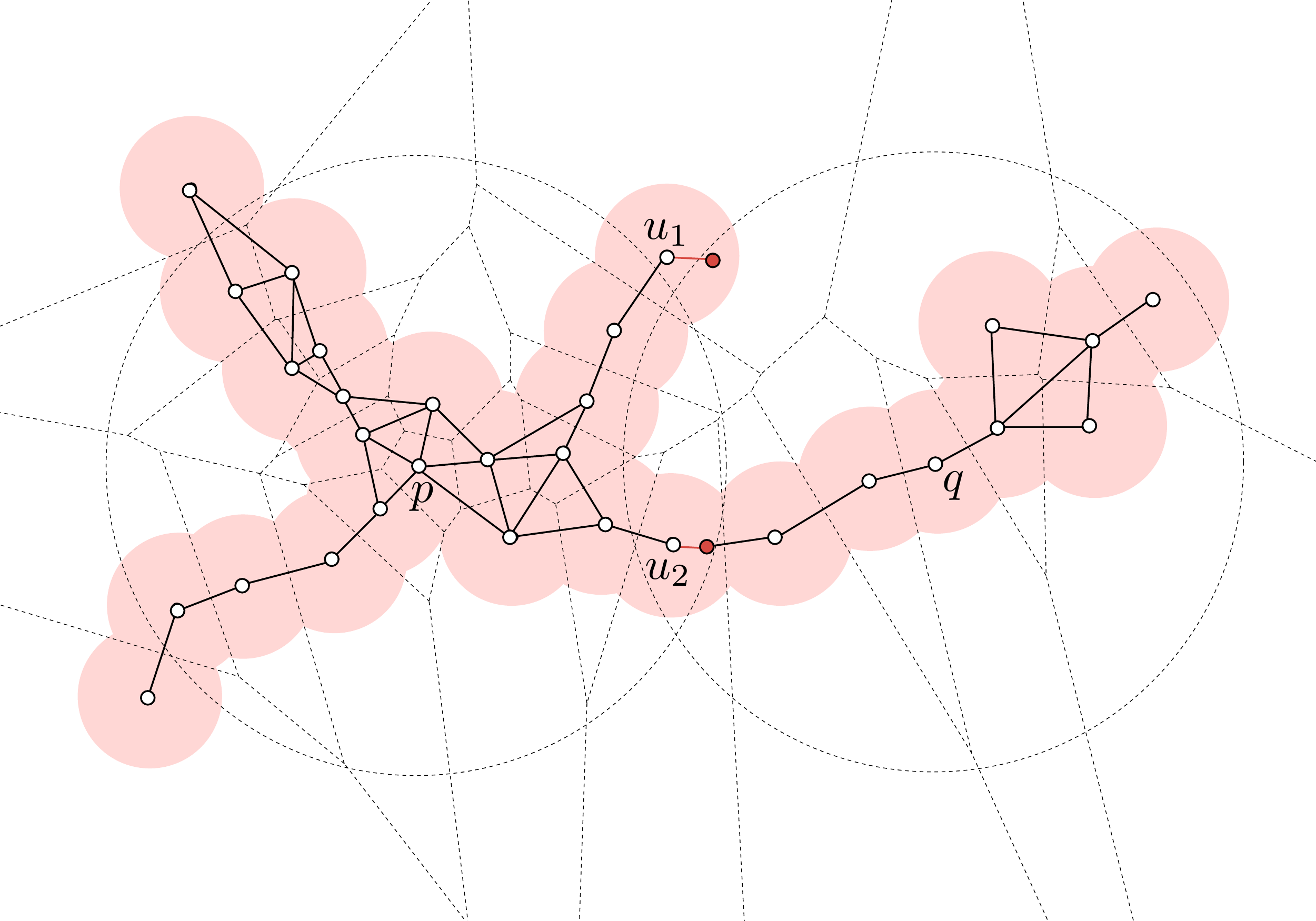} 
\end{center}
\caption[Illustration of the simplicial complexes constructed.]
{Illustration of the simplicial complexes constructed around two points $p$ and $q$. 
The underlying Voronoi decomposition of the space is shown in thin dotted lines.
$u_1$ and $u_2$ in $U$ are the points whose restricted Voronoi regions intersect with the lune at non-convex regions.}
\label{fig:vd2}
\end{figure}

An example of the simplicial complexes constructed in $\Rspace^2$ for a given $U$ are illustrated in Figure \ref{fig:vd2}.
A direct approach to construct these simplicial complexes runs into difficulties as the corners of the convex sets created by the bisector
can be shared by many sets; we defer the technicalities to Appendix \ref{app:purburbation}.

\subsubsection{Maps}

We now construct simplicial analogues
\begin{align*}
\psi_{\alpha}: \Hgroup(L(\alpha),L_0(\alpha)) \to \Hgroup(K(\alpha),K_0(\alpha).
\end{align*}
of the maps $\phi_{\alpha}^U$.

The containments $L_0(\alpha) \subseteq L(\alpha)$ and $K_0(\alpha) \subseteq K(\alpha)$
are obvious.
In order to define $\psi_{\alpha}$, we first need the following technical lemma:
\begin{lemma}[Containment Lemma]
 Assume that a simplex $\sigma$ is in $L_0(\alpha)$. If $\sigma$ is also in $K(\alpha)$, 
then $\sigma$ is in $K_0(\alpha),$ as well.
\label{result:ContL}
\end{lemma}
\begin{proof}.
By definition, $\sigma \in L_0(\alpha)$ iff there exists some point $x \in V^{\sigma} \cap P_0(\alpha).$
We must show that the set $V^{\sigma} \cap Z_0(\alpha)$ is non-empty.
Note that $x \in P_0(\alpha)$ implies that $x \in B_r(p)$, while $x \not \in \interior P(\alpha)$ implies
that $x \not \in \interior Z(\alpha)$.
If $x \in B_r(q)$, then we are done, since $Z_0(\alpha) = B_r(p) \cap B_r(q) - \interior Z(\alpha)$.

Otherwise, choose some point $y \in V^{\sigma} \cap U_{\alpha} \cap B_r(p) \cap B_r(q)$, which is possible
since $\sigma \in K(\alpha)$.
Since both $x$ and $y$ belong to the same convex set $V^{\sigma} \cap U_{\alpha} \cap B_r(p)$, there exists a directed line segment $\gamma$ from $x$ to $y$ within this set connecting them. We imagine moving along $\gamma$ and first we suppose that $\gamma$ intersects $B_r(q)$ before it intersects $\interior Q(\alpha)$. 
Let $z$ be the first point of intersection. 
Then $z \in B_r(p) \cap B_r(q)$, $z \notin \interior Q(\alpha)$. Therefore $z \in V^{\sigma} \cap Z_0(\alpha)$.
On the other hand, we may prove by contradiction that it is impossible for $\gamma$ to intersect $Q(\alpha)$ before it intersects $B_r(q)$.  
Let $z'$ be the first point of such an intersection.
Since $z' \in Q(\alpha)$, by definition $||z' - q||^2 - r^2 \leq ||z' - u_i||^2 - \alpha^2$, $\forall u_i \in U$.
Since $z' \in U_{\alpha}$, $\forall u_i \in \sigma$, $||z' - u_i||^2 - \alpha^2 \leq 0$.
Therefore $||z' - q||^2 - r^2 \leq ||z' - u_i||^2 - \alpha^2 \leq 0$, $\forall u_i \in \sigma$.
Since $z'$ is outside $B_r(q)$, $||z' - q||^2 - r^2 > 0$. This is a contradiction. \eop
\end{proof}

To define $\psi_{\alpha}$, we first construct a chain map $g = g_{\alpha}: \Cgroup(L(\alpha)) \to \Cgroup(K(\alpha))$
as follows.
Given a simplex $\sigma \in L(\alpha)$, we define $g(\sigma) = \sigma$ if $\sigma \in K(\alpha)$,
and $g(\sigma) = 0$ otherwise; we then extend $g$ to a chain map by linearity.
Using the Containment Lemma, we see that $g(\Cgroup(L_0(\alpha))) \subseteq \Cgroup(K_0(\alpha))$,
and thus $g$ descends to a relative chain map
$f = f_{\alpha}: \Cgroup(L(\alpha),L_0(\alpha)) \to \Cgroup(K(\alpha),K_0(\alpha))$.
Since $f$ clearly commutes with all boundary operators, it induces
a map on relative homology, and this is our $\psi = \psi_{\alpha}$.


\subsection{Correctness}
\label{subsec:correctness}

We show that our algorithm is correct by proving the following theorem. A sketch of the proof is given here, with the details deferred to Appendix \ref{app:correctness}.

\begin{theorem}[Correctness Theorem]
\label{equaldgm}
The persistence diagrams involving simplicial complexes  are equal to the persistence diagrams 
involving the point cloud, that is,  $\Ddgm{}{\kernel \phi^U} = \Ddgm{}{\kernel \psi}$ and $\Ddgm{}{\cokernel \phi^U} = \Ddgm{}{\cokernel \psi}$.
\end{theorem}

\paragraph{Proof sketch.}
To prove Theorem \ref{equaldgm}, we will prove, for each $\alpha \leq \beta$, that the following diagram (as well as a similar diagram involving cokernels) commutes, with
the vertical maps being isomorphisms.
\begin{align}
\label{diag:KEE}
\ldots \to  &\kernel \phi_{\alpha}^U   \to \kernel \phi_{\beta}^U  \to \ldots \notag\\
   &~\uparrow  \cong  ~~~~~~~~~~\uparrow \cong \notag\\
   \ldots \to & \kernel \psi_{\alpha}  \to  \kernel \psi_{\beta} \to \ldots.
\end{align}
Applying Theorem \ref{result:PET} then finishes the claim.
$\Ddgm{}{\kernel \phi^U} = \Ddgm{}{\kernel \psi}$ and $\Ddgm{}{\cokernel \phi^U} = \Ddgm{}{\cokernel \psi}$.

\section{Simulations}
\label{subsec:ER}

\begin{figure}[tbp]
\begin{center}
\includegraphics[scale=0.35]{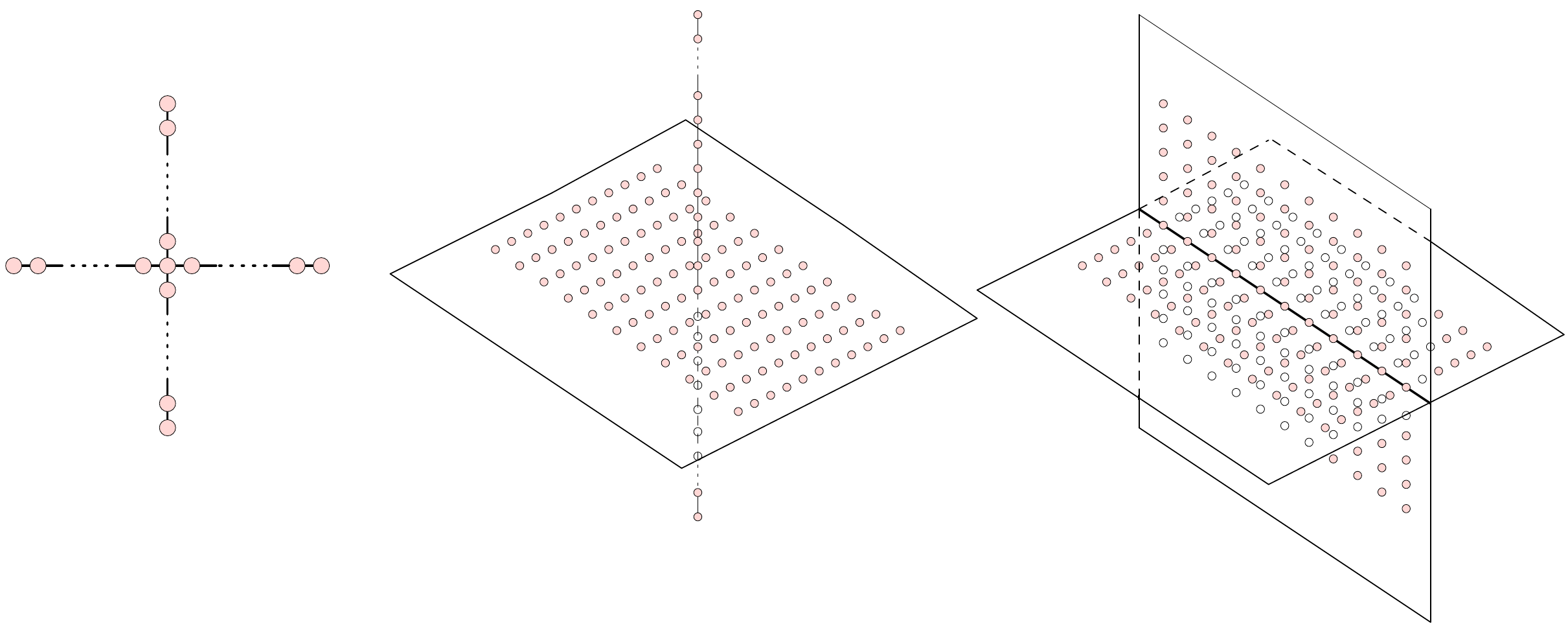} 
\end{center}
\caption[Examples of simple synthetic data.]{From left to right: points sampled from a cross; points sampled from a plane intersecting a line; points sampled from two intersecting planes. All points are located on the grid.}
\label{fig:simple-examples}
\end{figure}

We use a simulation on simple synthetic data with points sampled from grids to illustrate how the algorithm
performs. In these simulations we assume we know $\ep$, and we run our algorithm for $0 \leq \alpha \leq 2\ep$. The data sets are shown in Figure \ref{fig:simple-examples}.

\begin{figure}[tbp]
\begin{center}
\includegraphics[scale=0.35]{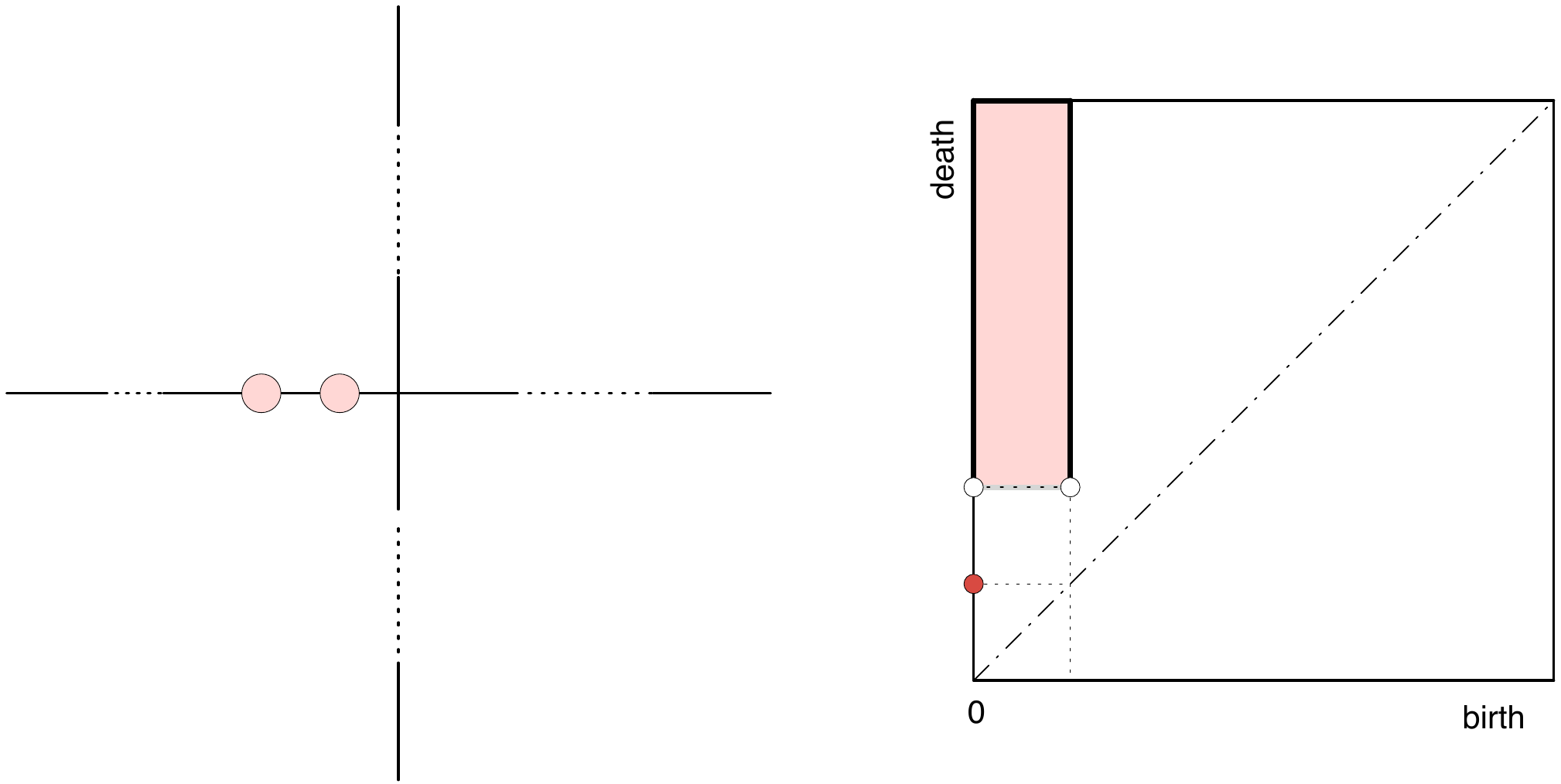} 
\includegraphics[scale=0.35]{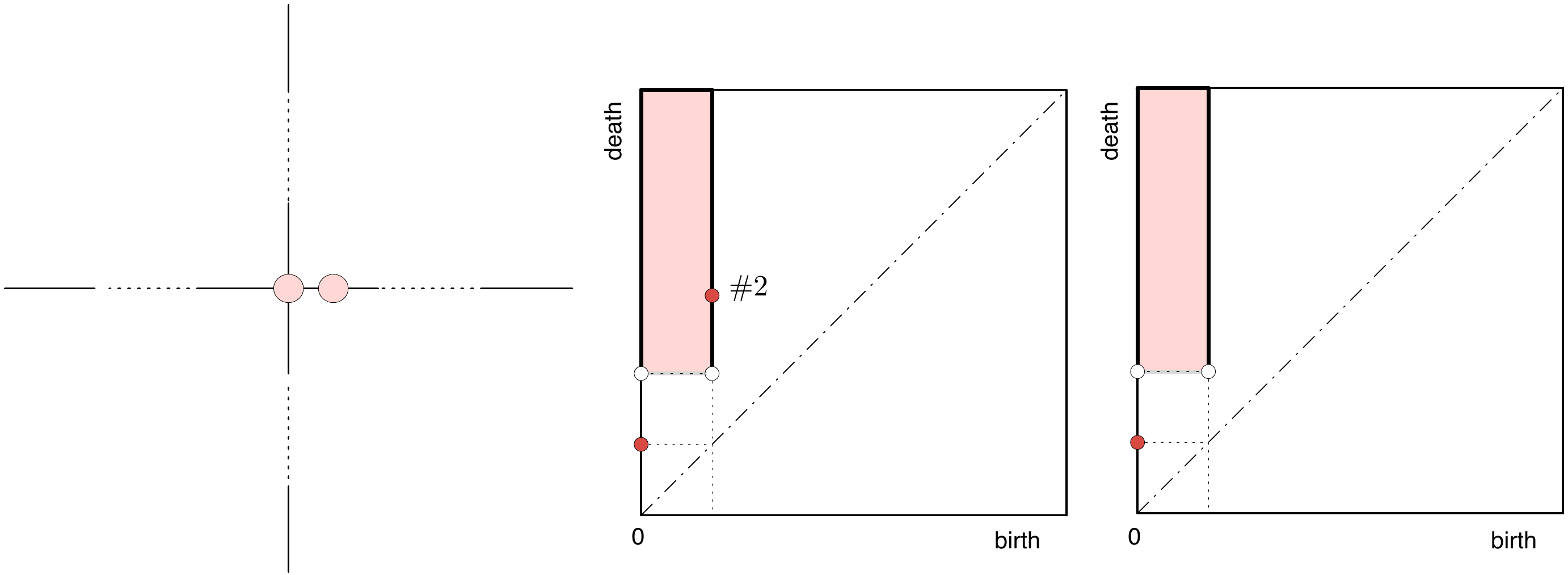} 
\end{center}
\caption[Points sampled from a cross and their corresponding ker/cok persistence diagrams.]{Top: both points are from $1$-strata. Bottom: one point from $0$-strata, one point from $1$-strata. Left part shows the locations of the points. Right part shows the ker/cok persistence diagram of two points respectively, if the diagrams are the same, only one is shown.
A number labeling a point in the persistence diagram indicates its multiplicity. }
\label{fig:crossA}
\end{figure}

\begin{figure}[tbp]
\begin{center}
\includegraphics[scale=0.35]{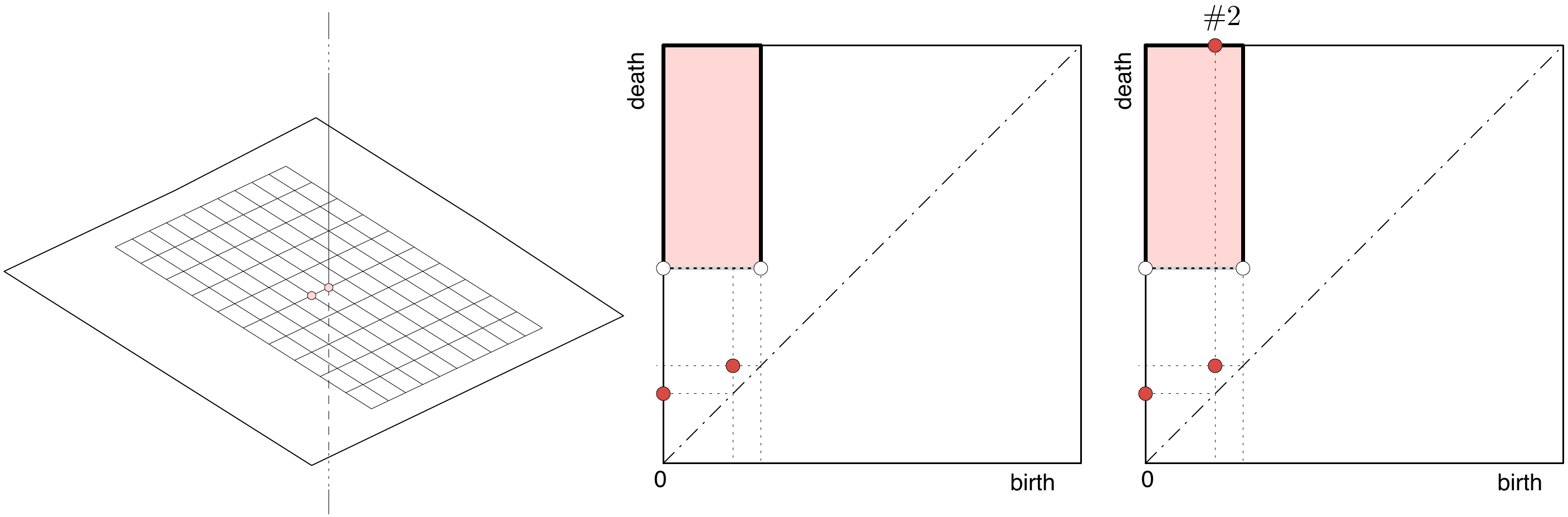} 
\includegraphics[scale=0.35]{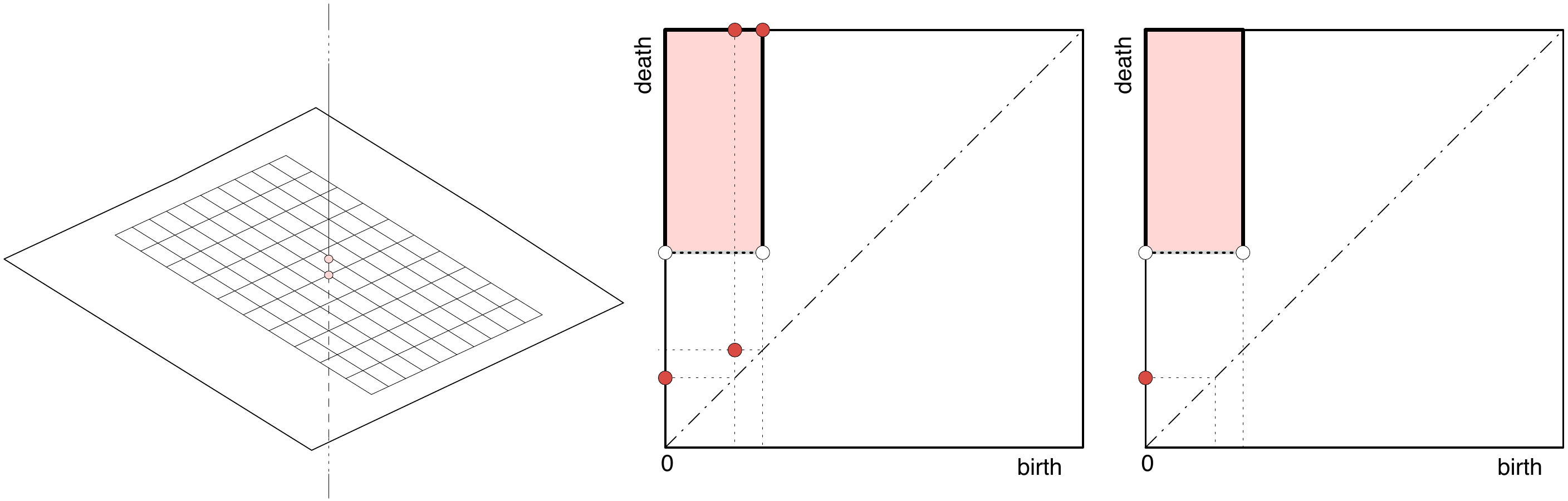} 
\includegraphics[scale=0.35]{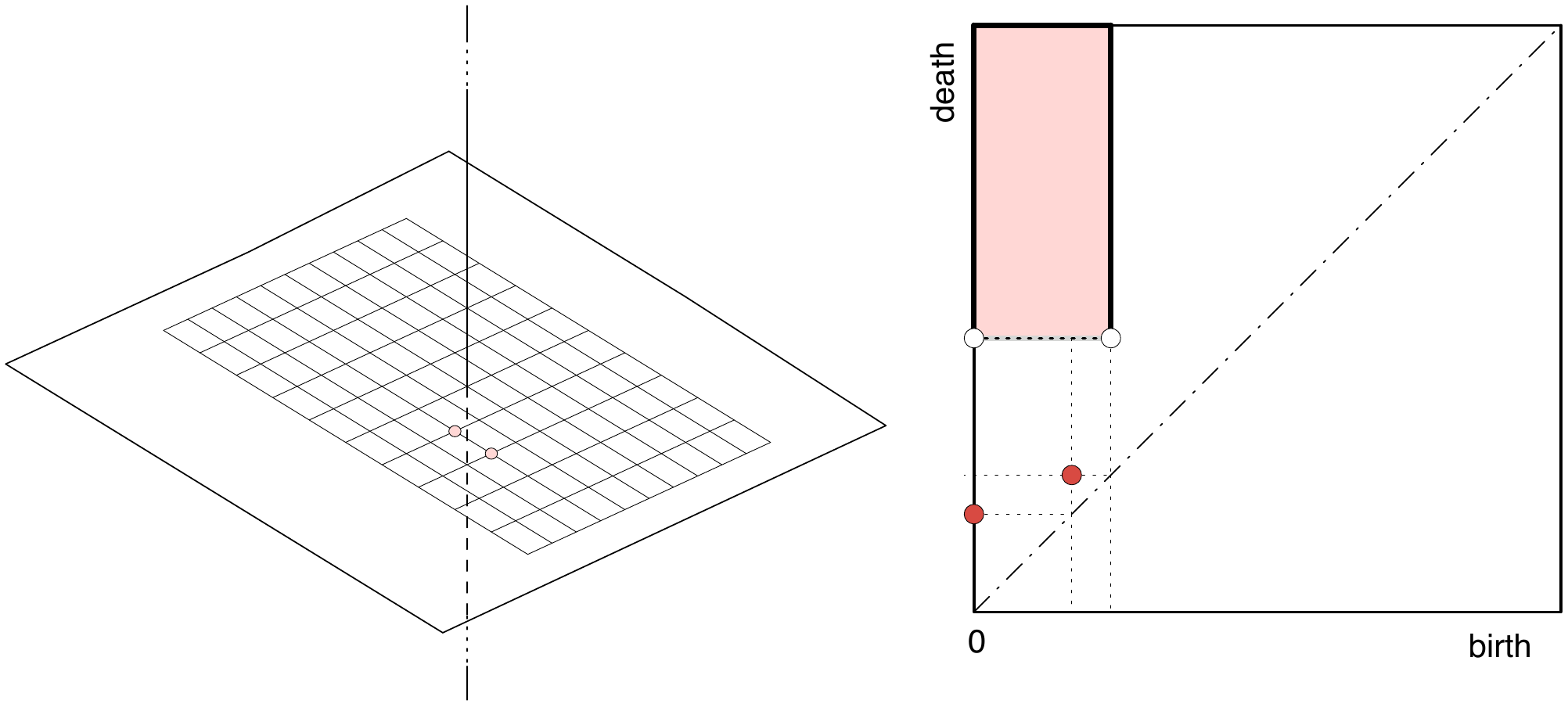} 
\end{center}
\caption[Points sampled from a plane intersecting a line, and their corresponding ker/cok persistence diagrams.]{Top: one point from $0$-strata, one point from $2$-strata. Middle: one point from $0$-strata, one from $1$-strata. Bottom: both points are from $2$-strata.
A number labeling a point in the persistence diagram indicates its multiplicity. }
\label{fig:squareC}
\end{figure}

\begin{figure}[tbp]
\begin{center}
\includegraphics[scale=0.35]{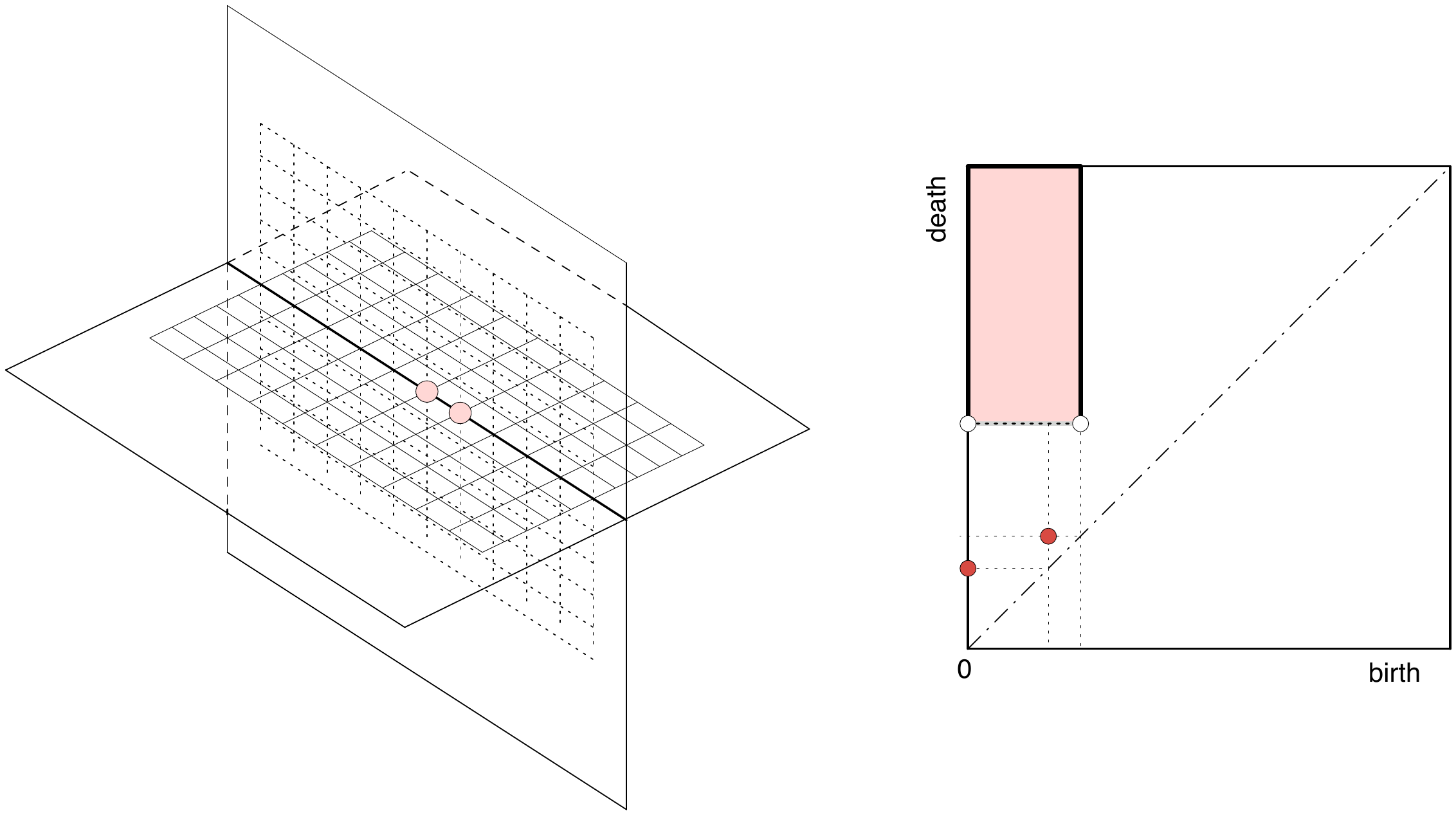} 
\includegraphics[scale=0.35]{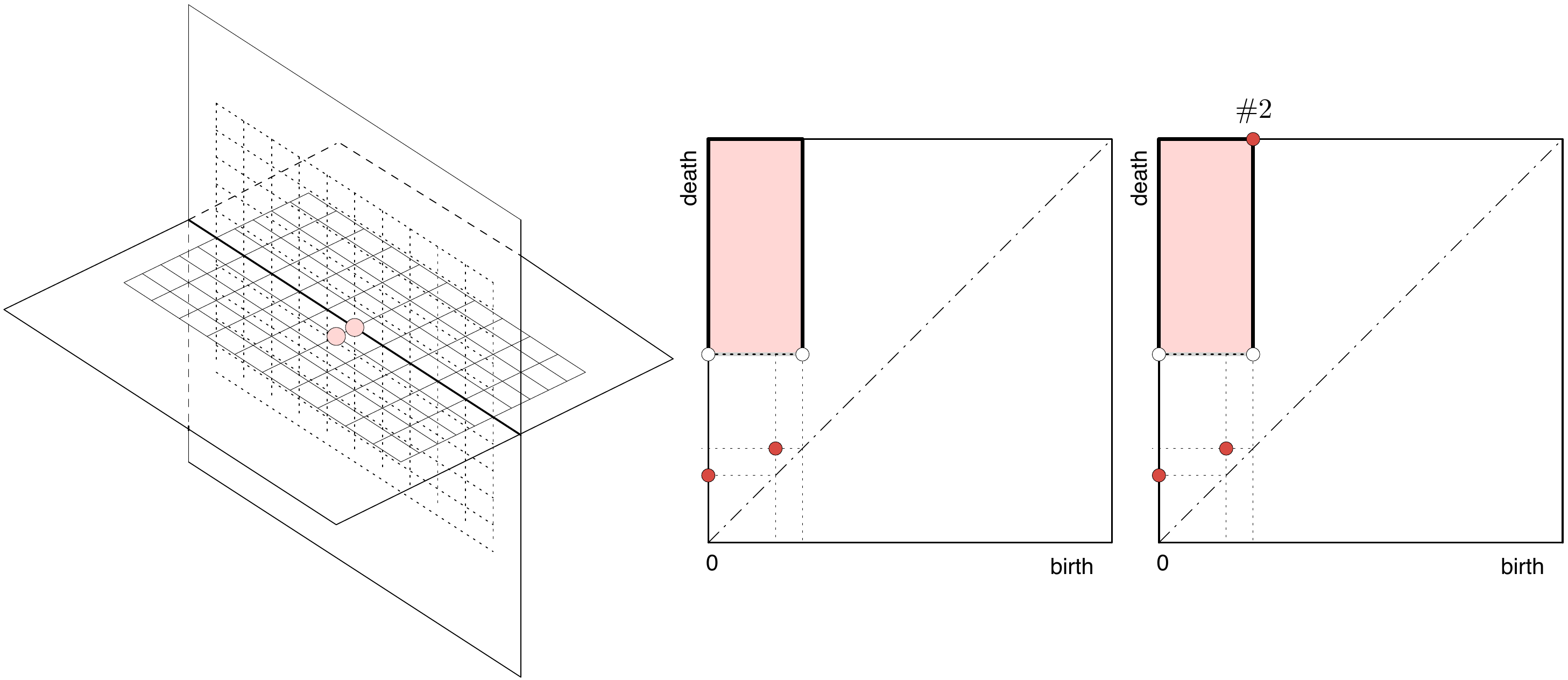}
\includegraphics[scale=0.35]{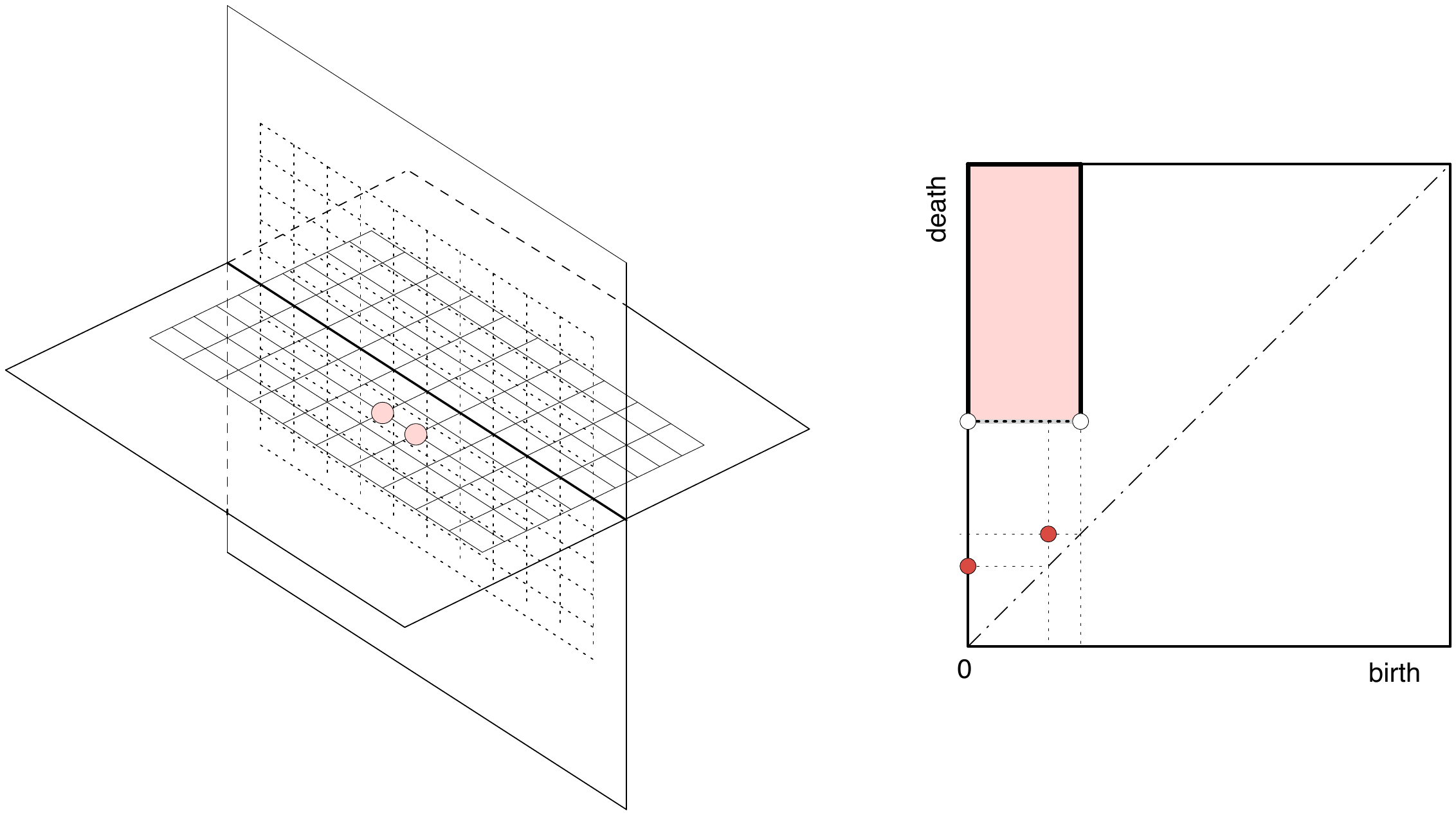}
\end{center}
\caption[Points sampled from two intersecting planes, and their corresponding ker/cok persistence diagrams.]{Top: both points from $1$-strata. Middle: one point from $1$-strata, one from $2$-strata. Bottom: both points are from $2$-strata. A number labeling a point in the persistence diagram indicates its multiplicity. }
\label{fig:twosheetA}
\end{figure}

We use the following results to demonstrate that the inference on local structure, at least for these very simple examples, is correct.
As shown in Figure \ref{fig:crossA} top, if two points are locally equivalent, their corresponding ker/cok persistence diagrams contain the empty quadrant
prescribed by our theorems, while
in Figure \ref{fig:crossA} bottom, the diagrams associated to two non-equivalent points do not contain such empty quadrants.
Similar results are shown for other data sets in Figure \ref{fig:squareC} and Figure \ref{fig:twosheetA}.

\section{Discussion}
\label{sec:Disc}

As the title of the paper suggests we have presented a first step towards learning stratified spaces. In the discussion we state some
future problems and extensions of interest. \\

\noindent{\bf Algorithmic efficiency:}  The algorithm to compute the (co)kernel diagrams from the thickened point cloud is based on an adaption of Delaunay triangulation and 
the power-cell construction. This algorithm should be quite slow when the dimensionality of the ambient space is high due to
the runtime complexity of Delaunay triangulation. One idea to address this bottleneck is to use Rips or Witness complexes \cite{SilCar2004}; at the moment, we
are not sure how to approach a proof of correctness, due to problems presented by the boundaries of the $r$-balls.
Another approach is to use  dimension reduction techniques such as principal components analysis (PCA)  or random projection
 that approximately preserve distance \cite{Cla2008} as a preprocessing step. Another idea that may work if the ambient
dimension is not too high is using faster algorithms to construct Delaunay triangulations \cite{BoiDevHor2009}. \\

\noindent{\bf Weighting local equivalence:} Currently we use a graph with $0/1$ weights based on the local
equivalence between two points. Extending this idea to assign fractional weights between points 
is appealing as it suggests a more continuous metric for local equivalence. This may also allow for
greater robustness when using spectral methods to assign points to strata.  \\

 \noindent{\bf Curvature moderated tubes:} Markus J.\ Pflaum  \cite{Pfl2001} introduced a concept called \emph{curvature moderation} 
 that regulates the behavior of the tangent spaces of a stratum 
near the boundary. In other words, a stratum is curvature moderate if it curves near the boundary in a controlled way, this includes the higher derivatives 
of the curvature. This is yet another way to describe how strata and their tubular neighborhood  are ``glued together nicely". It would be interesting to 
connect this concept to our idea of  ``local reach".     \\

 \noindent{\bf Noisy data:} Our sampling models draw points from the underlying topological space. A more general model would sample
 points that are concentrated on the topological space. A version of this type of sampling model is discussed in \cite{NiySmaWei2008}.
 It would be of interest to study how well our approach is suited to such a model. \\
 
 \noindent{\bf Adaptive sampling conditions:}  Throughout this paper we use $\ep$-approximation to characterize the similarity
of the point sample to the topological space. There are other approximation criteria that may be interesting to study and may provide
better sampling estimates. One such criterium is the $\ep$-sample \cite{Dey2007} which is adaptive in that it is proportional to the
local feature size. Another criterium of possible interest is the weak feature size \cite{ChaLie2005}. \\

\section*{\LARGE Appendices}
\setcounter{section}{0}
\def\thesection{\Alph{section}}
\section{Defining the Map $\phi$}
\label{app:phimap}

We give a more precise description of the map
\begin{align*}
\phi = \phi^U_{\alpha}: \Hgroup(B_{p}^U(\alpha), \bdr B_{p}^{U}(\alpha)) \to \Hgroup (B_{pq}^{U}(\alpha), \bdr B_{pq}^{U}(\alpha)).  
\end{align*}
The definition will be made on the chain level and will be given in terms of singular chains.

\subsection{Background}

We give here some necessary background as well as some material from algebraic topology and homological
algebra which will be needed in Appendix \ref{app:correctness}.
Most of the descriptions are adapted from \cite{Hat2002} and \cite{Mun1984}.

\paragraph{Chain homotopies.}
For our purposes, a \emph{chain complex} $\Ccal$ is a sequence of $\mathbb{Z}/2 \mathbb{Z}$- vector spaces $C_p$, one for each integer $p$, connected by boundary 
homomorphisms $\bdr_p^C: C_p \to C_{p-1}$
such that $\bdr_{p-1} \circ \bdr_p = 0$ for all $p$.
The $p$-th homology group of such a chain complex is defined by $\Hgroup_p = \kernel{\bdr_p} / \image{\bdr_{p+1}}$.

A \emph{chain map} $\eta: \Ccal \to \Dcal$ between two chain complexes is a sequence of homomorphisms $\eta_p : C_p \to D_p$ 
which commute with the boundary homomorphisms: $\bdr_{p}^D \circ \eta_p = \eta_{p-1} \circ \bdr_p^C$.
Every chain map induces a map $\eta_{*}$ between the homology groups of the two complexes.

A \emph{chain homotopy} $F$ between two chain maps $\eta,\eta': \Ccal \to \Dcal$ is a sequence
of homomorphisms $F_p: C_p \to D_{p+1}$ which satisfy the following formula for each $p$:
$\eta - \eta' = \bdr_{p+1}^D \circ F_p - F_{p-1} \circ \bdr_p^C$.
We say that $\eta$ and $\eta'$ are \emph{chain homotopic} and note that they must then induce
the same maps on homology: $\eta_*= \eta'_*$.
Finally, $\eta$ is called a \emph{chain homotopy equivalence} if there exists a chain map $\rho: \Dcal \to \Ccal$
such that $\eta \circ \rho$ and $\rho \circ \eta$ are both chain homotopic to the identify.
In this case $\eta$ and $\rho$ will both induce homology isomorphisms.

\paragraph{Singular homology.}
The \emph{standard $p$-simplex} is the subset of $\Rspace^{p+1}$ given by
\begin{align*}
\Delta_p = \{(t_0, ..., t_p) \in \Rspace^{p+1} | \sum_{i=0}^{p} t_i = 1, \forall i, t_i \geq 0\}.
\end{align*}
The $p+1$ vertices of $\Delta_p$ are points $\{e_i\} \subset \Rspace^{p+1}$ where
\begin{align*}
e_0 & = (1, 0, 0, ..., 0), \\
e_1 & = (0, 1, 0, ..., 0), \\
& ... \\
e_p & = (0, 0, 0, ..., 1).
\end{align*}
A \emph{singular $p$-simplex} of a topological space $X$ is a continuous map
$\delta: \Delta_p \rightarrow X.$
By taking formal sums of singular simplices (using binary coefficients for our purposes)
one forms $C_p(X)$, the \emph{singular chain group} of $X$ in dimension $p$.
Given points $a_0, ..., a_p$ in some Euclidean space, which need not be independent, there is a unique affine map $l$ of $\Delta_p$ that maps the vertices $e_i$ of $\Delta_p$ to $a_i$. This map defines the \emph{linear singular simplex} determined by $a_0,...,a_p$, denoted as $l(a_0, ..., a_p)$.
One then defines a boundary homomorphism $\bdr_p: \Delta_p(X) \to \Delta_{p-1}(X)$ by:
\begin{align*}
\bdr_p(\delta) = \Sigma_{i=0}^p (\delta \circ l(\ep_0,...,\hat{\ep}_i,...,\ep_p)),
\end{align*}
and defines the singular homology groups $\Hgroup_p(X)$ as above.
A continuous map $f$ from $X$ to another topological space $Y$ induces a chain map $f_{\#}: C_p(X) \to C_p(Y)$
given by the formula $f_{\#}(\delta) = f \circ \delta$, and thus also a homology map $f_*: \Hgroup_p(X) \to \Hgroup_p(Y)$.

The \emph{minimal carrier} of a singular simplex $\delta$ is its image $\delta(\Delta_p)$, and
the minimal carrier of a singular $p$-chain $\sum \delta_i$ is the union of the minimal carriers of the $\delta_i$.

\paragraph{Isomorphism between simplicial and singular homology.}
The (simplicial) homology groups of a simplicial complex $K$ and the singular homology groups of its realization $|K|$ are isomorphic.
To show an explicit isomorphism (\cite{Mun1984}), we first define a chain map
\begin{align*}
\eta: C(K) \to C(|K|)
\end{align*}
as follows \cite{Mun1984}: choose a partial ordering of the vertices of $K$ that induces a linear ordering on the vertices of each simplex of $K$. 
Orient the simplices of $K$ by using this ordering, and define
\begin{align*}
\eta([v_0, ..., v_p]) = l(v_0, ..., v_p),
\end{align*}
where $v_0 < ... < v_p$ in the given ordering.
We refer to the linear singular simplex $l(v_0, ..., v_p)$ as a \emph{simplicial linear singular simplex} and it is important in the subsequent sections. 
The chain map $\eta$ is in fact a chain equivalence as it has a \emph{chain-homotopy inverse}, for which a specific formula can be found in \cite{EilSte1952}.  
Hence the induced homology map $\eta_{*}$ provides an isomorphism of simiplicial with singular homology.

\subsection{Intersection Map Details}

We now give the full and formal definition of the homology map $\phi = \phi_{\alpha}^U,$ starting on the chain level.
For compactness, we will use the following shorthand:
\begin{align*}
X &= B_{p}^{U}(\alpha) = U_{\alpha}  \cap B_r(p),\\
B &= \bdr B_{p}^{U}(\alpha) = U_{\alpha} \cap  \bdr B_r(p),\\
S &= B_{pq}^U(\alpha) = U_{\alpha}  \cap B_r(p) \cap B_r(q),\\
A &=  U_{\alpha} \cap B_r(p) - \interior(S),\\ 
U &= U_{\alpha} \cap B_r(p) - S.
\end{align*}
Note that $X - U = S = B_{pq}^{U}(\alpha)$ and $A - U = \bdr S = \bdr B_{pq}^{U}(\alpha)$.
So to define $\phi$ we need only define a chain map $j: C(X, B) \to C(X - U, A - U)$ and
then take $\phi$ as the map induced on homology.
The map $j$ is defined as the composition $j = k \circ i$. The chain map $i: C(X,B) \to C(X,A)$ is induced
by inclusion on the second factor, while the chain map $k: C(X,A) \to C(X - U, A - U)$ is an excision, although this latter
statement requires further elaboration.

\paragraph{Excisions.}

The inclusion map of pairs $(X - U, A - U) \to (X,A)$ is called an excision if it induces a homology isomorphism;
in this case one says that $U$ can be excised.
We will make use of the following two results about excision (see, e.g., \cite{GreHar1981}):
\begin{theorem}(Excision Theorem)
\label{thm:excision2}
If the closure of $U$ is contained in the interior of $A$, then $U$ can be excised.  
\end{theorem}
\begin{theorem}(Excision Extension)
\label{thm:extension}
Suppose $V \subset U \subset A$ and 
\begin{itemize}
\item[(i)] $V$ can be excised.
\item[(ii)] $(X - U, A - U)$ is a deformation retract of $(X - V, A - V)$.
\end{itemize}
Then $U$ can be excised.
\end{theorem}
In our context, the closure of $U$ need not be contained in the interior of $A$, and so we must define a suitable $V \subset U$.
Although there are many ways to do this, one direct way is to choose some small enough positive $\delta$.
\begin{align*}
I & = U_{\alpha} \cap \bdr (B_r(p) \cap B_r(q)) \cap \closure(U),\\
I_\delta & = \{ x \in \closure(U) | d_I(x) \leq \delta\},\\
V & = U - I_{\delta},
\end{align*} 
where $d_I(x)$ is the distance from $x$ to the set $I$.
It is straightforward to verify that $V \subset U \subset A$ satisfies the hypotheses of Theorem \ref{thm:extension}.
In other words, the inclusion of pairs $(X - V, A - V) \to (X,A)$ is an excision; its induced chain map has a chain-homotopy
inverse, which we denote as $s: C(X,A) \to C(X - V, A- V)$.
Finally, we define $k = r_{\#} \circ s$, where $r_{\#}$ is the chain map induced by the retraction $r: (X - V, A - V) \to (X - U, A - U)$.

\paragraph{Subdivision.}

In order to fully carry out the analysis in Appendix B, we must first decompose the maps $i$ and $k$ through subdivision.
Given a topological space $X$ and a collection $\A$ of subsets of $X$ whose interiors form an open cover of $X$,
a singular simplex of $X$ is said to be \emph{$\A$-small} if its image set is entirely contained in a single element of $\A$.
For each dimension $p$, the chain group $C_p^{\A}(X)$ is the subgroup of $C_p(X)$ spanned by the $\A$-small singular $p$-simplices.
These groups form a chain complex, with homology $\Hgroup^{\A}(X)$.
Of course, any singular simplex on $X$ can be subdivided into a sum of $\A$-small simplices, so it is plausible,
and in fact true (\cite{Hat2002}), that
the inclusion $C^{\A}(X) \to C(X)$ is a chain homotopy equivalence.

Returning to our context, we set $\A = \{X - V, A\}$ and denote by $l$ the chain inclusion $C^{\A}(X,A) \to C(X,A)$.
We also let $\rho: C(X,B) \to C^{\A}(X,B)$ be the chain homotopy inverse of the chain inclusion $C^{\A}(X,B) \to C(X,B)$,
and let $t: C^{\A}(X,B) \to C^{\A}(X,A)$ be the chain map induced by inclusion on the second factor.
Finally we note that $i = l \circ t \circ \rho$.

We also decompose $k$ as $k = r_{\#} \circ \eta \circ \rho$, where $\eta$ is the chain homotopy
inverse of the chain map $C(X - V, A - V) \to C(X,A) \to C^{\A}(X,A)$.

\paragraph{Summary.}

To summarize, our map $\phi = j_*$, where $j$ is the chain map defined by the following sequence of chain maps
\begin{align*}
j = k \circ i = (r_{\#} \circ \eta \circ  \rho) \circ (l \circ  t \circ \rho).
\end{align*}
Following the same framework as above, we also define a chain map $j'$ and its induced homology map 
$\phi' = j'_{*}: \Hgroup(B_p^U(\alpha),P_0(\alpha)) \to \Hgroup(B_{pq}^U(\alpha),Z_0(\alpha))$
simply by adopting the notation:
\begin{align*}
X &= B_{p}^{U}(\alpha) = U_{\alpha}  \cap B_r(p),\\
B' &= P_0(\alpha),\\
S &= B_{pq}^U(\alpha) = U_{\alpha}  \cap B_r(p) \cap B_r(q),\\
A' &=  U_{\alpha} \cap B_r(p) - S + Z_0(\alpha),\\ 
U &= U_{\alpha} \cap B_r(p) - S,\\
I & = U_{\alpha} \cap \bdr (B_r(p) \cap B_r(q)) \cap \closure(U),\\
I_\delta & = \{ x \in \closure(U) | d_I(x) \leq \delta\},\\
V & = U - I_{\delta},
\end{align*}
defining our open cover to be $\A' = \{ X - V, A'\}$, and otherwise proceeding exactly as before.

Similarly, we create a chain map $f'$ which induces $\psi' = f'_{*}: \Hgroup(|\sd{L}|,|\sd{L}_0|) \to \Hgroup(|\sd{K}|,|\sd{K}_0|)$,
using the notation
\begin{align*}
X^{''} & = |\sd L|,\\
B^{''} & = |\sd L_0|,\\
A^{''} & = (|\sd L| - |\sd K|) \cup |\sd K_0|,\\
U^{''} & = |\sd L| - \interior |\sd K|,\\
I & = |\sd K| \cap \closure(U'),\\
I_\delta & = \{ x \in \closure(U') | d_I(x) \leq \delta\},\\
V^{''} & = U' - I_{\delta},
\end{align*}
with $\A'' = \{A'', X'' - V''\}$.

\section{Algorithm Details}
\label{app:purburbation}

We give the details in constructing the simplicial complexes described in our algorithm.
The various simplicial complexes, $L$, $L_0$, $K$ and $K_0$, are the nerves of collections of convex sets. 
Here we go through the construction of $L$; construction of the others is similar.

\paragraph{Implicit Perturbations.}

A direct approach to constructing $L$, the nerve of the collection $B$, runs into difficulties as the corners of the convex sets created by the bisector $\Bisector$ can be shared by many sets. 
To cope with this difficulty, we perturb these convex sets ever so slightly so that they meet in general position.
Note that this is not done by perturbing the bisector; rather, it is done by decomposing the bisector into pieces. 

We are interested in the restricted Voronoi diagram of the sublevel sets inside the ball $B_r(p)$, 
which we denote as $\V = \voi{U|U_{\alpha} \cap B_r(p)}$. 
The \emph{restricted Voronoi cell} of $u_i$ is defined as $V(u_i|U_{\alpha} \cap B_r(p)) = V(u_i) \cap B_r(p)$.


Given $\V$, we create three sets of points.
Let $\T_{pq}$ be the set of points $u_i \in U$ whose restricted Voronoi cells have non-trivial intersection with the bisector: 
$V(u_i|U_{\alpha} \cap B_r(p)) \cap \Bisector \neq \emptyset$.
We impose an ordering of points in $\T_{pq}$, 
w.l.o.g., let the ordered set be $\T_{pq} = \{x_1, x_2, ..., x_m\}$.
$\T_{p}$ is the set of points $u_i \in U$ which are not in $\T_{pq}$ and are closer to $p$ than they are to $q$.
Similarly, $\T_{q}$ is the set of points $u_i \in U$, which are not in $\T_{pq}$ and are closer to $q$.

\begin{figure}[tbp]
\begin{center}
\includegraphics[scale=0.3]{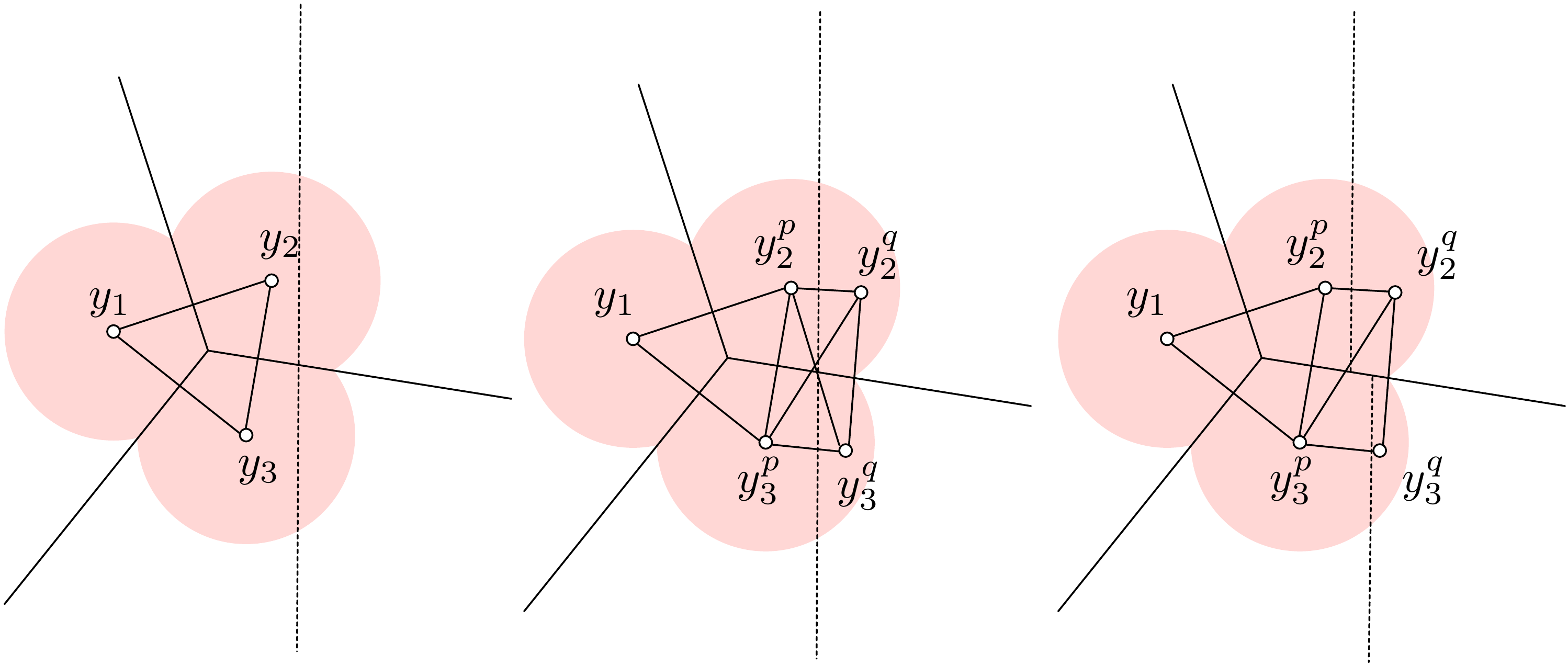} 
\end{center}
\caption{An example of the implicit perturbation. Dotted lines are the bisectors. A simplex $[y_1, y_2, y_3] \in L'$ is shown in the left. The simplices in $L$ and in $\tilde{L}$ are shown in the middle and right, respectively.}
\label{fig:strata-perturb}
\end{figure}

By construction, the bisector $\Bisector$ intersects the restricted Voronoi cells of the points in $\T_{pq}$. 
We denote these corresponding intersections as $\{\Bisector_1, \Bisector_2, ..., \Bisector_m\}$.
We perturb each $\Bisector_i$ slightly such that no two pieces are collinear.
Note that $\Bisector_i$ is perpendicular to the direction $q - p$.
One possible choice of perturbation would move each $\Bisector_i$ within the restricted Voronoi cell along the direction $q - p$
for $i\ep$ distance, where $\ep$ is sufficiently small. An example in $\Rspace^2$ is shown in Figure \ref{fig:strata-perturb}.

Given such a perturbation, we let $\tilde{\Acal}$ be the resulting collection of perturbed convex sets, and we compute 
$\tilde{L} = \Nerve(\tilde{\Acal})$ instead of $L = \Nerve(\Acal)$.
By the properties of nerve construction, $\Nerve(\tilde{\Acal}) \simeq \bigcup \tilde{\Acal}$, $\Nerve(\Acal) \simeq \bigcup \Acal$.
Since $\bigcup \tilde{\Acal} = \bigcup \Acal$, then we have $\tilde{L} \simeq L$.
We now describe how we construct $\tilde{L}$.

\paragraph{Case Analysis.}
Let $L'$ be the restricted Delaunay triangulation, $L' = \del{U | U_{\alpha} \cap B_r(p)}$.
We read the simplicies from $\tilde{\Acal}$ without explicit perturbations. 
Specifically, we follow a set of rules, described below, to construct $\tilde{L}$ from $L'$.

The bisector divides the restricted Voronoi cell of a point $x \in \T_{pq}$ into two convex sets. 
Let $x^p$ represent the perturbed convex set closer to $p$ in the nerve construction, and let $x^q$ represent the other set.
Let $\sigma$ be a simplex in $L'$ with $k$ vertices, that is, $\sigma = [y_1,y_2, ..., y_k]$.
There are seven cases regarding the membership of the points $\{y_1, y_2, ..., y_k\}$.
\begin{itemize}
\item[1.] 
All $y_i \in \sigma$ belong to $\T_p$.
We add the simplex $[y_1, y_2, ..., y_k]$ to $\tilde{L}$.
\item[2.]
All $y_i \in \sigma$ belong to $\T_q$.
Same as case 1.
We add the simplex $[y_1, y_2, ..., y_k]$ to $\tilde{L}$.
\item[3.]
All $y_i \in \sigma$ belong to $\T_{pq}$.
Suppose $\{y_1, y_2, ..., y_k\}$ are sorted according to the ordering in $\T_{pq}$.
We add the following simplicies and their faces to $\tilde{L}$:
\begin{align*}
&[y_1^p, ..., y_m^p, y_1^q]\\
&[y_2^p, ..., y_m^p, y_1^q, y_2^q]\\
&[y_3^p, ..., y_m^p, y_1^q, y_2^q, y_3^q]\\
&...\\
&[y_m^p, y_1^q, y_2^q, ..., y_m^q]
\end{align*}
\item[4.] 
Some $y_i$ are in $\T_p$, the rest are in $\T_{pq}$.
Suppose $\{y_{i_1}, ..., y_{i_n}\} \subseteq T_{p}$ and $\{y_{j_1}, ..., y_{j_l}\} \subseteq T_{pq}$.
We add $[y_{i_1}, ..., y_{i_n}, y^p_{j_1}, ..., y^p_{j_l}]$ to $\tilde{L}$.
\item[5.] 
Some $y_i$ are in $\T_q$, the rest are in $\T_{pq}$.
Similar to case 4, 
suppose $\{y_{i_1}, ..., y_{i_n}\} \subseteq T_{q}$ and $\{y_{j_1}, ..., y_{j_l}\} \subseteq T_{pq}$.
We add $[y_{i_1}, ..., y_{i_n}, y^q_{j_1}, ..., y^q_{j_l}]$ to $\tilde{L}$.
\item[6.] 
Some $y_i$ are in $\T_p$, the rest are in $\T_q$.
We show that Case 6 is impossible. 
Choose $y_i \in \T_p$ and $y_j \in \T_q$ such that $y_i$ and $y_j$ are connected by an edge.
Since $y_i$ and $y_j$ are on the opposite sides of $\Bisector$, this edge must intersect $\Bisector$ at some point $z$.
Then their corresponding restricted Voronoi cells, $V(y_i | U_{\alpha} \cap B_r(p))$ and $V(y_j | U_{\alpha} \cap B_r(p))$,
must meet at a Voronoi face, which contains the point $z$.
Suppose that the Voronoi face is in general position, that is, it is not parallel to $\Bisector$.
Then $\Bisector$ intersects the Voronoi face, by definition, $y_i$ and $y_j$ must belong to $\T_{pq}$.
This is a contradiction.
\item[7.] 
Some $y_i$ are in $\T_p$, some are in $\T_q$, and the rest are in $\T_{pq}$.
We show that case 7 is impossible using the same proof in case 6.
\end{itemize}

A simple example is shown in Figure \ref{fig:strata-perturb}.
Given $[y_1, y_2, y_3] \in L'$, simplex $[y_1, y_2^p, y_3^p]$ is added to $\tilde{L}$ according to case 4.
Given $[y_2, y_3] \in L$, simplices $[y_2^p, y_3^p, y_2^q]$, $[y_3^p, y_2^q, y_3^q]$ and their faces are added to $\tilde{L}$ according to case 3.

In summary, we construct $\tilde{L}$ by iterating through all simplices $\sigma$ in $L'$, adding new simplicies to $\tilde{L}$ constructed from $\sigma$ following the above cases.


\section{Algorithmic Correctness}
\label{app:correctness}

\begin{figure}[tbp]
\begin{center}
\includegraphics[scale=0.5]{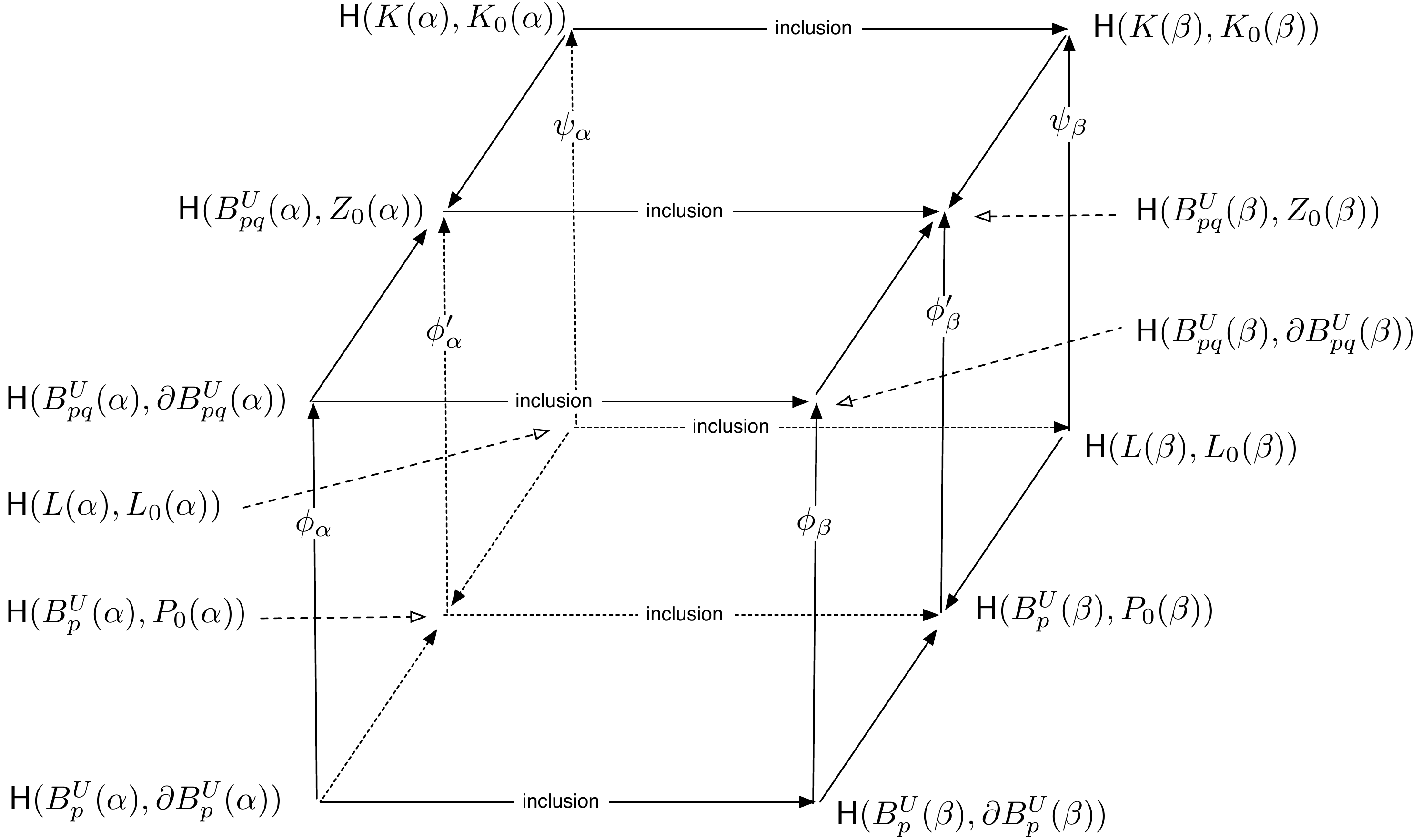} 
\end{center}
\caption{Two adjacent commuting cubes.}
\label{fig:strata-doublecube}
\end{figure}

We prove the correctness of the algorithm described in Section \ref{subsec:DC} by proving Theorem \ref{equaldgm}.
More precisely, we will prove that diagram \ref{diag:KEE} commutes, with the vertical arrows being isomorphisms, for some arbitrary
but fixed choice of $\alpha < \beta$; we will omit the very similar argument about cokernels.
The proof is unfortunately lengthy, and at times a bit technical, for 
in order to prove our statements about diagram \ref{diag:KEE}, we must also prove similar statements about several
other interlocking diagrams.
For sanity and clarity of presentation, we first exhibit all the diagrams at once in the form of the following double-cube (Figure \ref{fig:strata-doublecube}).

\subsection{Bottom Face}

The bottom face of the double-cube has been detached and drawn in diagram \ref{dgm:bottom}.
The horizontal maps in the upper square are induced by inclusion of pairs, and so the upper square
certainly commutes.
\begin{align}
\label{dgm:bottom}
 \Hgroup(B_{p}^{U}(\alpha), \bdr B_{p}^{U}(\alpha))
  &\xrightarrow{j_{\alpha}^{\beta}} 
 \Hgroup(B_{p}^{U}(\beta), \bdr B_{p}^{U}(\beta))\notag\\
 \downarrow {i_{\alpha}}~~~~~~~~~~~~&~~~~~~~~~~~~\downarrow{i_{\beta}}~~~~~~~\notag\\
\Hgroup(B_{p}^{U}(\alpha), P_0(\alpha))
&\xrightarrow{j_{\alpha}^{\beta}}
\Hgroup(B_{p}^{U}(\beta), P_0(\beta))\notag\\
 \uparrow{h_{\alpha}}~~~~~~~~~~~~&~~~~~~~~~~~~\uparrow{h_{\beta}}~~~~~~~\notag\\
\Hgroup(L(\alpha), L_0(\alpha))
 &\xrightarrow{g_{\alpha}^{\beta}}
\Hgroup(L(\beta), L_0(\beta)).
\end{align}
We have already shown that the two vertical maps in the upper square are isomorphisms; this was the content of
the Power Cell Lemma in Section \ref{subsec:DC}.
To show that the vertical maps in the lower square are isomorphisms requires a bit more work.
We make use of the following lemma, proven in \cite{BenCohEde2007}.
\begin{lemma}[General Nerve Subdivision Lemma (GNSL)]
Let $\C$ be the collection of maximal cells of a CW complex, each a convex set in $\Rspace^k$.
Define $f: |\sd{N}| \rightarrow \cup \C$ by piecewise linear interpolation of its values at the vertices. 
If $f(\hat{\sigma})$ is contained in the intersection of the cells that correspond to the vertices of $\sigma$, for each simplex $\sigma \in N$, then $f$ is a homotopy equivalence.
\end{lemma}
The vertical isomorphisms in the bottom square then follow from this next lemma, where we may of course replace $\alpha$ with $\beta$ if we wish.
\begin{lemma}[Nerve Subdivision Lemma]
Define $h = h_{\alpha} : |\sd{L(\alpha)}| \rightarrow B_p^U(\alpha)$ on the vertices $\hat{\sigma}$ of $\sd{\L(\alpha)}$ by the formula  
\begin{align*}
h_{\alpha}(\hat{\sigma}) & = \arg \min_{x \in V^{\sigma} \cap U_{\alpha} \cap B_r(p)} d_U^2(x) - d_p^2(x),
\end{align*}
and extend it by linear interpolation.
Then $h_{\alpha}$ is a homotopy equivalence of pairs from $(|\sd{L(\alpha)}|,|\sd{L_0(\alpha)}|)$ to $(B_p^U(\alpha),P_0(\alpha))$.
\end{lemma}

\begin{proof}
By construction, 
$h(\hat{\sigma})$ is contained in the intersection of the cells that correspond to the vertices of $\sigma$.
By the GNSL then, $h$ is a homotopy equivalence. 

Now we need to prove that the restriction of $h$ to $\sd{L_0(\alpha)}$ is also a homotopy equivalence. 
Let $\sigma \in L_0(\alpha)$ and put 
$h(\hat{\sigma}) = z$. 
For purposes of contradiction, suppose $z \notin P_0(\alpha)$. This means that $z \in \interior P(\alpha)$, by definition, 
and hence $d_U(z) ^ 2  - d_p(z)^2  > \alpha^2 - r^2$.

Now choose some $z' \in V^{\sigma} \cap P_0$, which must exist since $\sigma \in L_0(\alpha)$.
Then by definition we have $d_p(z')^2 - r^2 \geq d_U(z')^2 - \alpha^2$, or $d_U(z') ^ 2  - d_p(z')^2  \leq \alpha^2 - r^2$. 
Combining the above inequalities, we have 
$d_U(z') ^ 2  - d_p(z')^2  \leq \alpha^2 - r^2 < d_U(z) ^ 2  - d_p(z)^2$, which contradicts the assumption that $h(\hat{\sigma}) = z$.
We conclude that $z \in V_{\sigma} \cap P_0(\alpha)$. Applying the GNSL once more finishes the proof. \eop
\end{proof}

To show that the lower square commutes, we put $e = j_{\alpha}^{\beta} \circ h_{\alpha}$ and $e' = h_{\beta} \circ g_{\alpha}^{\beta}$, and we consider the map 
$H:|L(\alpha)| \times [0,1] \rightarrow U_{\alpha} \cap B_r(p)$ defined by
$H(x, t) = h_{\alpha_t} \circ g_{\alpha}^{\alpha_t}(x)$, where $\alpha_t = (1-t)\alpha + t\beta$.
Since the maps $g$ and $j$ are inclusions and the maps $h$ vary continuously with $\alpha$,
$H$ is a homotopy between $e$ and $e'$. This implies that the induced homomorphisms between the corresponding
homology groups are the same, $e_{*} = e'_{*}$. 

\subsection{Top Face}

We detach the top face of Figure \ref{fig:strata-doublecube}, drawing it in diagram \ref{dgm:top}. 
As before, we prove that all vertical maps are isomorphisms. The commutativity of the two smaller squares follows from nearly identical arguments to
the ones used for the bottom face.
\begin{align}
\label{dgm:top}
\Hgroup(B_{pq}^{U}(\alpha), \bdr B_{pq}^{U}(\alpha))
 &\rightarrow  
\Hgroup(B_{pq}^{U}(\beta), \bdr B_{pq}^{U}(\beta))\notag\\
 \downarrow {i'_{\alpha}}~~~~~~~~~~~~&~~~~~~~~~~~~\downarrow{i'_{\beta}}~~~~~~~\notag\\
\Hgroup(B_{pq}^{U}(\alpha), Z_0(\alpha))
 &\rightarrow 
\Hgroup(B_{pq}^{U}(\beta), Z_0(\beta))\notag\\
 \uparrow{h'_{\alpha}}~~~~~~~~~~~~&~~~~~~~~~~~~\uparrow{h'_{\beta}}~~~~~~~\notag\\
\Hgroup(K(\alpha), K_0(\alpha))
 &\rightarrow
\Hgroup(K(\beta), K_0(\beta))
\end{align}
The Intersection Power Cell Lemma tells us that the vertical maps in the top square are isomorphisms.
As promised, we give the proof of this lemma here, repeating the statement for completeness.
\begin{lemma}[Intersection Power Cell Lemma] 
Assume $B_r(p) \cap B_r(q) - Z_0(\alpha) \neq 0$.
The identity $i'$ on $B_{pq}^U(\alpha)$ is a homotopy equivalence of pairs between 
$(B_{pq}^U(\alpha),\bdr B_{pq}^U(\alpha))$ and
$(B_{pq}^U(\alpha), Z_0(\alpha))$.
\end{lemma}
\begin{proof}
It suffices to show that the restriction of the identity,
$i' = i'_{\alpha}: \bdr_{pq}^U(\alpha) \rightarrow Z_0(\alpha),$
is a homotopy equivalence. 
To do this, we first define a retraction
$j: Z_0(\alpha) \to \bdr_{pq}^U(\alpha)$
as follows.
Fix a point $y \in \interior Z(\alpha),$ recalling that this set is nonempty by assumption.
For each point $x \in Z_0(\alpha)$, we consider the unique ray starting at $y$ and passing through
$x$, and we let $x' = j(x)$ denote its intersection with $\bdr(B_r(p) \cap B_r(q))$. 
Note that $x' \in Z_(\alpha) \subseteq U(\alpha)$, and so $j$ is certainly well-defined.
That $j$ is a retraction, meaning $j \circ i'$ is the identity on $\bdr B_{pq}^U(\alpha)$, is obvious.
On the other hand, the map
$$\lambda: Z_0(\alpha) \times [0,1] \times Z_0(\alpha)$$
defined by $\lambda(x,t) = (1-t)x + tx'$ is a homotopy between
$i' \circ j$ and the identity map on $Z_0(\alpha)$, and so the claim follows. \eop
\end{proof}

To prove that the vertical maps in the lower square are isomorphisms, we again make use of the GNSL.
\begin{lemma}[Intersection Nerve Subdivision Lemma (INSL)]
Define $h' = h'_{\alpha} : |\sd{K(\alpha)}| \rightarrow B_{pq}^U(\alpha)$  by setting
\begin{align*}
h'_{\alpha}(\hat{\sigma}) & = \arg \min_{x \in V^{\sigma} \cap U_{\alpha} \cap B_r(p) \cap B_r(q)} 
\min \{d_U^2(x) - d_p^2(x), d_U^2(x) - d_q^2(x)\},
\end{align*}
where $\hat{\sigma}$ is the barycentre of $\sigma \in K{\alpha}$, and then extending by linear interpoloation.
Then $h'$ is a homotopy equivalence of pairs between $(|\sd{K(\alpha)}|,|\sd{K_0(\alpha)}|)$ and $(B_{pq}^U(\alpha),Z_0(\alpha))$. 
\end{lemma}
\begin{proof}
The proof is quite similar to that of the NSL.
By construction, 
$h'(\hat{\sigma})$ is contained in the intersection of the cells that correspond to the vertices of $\sigma$, and so we need only
prove that the restriction of $h'$ to the barycentric subdivision of $K_0(\alpha)$ is also a homotopy equivalence. Let $\sigma \in K_0(\alpha)$ and put
$h'(\hat{\sigma}) = z$.  

Suppose $z \notin Z_0(\alpha)$, and thus $z \in \interior Z(\alpha)$.  
By definition then, $d_p(z)^2 - r^2  < d_U(z)^2 - \alpha^2$ and $d_q(z)^2 - r^2  < d_U(z)^2 - \alpha^2$. 
In other words, $\min\{d_U^2(x) - d_p^2(x), d_U^2(x) - d_q^2(x)\} > \alpha^2 - r^2$.

Choose some point $z' \in V_{\sigma} \cap Z_0(\alpha)$.
Then one of the following inequalities must hold:
$d_p(z')^2 - r^2 \geq d_U(z')^2 - \alpha^2$, 
or $d_q(z')^2 - r^2 \geq d_U(z')^2 - \alpha^2$. 
That is, $\min\{d_U^2(z') - d_p^2(z'), d_U^2(z') - d_q^2(z')\} \leq \alpha^2 - r^2.$

Therefore, combining both inequalities,
$\min\{d_U^2(z') - d_p^2(z'), d_U^2(z') - d_q^2(z')\} 
\leq \alpha^2 - r^2 
<  \min\{d_U^2(z) - d_p^2(z), d_U^2(z) - d_q^2(z)\} 
$,
which contradicts the definition of $z$. \eop
\end{proof}

\subsection{Left and Right Faces}

We now come to the final and most complicated part of the correctness proof, involving the left face (diagram \ref{dgm:left}) of the double-cube; of course, everything we prove
here will also hold for the right face. We have already established that all vertical maps are isomorphisms, and now must show that both squares commute.
\begin{align}
\label{dgm:left}
 \Hgroup(B_{p}^{U}, \bdr B_{p}^{U})
  &\xrightarrow{\phi}
 \Hgroup(B_{pq}^{U}, \bdr B_{pq}^{U})\notag\\
 \downarrow {i_{*}}~~~~~~&~~~~~~~\downarrow{i'_{*}}~~~~~~~\notag\\
\Hgroup(B_{p}^{U}, P_0)
&\xrightarrow{\phi'}
\Hgroup(B_{pq}^{U}, Z_0)\notag\\
 \uparrow{h_{*}}~~~~~&~~~~~~~\uparrow{h'_{*}}~~~~~~~\notag\\
\Hgroup(L, L_0)
 &\xrightarrow{\psi}
\Hgroup(K, K_0).
\end{align}
The top square will in fact commute even on the chain level.
The bottom square is a little more complicated, and we start by addressing this first.

In diagram \ref{dgm:leftbottomchain}, this bottom square has been expanded into two smaller squares of chain groups
connected by chain maps.
We show that the lower of these squares commutes on the chain level, and that the two choices of path across the upper
square are connected by a chain homotopy. 
\begin{align}
\label{dgm:leftbottomchain}
 C(B_{p}^{U}, P_0)
&\xrightarrow{j'}
C(B_{pq}^{U}, Z_0) \notag\\
\uparrow{h_{\#}}~~~~~~&~~~~~~~~~~~~\uparrow{h'_{\#}}~~~~~~~\notag\\
C(|\sd L|, |\sd L_0|)
 &\xrightarrow{f'}
C(|\sd K|, |\sd K_0|).\notag\\
 \uparrow{\eta}~~~~~~~~~&~~~~~~~~~~~~\uparrow{\eta}~~~~~~~\notag\\
C(\sd L, \sd L_0)
 &\xrightarrow{f}
C(\sd K, \sd K_0).
\end{align}

\subsubsection{Map Details}

First we need to discuss two of the horizontal chain maps from diagram \ref{dgm:leftbottomchain} in more explicit detail.

\paragraph{Upper map.}

\begin{figure}[tbp]
 \begin{center}
  \includegraphics[scale=0.3]{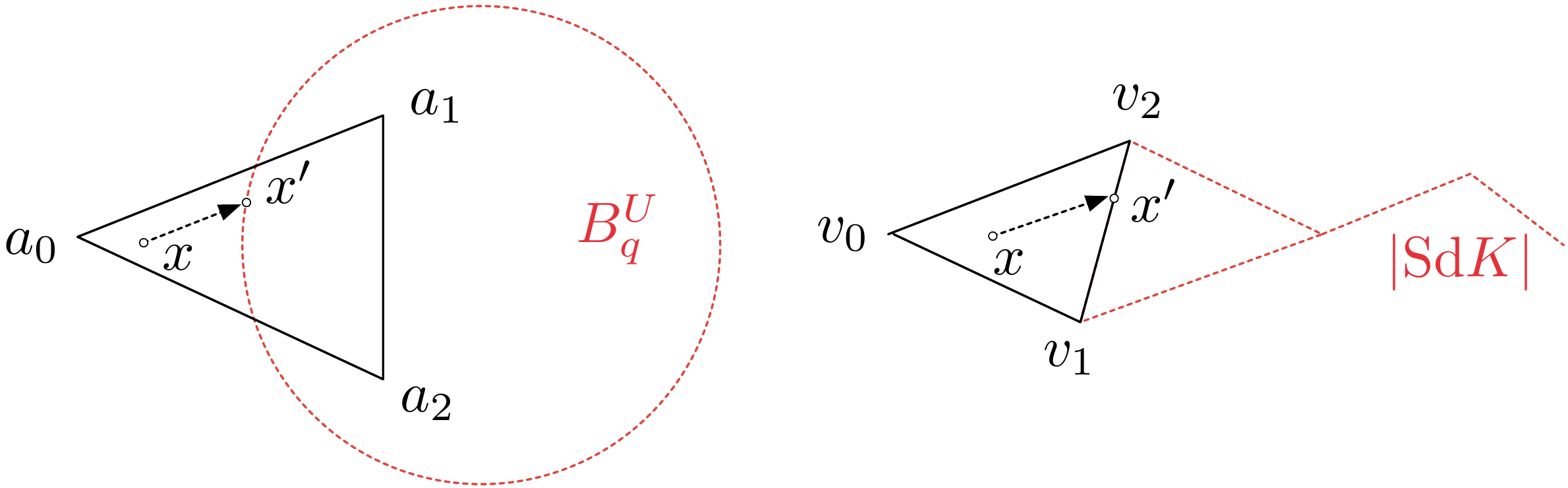}
 \end{center}
\caption[Map $j'$ and $f'$ for a linear singular simplex.]
{Left: map $j'$ for a linear singular simplex $l(a_0, a_1, a_2)$. Right: map $f'$ for a simplicial linear singular simplex $l(v_0, v_1,v_2)$.}
\label{fig:jdetail}
\end{figure}

\begin{figure}[tbp]
 \begin{center}
  \includegraphics[scale=0.3]{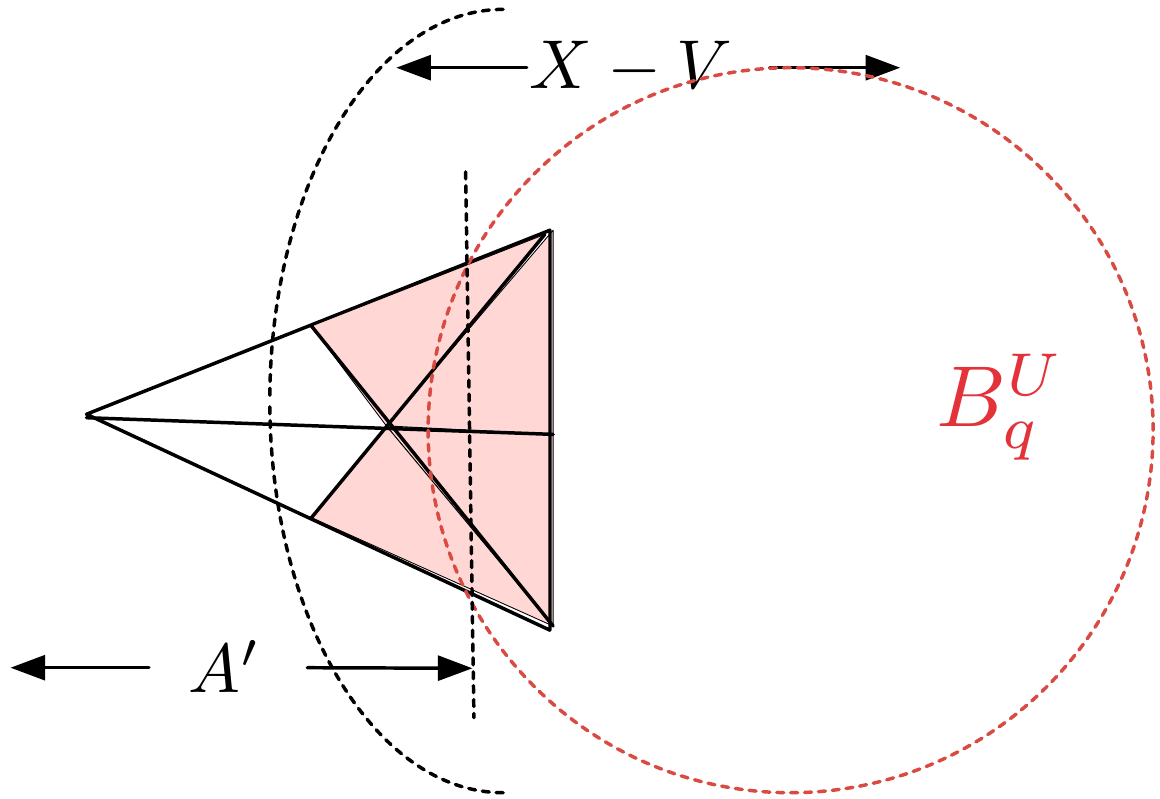}
 \end{center}
\caption[Map $j'$ for a linear singular simplex that requires barycentric subdivision.]
{Map $j'$ for a linear singular simplex that requires barycentric subdivision. In this illustrated example, all four shaded regions are the images of the four singular simplexes in the first barycentric subdivision which are $\A'$-small and have their images in $X - V$. Their formal sum gives a singular chain in $X - V$. Their retraction result in a singular chain in $B_{pq}^{U}$.}
\label{fig:jdetail2}
\end{figure}

We analyze the effect of $j'$ on an arbitrary linear singular simplex $\omega: \Delta_p \to B_p^U$, where
$\omega = l(a_0, ..., a_p)$ for some points $a_i$ in Euclidean space. 
The analysis can be broken up into three main cases:
\begin{itemize}

\item[(A.1)] $\omega(\Delta_p) \subseteq B_q^{U}$: Then $j'$ maps $\omega$ through unchanged, meaning:
\begin{align*}
[\omega: \Delta_p \to B_p^{U}] \xmapsto{j'} [\omega: \Delta_p \to B_{pq}^{U}].
\end{align*}
From now on we simplify notation by omitting the domain and range of the singular simplex, writing instead:
$\omega \xmapsto{j'} \omega.$

\item[(A.2)] $\omega(\Delta_p) \cap B_q^{U} = \emptyset$: Then
$\omega \xmapsto{j'} 0.$

\item[(A.3)] $\omega(\Delta_p) \nsubseteq B_q^{U}$ and $\omega(\Delta_p) \cap B_q^{U} \neq \emptyset$: here we have two sub-cases:

\begin{itemize}

\item[(A.3.a)] $\omega$ is $\A'$-small: This implies that $\omega(\Delta_p) \subseteq X - V$.  
Map $j'$ can be interpreted as a retraction. 
That is, letting $T = \omega(\Delta_p)$, $S = \omega(\Delta_p) \cap B_{q}^{U}$ and $R = \omega(\Delta_p) \cap \bdr B_{q}^{U}$, we define  
$r: T \to S$ by: for $x \in S$, $r(x) = x$; for $x \in T - S$, $r(x) = x'$, where $x' \in R$, as shown in the left of Figure \ref{fig:jdetail}. 
Then
$\omega \xmapsto{j'} \tau,$
where $\tau: \Delta_p \to B_{pq}^U$ is defined by: for $\ep \in \Delta_p$ where $\omega(\ep) \in S$, $\tau(\ep) = \omega(\ep)$;
otherwise for  $\ep \in \Delta_p$ where $\omega(\ep) \notin S$, $\tau(\ep) = r \circ \omega (\ep)$.

\item[(A.3.b)] $\omega$ is not $\A'$-small: We barycentrically subdivide $\omega$ enough times $m$ until $sd^m \omega$ is a $\A'$-small singular chain. 
Then each $\A'$-small singular simplex in $sd^m \omega$ that has its image in $X - V$ follows the pattern of (A.3.a), resulting in a singular simplex $\tau_i: \Delta_p \to B_{pq}^{U}$.
In the end we have,
$\omega \xmapsto{j'} c_{\tau},$
where $c_{\tau}$ is the singular chain, $c_{\tau} = \sum \tau_i$.
This is shown in Figure \ref{fig:jdetail2}.
\end{itemize}

\end{itemize}

\paragraph{Middle map.}

We now describe the action of $f'$ on an arbitrary simplicial linear singular simplex.
Let $\delta:\Delta_p \to |\sd L|$ be such a simplex with $\delta = \eta(\sigma) = l(v_0, ..., v_p)$, for some simplex $\sigma = [v_0, ..., v_p] \in \sd K$.
As above, we have three cases to consider::
\begin{itemize}

\item[(B.1)] $\delta(\Delta_p) \subseteq |\sd K|$: then 
$\delta \xmapsto{f'} \delta.$

\item[(B.2)] $\delta(\Delta_p) \cap |\sd K| = \emptyset$: 
$\delta \xmapsto{f'} 0.$

\item[(B.3)] $\delta(\Delta_p) \nsubseteq |\sd K|$ and $\delta(\Delta_p) \cap |\sd K| \neq \emptyset$:
From Lemma \ref{lemma:boundary} below, we know that $(\delta(\Delta_p) \cap |\sd K|)  \subseteq |\sd K_0|$, and so 
$\delta \xmapsto{f'} 0.$
\end{itemize} 

\begin{lemma}
\label{lemma:boundary}
Given a simplex $\sigma \in L$, if $\sigma \notin  K$ and there exists $\tau < \sigma$ such that $\tau \in  K$, 
then $\tau \in K_0$.
\end{lemma}
\begin{proof}
Suppose there exists $\omega < \tau$ such that $\omega \in K - K_0$. 
This implies that $V^{\omega}$ is completely contained in $B_{pq}^U$. 
Since $V^{\sigma}$ is the intersection of $V^{\omega}$ with the partial Voronoi cells of vertices in $\sigma$ that are not in $\omega$,
then $V^{\sigma}$ should be completely contained in $B_{pq}^U$.
This means that $\sigma$ is in $K$, which leads to a contradiction. \eop
\end{proof}

\subsubsection{Lower Square}

As promised, we now show that the lower square in diagram \ref{dgm:leftbottomchain} commutes.
Choose an arbitrary $\sigma = [v_0, ...,v_p] \in \sd L$, where each $v_i$ is a barycenter of some simplex $\sigma'$ in $L$; as always, we assume
that that the vertices are ordered by increasing dimension of their defining simplices. 
We have two cases:
\begin{itemize}

\item[(C.1)] $\sigma \in \sd K$: by definition, $\eta(\sigma) = l(v_0, ..., v_p)$ has its image in $|\sd K|$, and $f$ is the identify map, that is,
$\sigma \xmapsto{f} \sigma \xmapsto{\eta} \eta(\sigma).$
Meanwhile, by case (B.1),
$\sigma \xmapsto{\eta} \eta(\sigma) \xmapsto{f'} \eta (\sigma).$
Therefore $(\eta \circ f)(\sigma) = f' \circ \eta(\sigma)$.

\item[(C.2)] $\sigma \notin \sd K:$ then
$ \sigma \xmapsto{f} 0 \xmapsto{\eta} 0.$
On the other had, since $\sigma \notin \sd K$, we know that the image of $\delta = \eta(\sigma) = l(v_0, ..., v_p)$ cannot be entirely contained within $|\sd K|$.
There are then two sub-cases to consider:

\begin{itemize}

\item[(C.2.a)] $\delta(\Delta_p) \cap |\sd K| = \emptyset:$ this is case (B.2). We have
$\sigma \xmapsto{\eta} \delta \xmapsto{f'} 0.$
\item[(C.2.b)] $\delta(\Delta_p) \cap |\sd K| \subseteq |\sd K_0|:$ this is case (B.3). We have
$\sigma \xmapsto{\eta} \delta \xmapsto{f'} 0.$
\end{itemize}

\end{itemize}

\subsubsection{Upper Square}

Finally, we show that the upper square in diagram \ref{dgm:leftbottomchain} commutes up to chain homotopy; that is, we will
construct a chain homotopy $D$ between the two chain maps $e= j' \circ h_{\#}$  and $e' = h'_{\#} \circ f'$.
This will of course imply that $e_{*} = e'_{*}$; in other words, that the induced homology diagram commutes. 
For clarity, we zoom in on diagram \ref{dgm:leftbottomchain} and draw the relevant portion below as diagram \ref{dgm:subtopchain}.
\begin{align}
\label{dgm:subtopchain}
 C(B_{p}^{U}, P_0)
&\xrightarrow{j'}
C(B_{pq}^{U}, Z_0)\notag\\
~~~\uparrow{h_{\#}}~~~~~~~~&~~~~~~~~~~~~\uparrow{h'_{\#}}~~~~~~~\notag\\
C(|\sd L|, |\sd L_0|)
 &\xrightarrow{f'}
C(|\sd K|, |\sd K_0|).
\end{align}
For notational ease ,we set $X = |\sd L|$ and $Y = B_{pq}^U$.
To construct $D$, we will define for each $p$ a chain map $F_p: C_p(X \times I) \to C_p(Y),$ and then we will set
$D_p = F_{p+1} \circ G_p,$ where $G_p: C_p(X \times I) \to C_{p+1}(X \times I)$ is given by Lemma \ref{lemma:prism} below.

\paragraph{Construction of F.}

Let $\pi : X \times I \to X$ be projection on the first factor, and fix an arbitrary 
simplicial linear singular simplex $\kappa: \Delta_p \to X \times I$.
Then $\pi_{\#}(\kappa) = \delta = l(\hat{\sigma}_0, ..., \hat{\sigma}_p),$ for
some simplex $\sigma = [\hat{\sigma}_0, ..., \hat{\sigma}_p]$ in $\sd L$.
We define $F$ in stages, based on properties of $\delta$, as follows.

\begin{itemize}
\item[(D.1)] $\delta(\Delta_p) \subseteq |\sd K|$: 
following the $e'$-path and case (B.1), we have
$\delta \xmapsto{f'} \delta \xmapsto{h'_{\#}} \tau'.$
On the other hand, following the $e$-path results in
$\delta \xmapsto{h_{\#}} \omega.$
We now have three sub-cases, based on varying properties of $\omega$:
\begin{itemize}
\item[(D.1.a)] $\omega(\Delta_p) \subseteq B_q^{U}:$ following the $e$-path and case (A.1).
we have,
$\delta \xmapsto{h_{\#}} \omega \xmapsto{j'} \tau, $
where $\tau = \omega$ except for differing range.
We then can define $F(\kappa) = \iota$, where $\iota: \Delta_p \to Y$ is given by:
for every $\epsilon \in \Delta_p$, where $\kappa(\ep) = (x, t) \in X \times I$, 
$\iota(\ep) = (1 - t) \tau(\ep) + t \tau'(\ep)$.
This formula is illustrated in Figure \ref{fig:showmapF}.
\begin{figure}[tbp]
 \begin{center}
  \includegraphics[scale=0.65]{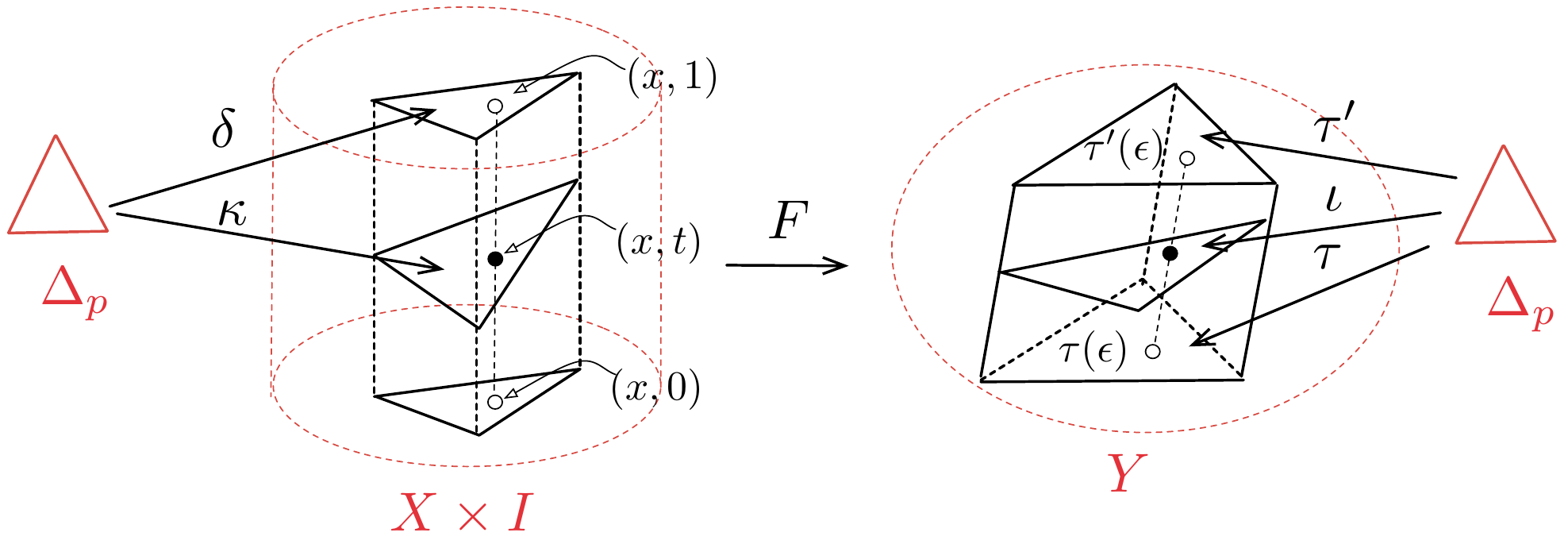}
 \end{center}
\caption[Case (D.1.a): illustration of $F$.]
{Case (D.1.a): illustration of $F$.}
\label{fig:showmapF}
\end{figure}

\item[(D.1.b)] $\omega(\Delta_p) \cap B_q^{U} = \emptyset$: This is case (A.2).
We branch further as follows:
\begin{itemize}

\item[(i)] $\delta(\Delta_p) \subseteq |\sd K_0|$: 
Following the $e'$-path,
$\delta \xmapsto{f'} 0 \xmapsto{h'_{\#}} 0.$
Similarly, following the $e$-path,
$\delta \xmapsto{h_{\#}} \omega \xmapsto{j'} 0.$
We define $F(\kappa) = 0$.

\item[(ii)] $\delta(\Delta_p) \nsubseteq |\sd K_0|$:
this is not possible.
Suppose it were.
This implies that there exists at least one vertex $\hat{\sigma}_i$ of $\sigma$ such that
$V^{\sigma_i} \cap B_{pq}^{U} \neq \emptyset$ and $V^{\sigma_i} \cap Z_0 = \emptyset$.
This means that $V^{\sigma_i}$ is completely contained in $B_r(q)$.
Therefore $h(\hat{\sigma}_i)$ is contained in $B_r(q)$, which contradicts our assumption.
\end{itemize}

\item[(D.1.c)] $\omega(\Delta_p) \nsubseteq B_q^{U}$ and $\omega(\Delta_p) \cap B_q^{U} \neq \emptyset$: 
we must consider two further sub-cases.
\begin{itemize}
\item[(i)] $\omega$ is $\A'$-small: this is case (A.3.a),
and we define $F(\kappa)$ similarly to case (D.1.a). 

\item[(ii)] $\omega$ is not $\A'$-small: this is case (A.3.b).
Then $\delta \xmapsto{h_{\#}} \omega \xmapsto{j'} c_\tau,$
where $c_{\tau} = \sum \tau_i$ for some collection of $\tau_i: \Delta_p \to B_{pq}^{U}$.
We now define $F(\kappa) = c_{\iota}$, where $c_{\iota} = \sum \iota_i$ and each singular simplex $\iota_i: \Delta_p \to Y$.
is defined as follows.
Let $m$ be the smallest integer such that $sd^m \omega$ is $\A'$-small.
For each singular simplex $\tau_i$ in $c_{\tau}$, there exists a singular simplex $\omega_i$ in $sd^m \omega$ such that $j'(\omega_i) = \tau_i$.
For each such $\omega_i$, there exists a singular simplex $\delta_i$ in $sd^m \delta$ such that $h_{\#}(\delta_i) = \omega_i$.
In other words, for each $\tau_i$ in $c_{\tau}$, there exist $\delta_i$ in $sd^m \delta$, such that following the $e$-path,
$\delta_i \xmapsto{h_{\#}} \omega_i \xmapsto{j'} \tau_i.$
Meanwhile, for each such $\delta_i$, following the $e'$-path gives
$\delta_i \xmapsto{f'} \delta_i \xmapsto{h'_{\#}} \tau'_i.$

On the other hand, for each such $\delta_i$, there exists a corresponding $\kappa_i$ in $sd^m \kappa$, such that $\delta_i = \pi(\kappa_i)$. 
We now define $\iota_i$ for each such $\kappa_i$.
For all $\ep \in \Delta_p$ where $\kappa_i(\ep) = (x, t) \in X \times I$,
$\iota_i(\ep) = (1-t) \tau_i(\ep) + t \tau'_i(\ep).$
This is illustrated in Figure \ref{fig:mapF2}.

\begin{figure}[tbp]
 \begin{center}
  \includegraphics[scale=0.5]{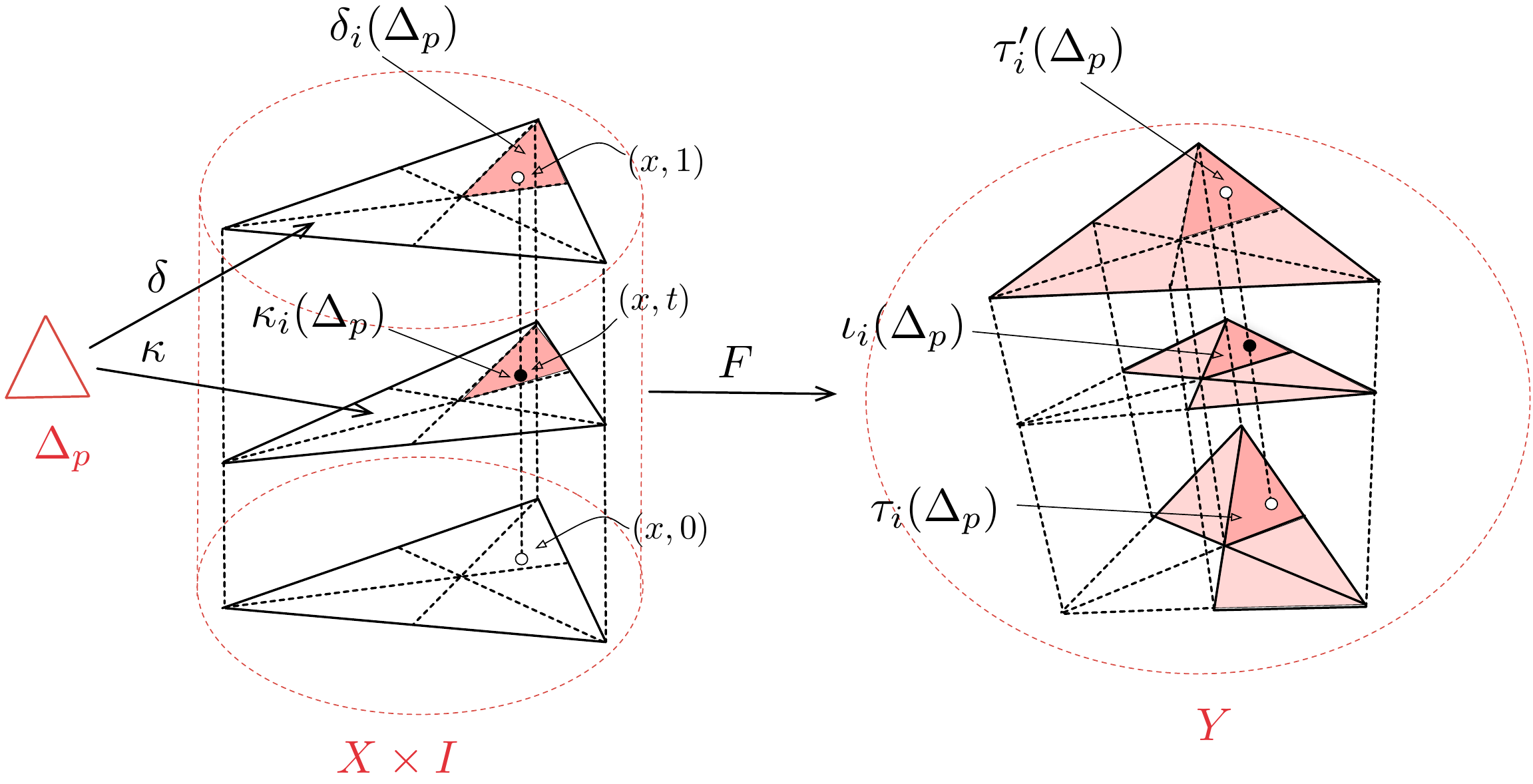}
 \end{center}
\caption[Case (D.1.c): illustration of $F$.]
{Case (D.1.c): illustration of $F$. Left: the dark shaded region is the minimal carrier of $\delta_i$ and $\kappa_i$. Right: the shaded regions from top to bottom are the minimal carriers of $\tau'$, $c_{\iota}$ and $c_{\tau}$ respectively; the dark shaded regions from top to bottom are the minimal carriers of $\tau'_i$, $\iota_i$ and $\tau_i$, respectively.
For simplicity, we illustrate the minimal carrier of the singular chain $c_{\tau}$ as the union of the minimal carriers of its simplexes before their retraction.}
\label{fig:mapF2}
\end{figure}

\end{itemize}

\end{itemize}

\item[(D.2)] $\delta(\Delta_p) \nsubseteq |\sd K|:$ we again have two sub-cases:

\begin{itemize}
\item[(D.2.a)] $\delta(\Delta_p) \cap |\sd K| = \emptyset:$ following the $e'$-path and case (B.2),
$\delta \xmapsto{f'} 0 \xmapsto{h'_{\#}} 0. $
Since $\delta(\Delta_p) \cap |\sd K| = \emptyset$, this implies that its corresponding $\sigma \notin \sd K$.
That is, for all $\hat{\sigma}_i$,
$V^{\sigma_i} \cap B_{pq}^{U} = \emptyset$, therefore
 all $h(\hat{\sigma}_i)$ lie outside of $B_r(q)$.
 Let $\omega = h_{\#}(\delta) = h \circ \delta$.
 This means $\omega$ has its image outside of $B_q^{U}$.
Then following the $e$-path, 
$\delta \xmapsto{h_{\#}} \omega \xmapsto{j'} 0.$
We define $F(\kappa) = 0$.

\item[(D.2.b)] $(\delta(\Delta_p) \cap |\sd K|) \subseteq |\sd K_0|:$
following the $e'$-path and case (B.3),
we have,
$\sigma \xmapsto{f'} 0 \xmapsto{h'_{\#}} 0.$
On the other hand, let $\omega = h_{\#}(\delta) = h \circ \delta$. Then $\omega(\Delta_p) \subseteq P_0$ and so 
following the $e$-path give
$\sigma \xmapsto{h_{\#}} \omega \xmapsto{j'} 0.$
We define $F(\kappa) = 0$.

\end{itemize}

\end{itemize}

\paragraph{Construction of D.}

\begin{figure}[tbp]
 \begin{center}
  \includegraphics[scale=0.3]{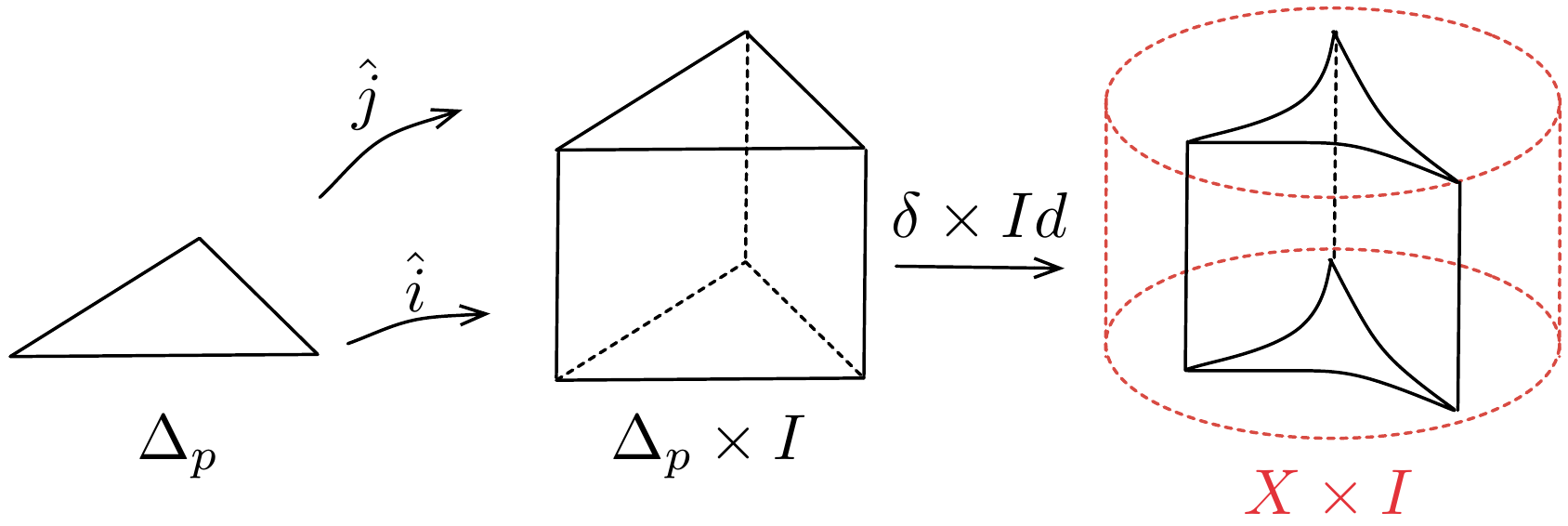}
 \end{center}
\caption[Illustration of $G$.]
{Illustration of $G$.}
\label{fig:mapDx}
\end{figure}

To define our chain homotopy $D$, we first need the following lemma:
\begin{lemma}(\cite{Mun1984}, page 171)
\label{lemma:prism}
There exists for each space $X$ and each non-negative integer $p$, a homomorphism
$G_p: C_p(X) \to C_{p+1}(X \times I),$
having the following property:
if $\delta: \Delta_p \to X$ is a singular simplex, 
then 
$\bdr G \delta + G \bdr \delta = j_{\#}(\delta) + i_{\#}(\delta),$
where the map $i: X \to X \times I$ carries $x$ to $(x,0)$, and the map $j: X \to X \times I$ carries $x$ to $(x,1)$.
\end{lemma}
This homomorphism is illustrated intuitively in Figure \ref{fig:mapDx}, where $\delta \times Id$ carries a singular $p+1$ chain that fills up the entire
prism $\Delta_p \times I$ to a singular chain on $X \times I$, and the maps $\hat{i}, \hat{j}: \Delta_p \to \Delta_p \times I$ carry each $x$ to $(x,0)$
and $(x,1)$ respectively.
Then, as promised, we set $D_p = F_{p+1} \circ G_p$.
To show that $D$ is a chain homotopy between $e$ and $e'$, we calculate
\begin{align*}
\bdr D & = \bdr (F G) \\
& = F \bdr G \\
& = F(j_{\#} + i_{\#} + G \bdr)\\
& = F j_{\#} + F i_{\#} + F G \bdr\\
& = F j_{\#} + F i_{\#} +  D \bdr
\end{align*}
Hence we need only show that $F j_{\#} = e'$ and $F i_{\#} = e$ to complete the argument.
In the case when $F(\kappa)$ is defined to be $0$, the corresponding $e(\delta)$ and $e'(\delta)$ are also $0$, so this is no problem.
In the case when $F(\kappa)$ is not defined to be $0$, as shown in Figure \ref{fig:showmapF} and Figure \ref{fig:mapF2}, $F j_{\#} (\delta) = e'(\delta)$, and $F i_{\#} (\delta) = e(\delta)$.
This concludes the proof that the upper square in diagram \ref{dgm:leftbottomchain} commutes up to chain homotopy,
and thus that the bottom square of diagram \ref{dgm:left} commutes.

\subsubsection{Top Square of Diagram \ref{dgm:left}}

As promised above, we now prove that the top square of diagram \ref{dgm:left} commutes, which will complete the proof that the left
face of Figure \ref{fig:strata-doublecube} commutes.
In fact, the top square commutes on the chain level, which we draw directly below.
\begin{align*}
C (B_{p}^{U}, \bdr B_{p}^{U})
  &\xrightarrow{j}
 C (B_{pq}^{U}, \bdr B_{pq}^{U})\notag\\
 \downarrow {i_{\#}}~~~~~~&~~~~~~~\downarrow{i'_{\#}}~~~~~~~\notag\\
C (B_{p}^{U}, P_0)
&\xrightarrow{j'}
C (B_{pq}^{U}, Z_0)\notag\\
\label{dgm:lastchain}
\end{align*}
Setting $e = j' \circ i_{\#}$ and $e' = i'_{\#} \circ j$, we show, once
again via an exhaustive case analysis, that $e = e'$.

First we need to understand the map $j$ for a linear singular simplex.
The interpretation of $j$ is almost the same as that of $j'$ (case (A)). 
More specifically, we let $\omega: \Delta_p \to B_p^U$ be an arbitrary linear singular simplex.
There are three cases:
\begin{itemize}

\item[(E.1)] $\omega(\Delta_p) \subseteq B_q^{U}:$ then
$\omega \xmapsto{j} \omega.$

\item[(E.2)] $\omega(\Delta_p) \cap B_q^{U} = \emptyset:$ then
$\omega \xmapsto{j} 0.$

\item[(E.3)] $\omega(\Delta_p) \nsubseteq B_q^{U}$ and $\omega(\Delta_p) \cap B_q^{U} \neq \emptyset:$
We have two sub-cases:
\begin{itemize}

\item[(E.3.a)] $\omega$ is $\A$-small: then
$\omega \xmapsto{j} \gamma,$
where $\gamma: \Delta_p \to B_{pq}^U$ is defined via the retraction-type arguments above.

\item[(E.3.b)] $\omega$ is not $\A$-small: then
$\omega \xmapsto{j'} c_{\gamma},$
where $c_{\gamma} = \sum \gamma_i$, with each $\gamma_i: \Delta_p \to B_{pq}^{U}$ described
by the subdivision and retraction arguments we have already given. 

\end{itemize}

\end{itemize}

To complete the proof, we fix an arbitrary singular simplex $\delta: \Delta_p \to B_p^U$, and again
argue by cases.
\begin{itemize}
\item[(F.1)]   $\delta(\Delta_p) \subseteq B_q^{U}:$ exploiting the analysis above, we note that following the $e'$-path
results in $\delta \xmapsto{j} \delta \xmapsto{i'_{\#}} \delta,$
while following the $e$-path gives
$\delta \xmapsto{i_{\#}} \delta \xmapsto{j'} \delta,$
as needed.

\item[(F.2)] $\delta(\Delta_p) \cap B_q^{U} = \emptyset:$ here both paths result in $0$.

\item[(F.3)] $\delta(\Delta_p) \nsubseteq B_q^{U}$ and $\delta(\Delta_p) \cap B_q^{U} \neq \emptyset:$
here we must analyze two sub-cases:

\begin{itemize}

\item[(F.3.a)] $\delta$ is $\A$-small: this implies that $\delta(\Delta_p) \subseteq X - V$. Following the $e$-path gives.
$\delta \xmapsto{j} \gamma \xmapsto{i'_{\#}} \gamma $.
On the other hand, $\delta$ is also $\A'$-small, since $\A$ and $\A'$ share the element $X - V$, and hence the $e'$ path
\begin{align*}
\delta \xmapsto{j} \delta \xmapsto{i'_{\#}} \tau.
\end{align*}
But really the fact that $X - V$ is part of $\A$ and $\A'$ means that $\tau$ and $\gamma$ follow the same retraction, and thus $\gamma = \tau$. 

\item[(F.3.b)] $\delta$ is not $\A$-small: the analysis here is the same as the last case, with some words about subdivision added. 

\end{itemize}
\end{itemize}

\subsection{Finale}
We are now ready to finish the proof of Theorem \ref{equaldgm},
which boils down to verifying that diagram \ref{diag:KEE} commutes, with the vertical maps
being isomorphisms. That is,
\begin{align*}
\ldots \to  &\kernel \phi_{\alpha}^U   \to \kernel \phi_{\beta}^U  \to \ldots \notag\\
   &~\uparrow  \cong  ~~~~~~~~~~\uparrow \cong \notag\\
   \ldots \to & \kernel \psi_{\alpha}  \to  \kernel \psi_{\beta} \to \ldots.
\end{align*}
But this is now just easy diagram-chasing.
Commutativity of diagram \ref{diag:KEE} follows directly from the commutativity of the bottom face of the double-cube in Figure \ref{fig:strata-doublecube}, 
and the leftmost (rightmost) vertical isomorphism derives from our statements about the left (right) face of the double-cube.
The commutativity of the top face implies that the cokernel analogue to diagram \ref{diag:KEE} commutes, after a little more algebra which we omit.

\section{Proof of Theorem \ref{result:GIIT}}
\label{app:geoproof}

In this Appendix, we give a proof for Theorem \ref{result:GIIT}.
First we need a technical lemma involving some simple algebraic topology.

\subsection{Absolute Homology Modules}

Recall from before that $\sigma(p,r)$ is the feature size of the relative homology persistence module
$\{\Hgroup(B_{p}^{\Xspace}, \bdr B_{p}^{\Xspace})\}$.
On the other hand, the same thickening process also defines two absolute homology persistence modules,
$\{\Hgroup(B_{p}^{\Xspace})\}$ and $\{\Hgroup(\bdr B_{p}^{\Xspace})\}$.
We let $\sigma_i(p,r)$ and $\sigma_b(p,r)$ denote the feature sizes of these modules.
Similarly, we define $\sigma_i(p,q,r)$ and $\sigma_b(p,q,r)$, respectively, to be the feature sizes of the absolute homology persistence modules
$\{\Hgroup(B_{pq}^{\Xspace})\}$ and $\{\Hgroup(\bdr B_{pq}^{\Xspace})\}$.
\begin{theorem}[Relative/Absolute Lemma]
The feature size of each relative module is at least the minimum of those of its two associated absolute modules:
\begin{eqnarray*}
 \sigma(p,r) &\geq& \min\{\sigma_i(p,r), \sigma_b(p,r)\},\\
 \sigma(p,q,r)& \geq& \min\{\sigma_i(p,q,r), \sigma_b(p,q,r)\}.
 \end{eqnarray*}
\label{result:RAL}
\end{theorem}
\begin{proof}
We prove the first equality; the second can then be proven with only minor notational adjustment.
For any two non-negative reals $\alpha < \beta$, and for each homological dimension $i \geq 0$, consider the following commutative diagram:
\begin{align}
\Hgr_i(\bdr B_{p}^{\Xspace}(\alpha))  \rightarrow  \Hgr_i(B_{p}^{\Xspace}(\alpha))  \rightarrow  \Hgr_i(B_{p}^{\Xspace}(\alpha), \bdr B_{p}^{\Xspace}(\alpha))  \rightarrow 
 \Hgr_{i-1}(\bdr B_{p}^{\Xspace}(\alpha))  \rightarrow  \Hgr_{i-1}(B_{p}^{\Xspace}(\alpha)) \notag \\
 \downarrow ~~~~~~~~~~~~~~~~~~~~~ \downarrow ~~~~~~~~~~~~~~~~~~~~~~~~~~~~ \downarrow ~~~~~~~~~~~~~~~~~~~~~~~~ \downarrow ~~~~~~~~~~~~~~~~~~~~~ \downarrow ~~~~~~~~~~~~~~~\notag\\  
 \Hgr_i(\bdr B_{p}^{\Xspace}(\beta))  \rightarrow  \Hgr_i(B_{p}^{\Xspace}(\beta))  \rightarrow  \Hgr_i(B_{p}^{\Xspace}(\beta), \bdr B_{p}^{\Xspace}(\beta))  \rightarrow 
 \Hgr_{i-1}(\bdr B_{p}^{\Xspace}(\beta))  \rightarrow  \Hgr_{i-1}(B_{p}^{\Xspace}(\beta))
\end{align}
where the vertical maps are induced by the inclusion $X_{\alpha} \hookrightarrow X_{\beta}$ and the two rows
come from the long exact sequences of the pairs $(B_{p}^{\Xspace}(\alpha), \bdr B_{p}^{\Xspace}(\alpha))$ and $(B_{p}^{\Xspace}(\beta), \bdr B_{p}^{\Xspace}(\beta))$ (\cite{Mun1984}).

Suppose that the middle vertical map fails to be an isomorphism.
Then the Five-Lemma (\cite{Mun1984}, p.140) tells us that at least one of the other four vertical maps will also fail
to be an isomorphism.
In other words, any change within the relative module must be accompanied by a simultaneous change in at least one of the two absolute modules.
The inequality follows. \eop
\end{proof}

This theorem together with the definition of $\rho(p,q,r)$ implies the following inequality
\begin{equation}
\rho(p,q,r) \geq \min\{\sigma_i(p,r),\sigma_b(p,r),\sigma_i(p,q,r),\sigma_b(p,q,r)\}.
\label{eqn:relabs}
\end{equation}

\subsection{Proof}

To prove Theorem \ref{result:GIIT}, we will further lower bound the $\sigma$-parameters above.
This is accomplished via one more lemma.
\begin{lemma}[Deformation Lemmas]
The following four claims all hold for every small enough $\delta > 0$. In each of
the claims, the homotopy equivalence is given by a deformation retraction:
\begin{align*}
& \forall \alpha < \min\{\tau(p,r), \eta(p,r)\},
(\Xspace_{\alpha} \cap B_r(p)) \simeq (\Xspace_{\delta} \cap B_r(p)),\\
& \forall \alpha < \min\{\tau_0(p,q,r), \eta(p,r)\},
(\Xspace_{\alpha} \cap \partial B_r(p)) \simeq (\Xspace_{\delta} \cap \partial B_r(p)),\\
& \forall \alpha < \min\{\tau(p,q,r), \eta(p,q,r)\},
(\Xspace_{\alpha} \cap B_r(p) \cap B_r(q)) \simeq (\Xspace_{\delta} \cap B_r(p) \cap B_r(q)),\\
& \forall \alpha < \min\{\tau_0(p,q,r), \eta(p,q,r)\},
(\Xspace_{\alpha} \cap \partial(B_r(p) \cap B_r(q))) \simeq (\Xspace_{\delta} \cap \partial(B_r(p) \cap B_r(q))).
\end{align*}
\label{result:DL}
\end{lemma}
\begin{proof}
All four claims follow from Stratified Morse Theory \cite{GorMac1988}.
We prove only the first claim; the other three can be proven with only slight modifications.
Consider the stratification of $B_r(p)$ with singular set $\Sigma = \Mcal(p,r) \cup \partial B_r(p)$ and
whatever further decomposition of $\Sigma$ is needed.
Setting $d = d_{\Xspace}|B_r(p) : B_r(p) \to \Rspace$, we note that the sets $X_{\alpha} \cap B_r(p)$ are simply
the sublevel sets of $d$ for various parameters $\alpha$.
Generically, $d$ will be a Stratified Morse function on $B_r(p)$ with its above stratification.
Consider the set $H$ of all critical points of $d$ which have positive $d$-value.

We claim $H \subset (\Mcal(p,r) \cup G(p,r)):$ to see this, we suppose $y \in H$ and we assume first
that $y$ is in the interior of $B_r(p)$. Then $y$ is also a critical point of the globally defined function $d_{\Xspace}$,
and since $d(x) = d_{\Xspace}(x) > 0$, we know that $y \in \Mcal$. Since $y$ is also in $B_r(p)$ by assumption,
we know in fact that $y \in \Mcal(p,r)$.
On the other hand, suppose that $y \in \partial B_r(p)$; we can also assume that $y \not \in \Mcal(p,r)$ or we are already done.
Then by definition $y$ is a critical point of the restriction of $d_{\Xspace}$ to $\partial B_r(p)$.
Since the gradient of this latter function is simply the projection of $\nabla d_{\Xspace}$ onto $\partial B_r(p)$,
we can conclude $y \in G(p,r)$.

In other words, if $\alpha < \{\tau(p,r), \eta(p,r)\}$, then $(\Xspace_{\alpha} \cap B_r(p)) \cap H = \emptyset$, and
hence the interval $[\delta, \alpha]$ contains no critical values of $d$.
The claim then follows from the first fundamental theorem of Stratified Morse Theory \cite{GorMac1988}. \eop
\end{proof}

Finally, we finish the proof of Theorem \ref{result:GIIT}, which we restate here for convenience.

\begin{theorem}[Geometric lower bound]
If we define
$$\gamma = \gamma(p,q,r) = \min\{\tau(p,r), \tau(p,q,r), \eta(p,r), \eta(p,q,r)\},$$
then
$\rho(p,q,r) \geq \gamma(p,q,r)$.
\end{theorem}

\begin{proof}

Note that $\tau(p,r) \leq \tau_0(p,r)$ and $\tau(p,q,r) \leq \tau_0(p,q,r)$ so we need not consider
$ \tau_0(p,r)$ and  $\tau_0(p,q,r)$.

Recall $\sigma_i(p,r)$ and $\sigma_b(p,r)$ were defined to be the feature sizes of the persistence modules
$\{\Hgroup(B_{p}^{\Xspace}(\alpha))\}$ and $\{\Hgroup(\bdr B_{p}^{\Xspace}(\alpha))\}$,
respectively.

By the first and second of the Deformation Lemmas the following holds
$$\sigma_i(p,r), \sigma_b(p,r) \geq \min\{\tau(p,r), \eta(p,r)\}.$$
For the same reason
$$\sigma_i(p,q,r), \sigma_b(p,q,r) \geq \min\{\tau(p,q,r),\eta(p,q,r)\}.$$
These inequalities, together with (\ref{eqn:relabs})
$$\rho(p,q,r) \geq \min\{\sigma_i(p,r),\sigma_b(p,r),\sigma_i(p,q,r),\sigma_b(p,q,r)\},$$
prove the theorem, $\rho(p,q,r) \geq \gamma(p,q,r)$. \eop

\end{proof}

\section*{Acknowledgments}
All the authors would like to thank Herbert Edelsbrunner and John Harer for useful discussions and suggestions. PB would like to thank David Cohen-Steiner and Dmitrity Morozov
for helpful discussion, and SM would like to
thank Shmuel Weinberger for useful comments. SM and BW would like to acknowledge the support of NIH
Grant P50 GM 081883, and SM would like to acknowledge the support of NSF Grant DMS-07-32260
and Grant R01 CA123175-01A1. PB thanks the Computer Science Department at Duke University for 
hosting him during the Spring semester of 2010.

\bibliographystyle{plain}
\bibliography{BMW10-refs}

\end{document}